\documentclass[12pt]{article}

\usepackage{amsmath,amsthm,amssymb}

\theoremstyle{plain}
\newtheorem{Thm}{Theorem}
\newtheorem{Prop}{Proposition}[section]
\newtheorem{Cor}{Corollary}
\newtheorem{lem}{Lemma}[section]
\theoremstyle{remark}
\newtheorem{rem}{\indent \sc Remark}[section]

\makeatletter

\@addtoreset{equation}{section}

\begin{document}

\def\R{\mathbb{R}}
\def\Z{\mathbb{Z}}
\def\N{\mathbb{N}}
\def\H{\mathbb{H}}
\def\C{\mathbb{C}}
\def\D{\mathbb{D}}
\def\U{\mathcal{U}}
\def\Exp{{\bf E}}
\def\Prob{{\bf P}}
\def\V{{\bf V}}
\def\v2{\vskip2mm}
\def\n{\noindent}
\renewcommand{\Re}{{\rm Re}}
\renewcommand{\Im}{{\rm Im}}
\def\sgn{{\rm sgn}}
\def\hcap{{\rm hcap}}
\def\rad{{\rm rad}}
\def\diam{{\rm diam}}
\def\dist{{\rm dist}}
\def\textmc{\rm}
\def\({(\!(}
\def\){)\!)}
\def\F{{\cal F}}
\def\B{{\cal B}}
\def\a{\alpha}
\def\b{\beta}
\def\e{\varepsilon}
\def\de{\delta}
\def\ga{\gamma}
\def\k{\kappa}
\def\la{\lambda}
\def\fa{\varphi}
\def\th{\theta}
\def\si{\sigma}
\def\t{\tau}
\def\om{\omega}
\def\De{\Delta}
\def\Ga{\Gamma}
\def\La{\Lambda}
\def\Om{\Omega}
\def\Th{\Theta}
\def\lan{\langle}
\def\ran{\rangle}
\def\lbr{\left(}
\def\rbr{\right)}
\def\pf{{\it Proof.}}
\def\v2{\vskip2mm}
\def\n{\noindent}
\def\z{{\bf z}}
\def\x{{\bf x}}
\def\y{{\bf y}}
\def\bv{{\bf v}}
\def\be{{\bf e}}
\def\0{{\bf 0}}
\def\pr{{\rm pr}}
\def\cp{{\rm Cap}}
\def\tst12{{\textstyle \frac12}}
\def\1{{\bf 1}}
\def\sg{{\rm sgn\,}}

\def\n{\noindent}
\def\beq{\begin{eqnarray*}}
\def\eeq{\end{eqnarray*}}
\def\supp{\mbox{supp}}
\def\beqn{\begin{equation}}
\def\eeqn{\end{equation}}

\begin{center}
{\bf  Recurrent random walks on $\mathbb{Z}$ with infinite variance: \\
 transition probabilities of them killed on a finite set} \\
\vskip4mm
{K\^ohei UCHIYAMA} \\
\vskip2mm
Department of Mathematics, Tokyo Institute of Technology \\
Oh-okayama, Meguro Tokyo 152-8551\\
e-mail: \,uchiyama@math.titech.ac.jp
\end{center}

\vskip6mm

{\it running head}:   Recurrent random walks on $\Z$ 

{\it key words}:  one dimensional random walk; 
first passage time; killing at the origin; in a domain of attraction, transition probability,
tunneling


\vskip6mm

\begin{abstract}
In this paper 
we consider an irreducible  random walk on the integer lattice $\Z$ that is in the  domain 
of normal attraction of a strictly stable process with index $\alpha\in (1, 2)$ and  obtain the asymptotic form of the distribution of the hitting time of the origin and that of the transition probability for the  walk killed when it hits a finite set. The asymptotic forms obtained 
are valid uniformly in the natural  domain of the space and time variables. 
 \end{abstract}
\vskip6mm

\section{Introduction}

Let   $S_n=X_1+\cdots+ X_n$ be a  random walk on the integer lattice  $\Z$ started at   $S_0 \equiv 0$, where  the increments $ X_1,  X_2, \ldots$ are  independent and identically distributed  random variables  defined on some probability  space $(\Om, \F, P)$ and taking values in $\Z$. Let $E$  indicate the expectation under $P$  as usual and $ X$ be a random variable having the same law as $ X_1$.  
We suppose throughout the paper that the walk $S_n$ is 
\v2
1)  in the  domain of normal attraction of a strictly stable law of index $1<\a<2$ or, what amounts to the same thing (cf \cite{F}),  if $\phi(\th):= Ee^{i\th X}$, then 
\beqn\label{f_hyp}\lim_{\th \to \pm 0} \frac{1-\phi(\th)}{|\th|^\a} = c_\circ e^{\pm i\pi\ga/2}
\eeqn
with some real numbers  $c_\circ$  and $\ga$  such that $c_\circ >0$ and $|\ga|\leq 2-\a$. 
\v2\n
For simplicity we also suppose (except in Theorem \ref{thm7}) that 
\v2
2) the walk is  {\it strongly aperiodic} in the sense of  Spitzer \cite{S}, namely for any $x\in \Z$, $P[S_n =x] >0$ for all sufficiently large $n$. 

\v2
The condition  1)  entails $EX =0$ so that  the walk is recurrent. (See Section 7.2 for an equivalent  condition in terms of the tails of distribution function of $X$ and some related facts.) The condition 2)  gives rise to no loss of generality  (see Remark \ref{R8}(b)). 

Under these assumptions  
we  obtain in this paper  precise  asymptotic forms of the distribution of the hitting time of the origin and of  the transition probability for the  walk killed when it hits the origin.  
The estimates obtained are  uniform for the space variables within the natural space-time regime $x=O(n^{1/\a})$.  We extend the results  to the case when  the walk is killed on hitting a finite set instead of the origin. 
The corresponding  results are obtained  for the walks with finite variance by the present author \cite{U1dm}, \cite{U1dm_f}. In a classical paper   \cite{K} Kesten studied  similar problems  and obtained an exact asymptotic result for the ratio of transition probability and hitting time \lq density'   under
 a mild  assumption on the walk  where, however, the space variables are fixed (cf.  Remark  \ref{rem4C}  at the end of the next section).     Although we consider the problem for all admissible  $\ga$, our main  interest is in the extreme case $\ga =|2-\alpha|$ when the  limiting stable
process has jumps only in one direction: the other case is much simpler and  the asymptotic forms  obtained are   quite different between  the  two cases for large space variables.  The condition 1) is  restrictive  and it is desirable  to replace it by a weaker one, to that end  however we encounter a serious difficulty  for the present approach.  In any case it must be worth to reveal 
what kind of behaviour  of the transition probability of the killed process even under such a restrictive condition.

\section{Statements of results}

 We  first introduce fundamental objects  that appear in the description of our results and state some well known facts concerning them. 
Put $p^n(x)=P[S_n=x]$, $p(x)=p^1(x)$ ($x\in \Z$) and define the potential function 
 $$a(x)=\sum_{n=0}^\infty[p^n(0)-p^n(-x)];$$
 the series on the RHS is convergent  and   $a(x)/|x|\to 1/\sigma^2$ and $a(x+y)-a(x) \to \pm y/\sigma^2$ as $x\to \pm \infty$ (cf. Spitzer \cite{S}:Sections 28 and 29).  
To make  expressions concise we  use the notation
$$a^\dagger(x) = \1(x=0) + a(x),$$
where $\1({\cal S})$ equals  1 or 0 according as  a statement ${\cal S}$ is true or false.
The condition 3)  in Introduction entails that $a(x)>0$ whenever  $x\neq 0$, whereas 
if $S$ is left-continuous (i.e., $p(x)=0$ for $x\leq2$), then $a(x)= 0$ for all $x\geq 0$ (under $\sigma=\infty$), and similarly for  right-continuous walks.  (See Section 8.3 for additional facts related to $a$.)

 We write  $S^x_n$ for  $x+S_n$, the walk started at $x\in \Z$. For a subset $B\subset \R$,   put  
 $$\sigma^x_{B}=\inf \{n\geq 1: S^x_n\in B\},$$
  the time of the first entrance of  the walk $S^x$ into $B$. 
To avoid the overburdening of notation we write $S^x_{\sigma_B}$    for  $S^x_{\sigma^x_B}$ and 
 $S_{\sigma_B}$ for $S_{\sigma^0_B}$; sometimes $\sigma B$ is written for $\sigma_B$, e.g., $S^x_{\sigma [0,\infty)}$ for $S^x_{\sigma_{[0,\infty)}}$.

 When the  spatial variables  become indefinitely large the asymptotic results are naturally  expressed by means of the   stable process appearing in the scaling limit and we need to introduce relevant quantities.   
Let $Y_t$ be a stable process started at zero with characteristic exponent
$$\psi(\th) =  e^{ i (\sg \th) \pi\ga/2}|\th|^a  \qquad (|\ga| \leq 2-\a, \ga \;\mbox{is real})$$
so that $Ee^{i\th Y_t} = e^{- t\psi(\th)}$, where  $\sg \th =1$ if $\th >0$, $0$ if $\th=0$ and  $-1$ if $\th <0$.  ($\ga$ has the same sign as the skewness parameter so that  the extremal case $\ga=2-\a$ corresponds to the spectrally  positive case.)  Denote by 
 $\mathfrak{p}_t(x)$ and  $\mathfrak{f}^x(t)$  the  density of the distribution of $Y_t$ and of the  first hitting time to the origin  by $Y^x_t := x+Y_t$, respectively:
 $$ \mathfrak{p}_t(x)= P[Y_t\in dx]/dx,  \qquad  \mathfrak{f}^x(t) = (d/dt)P[ \exists s\leq t,   Y^x_s =0];$$
 there exist the jointly continuous versions of these densities  (for $t>0$) and we shall always choose such ones.  It follows that $S_{\lfloor nt\rfloor}/n^{1/\a} \Rightarrow Y_{c_\circ t}$ (weak convergence of distribution) and  by Gnedenko's local limit theorem \cite{GK}
 \beqn\label{eqLLT}
 \lim_{n\to\infty} \sup_{x\in \Z} | n^{1/\a}p^n(x) - \mathfrak{p}_{c_\circ} (x/n^{1/\a}) | = 0,
 \eeqn 
where $\lfloor b\rfloor$ denotes the integer part of a real number $b$.
 For real numbers $s, t$, 
 $s\vee t= \max\{s,t\}$ and  $s\wedge t=\max\{s,t\}$, $t_+=t\vee 0$, $t_- = (-t)_+$  and  $\lceil t\rceil$ denotes the smallest integer that does not less than  $t$;  for positive sequences $(s_n)$ and $(t_n)$, $s_n\sim t_n$ and $s_n\asymp t_n$  mean, respectively,  that the ratio  $s_n/t_n$ approaches unity and    that  $s_n/t_n$ is bounded away from zero and infinity.
 We use the letters $x,y,z$ and  $w$ to represent integers  which indicate points assumed by the walk when discussing  matters on the random walk, while the same letters may stand for real numbers when  the stable process is dealt with;   we  shall sometimes use the Greek  letters $\xi$, $\eta$ etc. to denote the real variables  the stable process may assume. 


 \v2
{\bf 2.1.  Hitting time distribution.}
\v2

Let $f^x(n)$ denote the probability that the walk started at $x$  visits the origin at  $n$ for the first time: 
$$f^x(n) = P[\sigma^x_{\{0\}} =n].$$
Put
$$\k_{\a,\ga} =\k_{\a, -\ga} = \frac{(\a-1)\sin \frac{\pi}{\a}}{\Ga(\frac1{\a})\sin \frac{ \pi(\a -\ga)}{2\a} } = \frac{(1-\frac1{\a})\sin \frac{\pi}{\a}}{\mathfrak{p}_1(0) \pi};$$
in particular if $\ga=|2-\a|$, $\k_{\a,\ga} =(\a-1)/\Ga(1/\a)$. 
\v2

\begin{Thm}\label{thm1} \,  For any admissible $\ga$,  as  $n\to\infty$
$$f^0(n) \sim   \frac{ \k_{\a,\ga} c_\circ^{1/\a}}{n^{2-1/\a}}.$$
\end{Thm}
\v2
 When $\ga=0$ (i.e., the limit stable process is symmetric), the above asymptotic form of  $f^0(n)$ is derived  by Kesten \cite{K} in which an asymptotic form  for $\a=1$ is also obtained,  which reads $f^0(n) \sim \pi c_\circ/n(\log n)^2$.

 We write  $x_n$ for $x/n^{1/\a}$.
\begin{Thm}\label{thm2} \,  Let  $|\ga| <2-\a$.  Then, for each $M>1$,  as  $n\to\infty$
\v2
\beqn \label{eq_thm1}
f^x(n) \sim \left\{\begin{array}{lc} {\displaystyle  a^\dagger(x) f^0(n) } \qquad  (x_n \to 0),\\[2mm]
c_\circ \mathfrak{f}^{x_n}(c_\circ)/n \qquad  (\mbox{uniformly for}\; \; 1/M \leq  |x_n| <M).
\end{array}\right.
\eeqn

\end{Thm}

\begin{Thm}\label{thm3}  \, Let  $|\ga|  = 2-\a$. Then as $n\to\infty$  (\ref{eq_thm1}) holds if $x \ga \leq 0$,  and uniformly for $0< \ga x_n <M$,
\beqn\label{eq_thm2}
f^x(n) \sim a^\dagger(x) f^0(n) +  \frac{ |x_n|\mathfrak{p}_{c_\circ}(-x_n)}{n}.
\eeqn
 \end{Thm}

In  case   $|x_n|\to\infty$ an upper bound is provided by the following proposition, where we include a reduced version of that  for the case  $|x_n|<1$ given above. 
\begin{Prop}\label{prop2.1}\, There exists a constant $C$ such that for all  $\gamma$ and $x$,
\[
f^x(n) \leq C( |x_n|^{\alpha-1} \wedge  |x_n|^{-\alpha})/n. 
\]
\end{Prop}

\v2
\begin{rem} \label{rem1}\, (a)
We shall see   (cf. Lemma \ref{lem3.1}(i)) that   as $|x|\to\infty$
\beqn \label{a(x)}
c_\circ a(x) =   \left\{\begin{array}{lll} o(|x|^{\a-1}) \quad  &\mbox{if} \quad  \ga x\to +\infty, |\ga| = 2-\a \\
\k^a_{\a,\ga,\, {\sg x}}\, |x|^{\a-1} \{1+o(1)\} \; \quad&\mbox{otherwise} 
\end{array}\right.
 \eeqn
where     $\k^a_{\a,\ga,\, {\sg x}}$ is a constant (depending on $\a, \ga$ and  ${\sg x}$) which is positive if $\ga\,{\sg x} \neq   2-\a$ and equals $1/\Ga(\a)$ if $\ga\,{\sg x} = -  2+\a$. In particular
\beqn\label{a(-y)}
\mbox{if} \;\; \ga=2-\a, \;\; c_\circ a(y) \sim (-y)^{\a-1}/\Ga(\a) \quad \mbox{as} \;\;  y\to -\infty.
\eeqn
  If $\ga \neq 2-\a$, then
$$f^x(n)\sim  \k_{\a,\ga}  c_\circ^{1/\a}\, a(x) / n^{2-1/\a}
 \sim c_\circ \mathfrak{f}^x(c_\circ n)$$
 as $x\to\infty$ and  $x/n^{1/\a}\to 0$ (Lemma \ref{lem7.1}), so that the two expressions on the RHS of (\ref{eq_thm1}) are asymptotically equivalent to each other in this regime. 
This is contrasted with the first  half of (\ref{eq_thm2})  which  implies that  if $\ga = 2-\a$,
as $x\to +\infty$ under  $x< M n^{1/\a}$ 
\beqn\label{crs_ov}
f^x(n) \sim \left\{\begin{array}{ll} \k_{\a,\ga}  c_\circ^{1/\a}\, a(x)/ n^{2-1/\a} \quad & (a(x)/x >\!\!>  n^{1-2/\a}), \\[2mm]
 \mathfrak{p}_{c_\circ}(- x_n) x/n^{1+1/\a}  \quad & (a(x)/x <\!\!<  n^{1- 2/\a}), 
 \end{array} \right.
 \eeqn
 where  $s <\!\!< t$  means  $t>0$ and $s/t\to 0$. It is noted that
   $a(x), x>0$ is positive if  $P[X\geq 2] >0$ and possibly bounded (see (\ref{a_bdd})).  
\v2
(b)\,   Whenever $|x_n|\to0$  (\ref{eq_thm2})  is valid  for all (admissible) $\ga$, for   if  either $|\ga|<2-\a$ or $x\ga <0$, then  in view of (\ref{a(x)}) the second term on the RHS of (\ref{eq_thm2}) is negligible as $x_n\to0$ in comparison to the first so that it reduces to the first case of (\ref{eq_thm1}).
\v2
(c)\,  If $\ga = 2 -\a$ (when the limiting  stable process has no negative jumps), then it holds that 
\beqn\label{eqR(b)}
\mathfrak{f}^x(t) = xt^{-1}\mathfrak{p}_t(-x)\quad \mbox{ for} \quad x>0
\eeqn
  (cf., e.g.,  \cite[Corollary 7.3]{Bt}), which shows that  in the regime $1/M<|x_n|<M$
 the asymptotic forms of $f^x(n)$ given  in Theorems \ref{thm3} and \ref{thm2} are equivalent to each other 
in view of   the scaling relation (\ref{scl_rl})  below.  Thus 
 for all $|\ga|\leq 2-\a$, as $n\to\infty$
\beqn\label{eq_0}
f^x(n) \sim c_\circ \mathfrak{f}^x(c_\circ n) \qquad \mbox{uniformly  for}\;\; 1/M\leq |x_n| \leq M.
\eeqn
\v2

(d) \, It seems hard to improve the estimate for $|x_n|>1$  given in Proposition \ref{prop2.1} under  (\ref{f_hyp}) only. However,  if we assume some additional regularity condition on  $p(x)$ as $x\to -\infty$ (resp. $+ \infty$) the upper bound of $f^x(n)$  for $x_n>1$ (resp. $x_n<-1$)   is improved to $|x_n|^{-\alpha -1}/n$. 
\end{rem}
\v2\v2

The density function $\mathfrak{f}^x(t)$ satisfies the scaling relation
\beqn\label{scl_rl} 
\mathfrak{f}^x(c_\circ t) = \mathfrak{f}^{x/t^{1/\a}}(c_\circ)/t = \mathfrak{f}^1(c_\circ t/x^\a)/x^\a.
\eeqn 
In case $\ga =|2-\a|$,  expansions of $\mathfrak{f}^x(t)t$ into power series of $x/t^{1/\a}$ are known. Indeed, if $\ga = 2 -\a$,  owing to (\ref{eqR(b)})   the power series expansion for $x>0$ 
is obtained from that of  $t^{1/\a}\mathfrak{p}_t(-x)$ which is found in   \cite{F}, while  for $x<0$, the  series expansion is recently derived by Peskir \cite{P}. Peskir's result  implies
\beqn\label{Psk}\mathfrak{f}^x(t) =[\Ga(\alpha -1) \Ga(1/\alpha)]^{-1} (-x)^{\alpha-1} t^{-2+1/\alpha}\{1+ O([-x/t^{1/\alpha}]^{2- \alpha})\} 
 \quad   (x<0)
 \eeqn
for $x=O(t^{1/\a})$. For $|\ga| <2-\a$ a corresponding asymptotic form is obtained as a by-product of the proof of Theorem \ref{thm2}.  As a consequence we have  the following  corollary. 

\begin{Cor}\label{cor1} \, As $t\to\infty$
$$  \mathfrak{f}^1(t) \sim   \left\{ \begin{array}{ll}  [-1/\Ga(-1/\a)]/t \quad &\mbox{if}\quad 
\ga=2-\a, \\[2mm]
 \k^{\mathfrak{f}}_{\a,\ga}/t^{2-1/\a} \quad &\mbox{if}\quad 
\ga\neq 2-\a, 
\end{array} \right.
$$ 
where
\beq
\k^\mathfrak{f}_{\a,\ga} =\frac{\sin \pi/\a}{ \pi \mathfrak{p}_1(0)}\int_0^\infty u^{1-\a}\mathfrak{p}_1'(-u)du 
= \frac{\Ga(2-\a)\sin (\pi/\a) \,\sin [\frac\pi2 (\a +\ga)]}{\a \pi^2 \mathfrak{p}_1(0)}.
\eeq
[The last expression shows  that  $\k^\mathfrak{f}_{\a,\ga}$    is positive if  $\ga < 2-\a$ and zero if $\ga =2-\a$ (cf. Lemma \ref{lem7.1}). In case $\ga=-2+\a$ the formula above yields the leading term in (\ref{Psk}).]
\end{Cor}

\v2
{\bf 2.2. Transition probability of the walk killed on $\{0\}$.}
\v2

For a non-empty subset $B\subset \Z$ put
\beqn\label{trans}
p^n_B(x,y) = P[S^x_n=y, \sigma^x_B >n]
\eeqn
(in particular $p^0_B(x,y)=\1(x=y)$ and $p^n_B(x,y)=0$ whenever $n\geq1, y\in B$) and similarly for a closed set  $\De \subset \R$
\[
 \mathfrak{p}^\De_t(\xi,\eta) = P[Y^\xi_t\in d\eta, \sigma_\xi^\De >t]/d\eta.
 \]
($Y^\xi_t =\xi +Y_t$ and $\sigma^\Delta_\xi$ is the first entrance time of $Y^\xi$ into $\De$.) By the scaling law for stable processes we have
\[
  \mathfrak{p}^\De_n(x,y) = n^{-1/\alpha} \mathfrak{p}^{\De/n^{1/\alpha}}_1(x_n,y_n).
\]  
In this subsection we give the results   for the  special case  $B=\{0\}$.  The results in the general case of finite sets  closely  parallel  to them and are given in the last subsection {\bf 2.4}.
\v2
We write $x_n$ (resp. $y_n$) for $x/n^{1/\a}$  (resp. $y/n^{1/\a}$) as before. From Theorems \ref{thm2} and \ref{thm3}  it follows   that for all $\ga$ 
$$
f^x(n)\sim a^\dagger(x)f^0(n) + \frac{|x_n|\mathfrak{p}_{c_\circ}(-x_n)}{n}   \quad(x_n \to 0), 
$$
where the second term  on the RHS
is redundant unless  $|\ga|=2-\a$ and $\ga x>0$ as noted in Remark \ref{rem1}(b).
%
\begin{Thm}\label{thm4}\, Let $|\ga|<2-\a$.  For any $M>1$, uniformly for $|x_n| \vee |y_n|<M$, as $n\to\infty$
\beqn\label{eq_thm3}
p^n_{\{0\}}(x,y) \sim \left\{\begin{array}{lr}{\displaystyle f^x(n)a(-y)} &  (y_n \to 0),
\\[2mm]
{\displaystyle  a^\dagger (x)f^{-y}(n) } \quad &( x_n \to 0, y\neq0),\\[2mm]
\mathfrak{p}^{\{0\}}_{c_\circ n} (x,y) & (|x_n| \wedge |y_n|\geq 1/M).
\end{array} \right.
\eeqn
[The first two formulae on the RHS   are asymptotically equivalent to each other as $x_n\vee y_n\to 0$ but not if $x_n\vee y_n>1/M$.]
\end{Thm}

Recalling Remark \ref{rem1}(a) it follows   that if $|\gamma|<2-\a$, then for  any $M>1$,
\beqn\label{eqT4}
p_{\{0\}}^n(x,y)\asymp |x_ny_n|^{\alpha-1}/n^{1/\alpha}\quad (|x_n\vee|y_n| < M)
\eeqn

\begin{Thm}\label{thm5}\, Let $\ga=2-\a$.  For any $M>1$, uniformly for $|x_n|<M$  and $0 <  y\leq M n^{1/\a}$, as $n\to\infty$
\beqn\label{eq_thm4}
p^n_{\{0\}}(x,y) \sim \left\{\begin{array}{lr}{\displaystyle  f^x(n) a(-y)  } \quad  & (y_n \to 0),
\\[1mm]
{\displaystyle   a^\dagger(x)f^{-y}(n) + \frac{(x_n)_+ K_{c_\circ}(y_n)}{n^{1/\a} }} \quad\;\;  & (x_n \to 0),\\[2mm]
\mathfrak{p}^{\{0\}}_{c_\circ n} (x,y) \quad & \;\;   (|x_n| \wedge y_n\geq 1/M).
\end{array} \right.
\eeqn
where  $K_{t}(\eta) =0 $ ($\eta\leq 0$) and 
\beqn\label{Kt}
K_{t}(\eta) =\lim_{\xi\downarrow 0} \frac1{\xi}\mathfrak{p}^{(-\infty,0]}_{t} (\xi,\eta) \quad (\eta>0).
\eeqn
\end{Thm}

The duality relation  $p^n_{\{0\}}(x,y)= p^n_{\{0\}}(-y,-x)$ ($xy \neq 0$)  gives another  apparently  different   statement  of Theorem \ref{thm5}.    Specializing to the case $|x_n|\wedge |y_n| \to 0$ and incorporating  Theorems \ref{thm2} and \ref{thm3}  we  here write down it  as the following corollary   for convenience of later citations.

\begin{Cor}\label{cor2}\, If  $\ga=2-\a$, 
uniformly for $-Mn^{1/\alpha} \leq x\leq 0$ and $|y_n|< M $,   as $n\to\infty$
\beqn\label{eqC2}
 p^n_{\{0\}}(x,y) 
\sim  \left\{\begin{array}{lr}
a^\dagger(x) c_\circ \mathfrak{f}^{-y_n}(c_\circ)/n           & (x_n\uparrow 0,\, y_n>1/M),   \\[2mm]
   a^\dagger (x)\Big[f^{0}(n)a(-y) + |y_n| \mathfrak{p}_{c_\circ}(y_n)n^{-1} \Big]  \quad  & (x_n \uparrow  0,\, y<0),
\\[2mm]
{\displaystyle   a(-y)f^{x}(n) + \frac{(y_n)_- K_{c_\circ}(-x_n)}{n^{1/\a} }} \quad\;\;  & (y_n \to 0).
\end{array} \right.
\eeqn
\end{Cor}

Note that  (\ref{eqC2}) includes the case $y<0$, $x<0$ that is excluded from (\ref{eq_thm4}). The case $x>0$ and $y<0$  excluded from the both will be discussed after Remark \ref{rem2} below.
If  the walk  is left-continuous in particular, namely if $P[X\leq -2]=0$ (possible for $\gamma=2-\alpha$), then  in case $y<0$,   $a(-y)=0$ and  (\ref{eq_thm4}) cannot hold, its right side  vanishing while the left side being positive for $x\leq 0$. This case however is  included in (\ref{eqC2}).   Similarly   (\ref{eqC2}) for the  case $a(x)=0$ is complemented by (\ref{eq_thm4}). 
If the walk is  not left-continuous,  (\ref{eq_thm4})  (resp. (\ref{eqC2}))  is extended to the case $-M< y<0$ (resp. $0<x <M$). The extension can be trivially made in the course of the proof, although we shall not mention  it.  The same comment applies to several places in the sequel where analogous situations occur.

For $\ga=-2+\a$, the result specialized  to the case $y> 0$ and $|x_n|\wedge |y_n| \to 0$,  is given as  follows:  uniformly for  $x\geq 0$ and $|x_n|\vee |y_n|<M$, as $n\to\infty$ 
 $$p^n_{\{0\}}(x,y) \sim \left\{\begin{array}{ll}
   f^{-y}(n)a(x)  \quad  & (x_n\to 0),\\
{\displaystyle   a^\dagger(x)a(-y)f^{0}(n) + \frac{(y_n)_+ \hat K_{c_\circ}(x_n)}{n^{1/\a} }} \quad\;\;  & (y_n \to 0),
\end{array} \right.
$$
which  is immediately  deduced from  (\ref{eq_thm4}) by  using duality relations:  $p^n_{\{0\}}(x,y) =\hat p^n_{\{0\}}(y,x)$,  $\hat a(x)= a(-x)$ and $\hat f^x(n) = f^{-x}(n)$ and
 $$\hat K_t(\eta) =\lim_{\xi\downarrow 0} \frac{\mathfrak{p}^{(-\infty,0]}_{t} (\eta,\xi)}{\xi}  =\lim_{\xi\downarrow 0} \frac{\mathfrak{p}^{[0,\infty)}_{t} (-\xi, -\eta)}{\xi} \quad(\eta>0),$$ where $\hat\,$ indicates the corresponding functions for the dual walk. 


\begin{rem} \label{rem2}
(a) \,   The same crossover  as described in Remark \ref{rem1}(a)   takes place in (\ref{eq_thm4}) plainly for the first case of it but also  in the second case:   in the both  the crossover occurs around  $a(x)/x \asymp n^{1-2/\alpha}$ as in (\ref{crs_ov}),
and similarly in (\ref{eqC2}) around $a(-y)/y \asymp n^{1-2/\a }$.
\v2
(b) \,  
If \, $\ga=2-\a,$ then  $\mathfrak{p}^{\{0\}}_{t} (x,y) =\mathfrak{p}^{(-\infty,0]}_{t} (x,y)\;\; (x, y>0).$ (Cf. e.g., \cite{Bt}.)


\v2
(c) \, Paralleling   Remark \ref{rem1}(c)  concerning  $f^x(n)$ it holds that for all admissible $\ga$ and for  $|x_n|, |y_n| \in [1/M, M]$,  
$p^n_{\{0\}}(x,y) \sim n^{-1/\a}\mathfrak{p}_{c_\circ}^{\{0\}}(x_n,y_n)$,   whenever $xy>0$. This remains  true in case $xy<0$  if  either  $|\ga|<2-\a$ or $x\ga <0$, but does not   anymore  otherwise,
namely if  either  $x>0, y<0$ and  $\ga=2-\a$ or $x<0, y>0$ and  $\ga=-2+\a$ (see Theorem \ref{thm6} and (\ref{eqL2.1})).

\v2

(d) \, Let $\gamma = 2-\alpha$. The first formula (\ref{eq_thm4})   implies that  if one takes the successive  limit as  first  $x_n\to \xi >0$, $y_n\to \eta >0$ as well as $n\to\infty$ and then   $\xi \vee \eta \to 0$, then
\[
 \frac{p^n_{\{0\}}(x,y)} {f^x(n)a(-y)}=   \frac{p^n_{\{0\}}(x,y) }{ \mathfrak{p}^{\{0\}}_{c_\circ n}(x, y)} \frac{\mathfrak{p}^{\{0\}}_{c_\circ n}(x, y)}{f^x(n)a(-y)} \to 1. 
\]
Since  the limit of the first ratio of the middle member equals 1 by virtue of the second relation of (\ref{eq_thm4}), it therefore follows from Theorem \ref{thm3} that   as $\xi \vee \eta \to 0$ and $n\to\infty$
\beqn\label{eqR(d)}
\frac{ \mathfrak{p}^{\{0\}}_{c_\circ}(\xi, \eta) n^{1-1/\a}}{\xi \mathfrak{p}_{c_\circ}(-\xi) a(-y)} 
\to  1.
\eeqn
On noting  $\mathfrak{p}^{\{0\}}_{c_\circ}(-\xi, -\eta) = \mathfrak{p}^{\{0\}}_{c_\circ}(\eta, \xi)$ and using 
(\ref{a(-y)}) this   shows that  
\beqn\label{R3eq}
\mathfrak{p}^{\{0\}}_{c_\circ}(\xi, \eta ) \sim  \frac{\mathfrak{p}_{c_\circ}(0)}{c_\circ \Ga(\a)} \times \left\{\begin{array}{lr}  \xi \eta^{\a-1} \quad  &(\xi \downarrow 0,  \eta \downarrow 0),\\
-\eta(-\xi)^{\a-1} \quad & (\xi \uparrow 0, \eta \uparrow 0).
\end{array} \right.
\eeqn
\v2

(e) \, 
In the same way as in  (d)
 we deduce   from Theorem \ref{thm4}  that if  $|\ga| <2-\a$,   
 $$\mathfrak{p}^{\{0\}}_{1}(\xi, \eta) 
\sim  \k'_{\a,\ga} \{\sin[\tst12 \pi(\a + (\sg \xi) \ga)]\}\{\sin[\tst12 \pi(\a - (\sg \eta) \ga)]\} |\xi\eta|^{\a-1} \qquad (|\xi|\vee |\eta|\to 0),$$
where $\k'_{\a,\ga}=\k_{\a,\ga}\{\Ga(1-\a)/\pi\}^2$. 
Similarly,  with $\xi>0$ fixed and   $\eta$  tending to zero,  noting   (\ref{eqR(b)})  we see  that if $\ga<2-\alpha$ or $y>0$ (i.e., when $a(-y)/y^{\alpha-1}$ tends to a positive constant),
$$
\frac{ \mathfrak{p}^{\{0\}}_{c_\circ}(\xi, y_n)}{c_\circ \mathfrak{f}^\xi(c_\circ) [a(-y)/n^{1-1/\a}]} 
\to  1  \qquad (y_n \to 0, y\to\infty, \xi>0).
$$

\v2
(f) \, If $\gamma = 2-\alpha$,  then  $\mathfrak{p}_{c_\circ}(0) = 1/c_\circ^{1/\alpha}\alpha \Gamma(1-1/\alpha)$ (see Lemma \ref{lem3.2}) and by (\ref{R3eq}) 
\beqn\label{R3f1}
K_{c_\circ}(\eta) \sim  \frac{\mathfrak{p}_{c_\circ}(0)}{c_\circ \Ga(\a)} \eta^{\alpha-1} \qquad (\eta \downarrow 0),
\eeqn
with the help of  which we deduce from Theorem \ref{thm5} and its corollary that for $|x_n|\vee |y_n|<M$,
\beqn\label{R3f2}
 p^n_{\{0\}}(x,y) 
\asymp  \left\{\begin{array}{lr}
f^x(n)a(-y) \asymp  |y_n|^{\alpha-1}\{a^\dagger(x)n^{-1} + |x_n| n^{-1/\alpha} \}  \quad  & (y >  0),
\\[2mm]
 (|x_n|^{\alpha-1}\vee1)\{a(-y)n^{-1} + |y_n| n^{-1/\alpha} \}  \quad  & (x\leq 0).
  \end{array} \right.
\eeqn
\end{rem}

\v2\v2
 From Theorem \ref{thm5}  is excluded
  the regime  $x>0, \, y<0, \, x\wedge (-y) \to +\infty$ (as noted previously), where     there arises  a difficulty in  estimating  $p^n_{\{0\}}(x,y)$ in general; in below we give a result under an extra assumption on the tail  as $t\to-\infty$ of the distribution function
$$F(t) := P[X\leq t].$$
   In  \cite[Theorem 2(iii)]{Uladd} a criterion for the limit
 $$C^+ :=  \lim_{x\to +\infty}  a(x) \leq \infty$$
(which exists) to be  finite  is obtained.   Under the present assumption on $F$ it say that
 \beqn\label{a_bdd}
 \int_{0}^\infty F(-t)t^{2\a-2}dt <\infty \quad \mbox{and} \quad F(-2)>0
 \eeqn
is  necessary and sufficient  for   $0<C^+<\infty$. Note that   (\ref{a_bdd}) entails $\ga=2-\a$ and the walk is not left-continuous.
 
  \v2
\begin{Thm}\label{thm6}\, Let (\ref{a_bdd}) hold. Then, given  $M>1$, uniformly for $-M<y_n<0< x_n <M$,
\v2
{\rm (i)} \quad $p^n_{\{0\}}(x,y) \sim 
{\displaystyle 
a^\dagger(x)a(-y)f^0(n) +\frac{ a^\dagger(x)|y_n|\mathfrak{p}_{c_\circ}(y_n) +  a(-y)x_n\mathfrak{p}_{c_\circ}(-x_n)} { n} }$ \\
\qquad \qquad \qquad \qquad \qquad\qquad\qquad\qquad\qquad \qquad\qquad \qquad\qquad \qquad\qquad $(x_n\wedge (-y_n) \to 0),$
\v2
{\rm (ii) }\quad 
$p^n_{\{0\}}(x,y) \sim  {\displaystyle  \frac{C^+(x_n-y_n)\mathfrak{p}_{c_\circ}(y_n -x_n)}{ n}=C^+ c_\circ} \mathfrak{f}^{\, x-y}(c_\circ n)
 \qquad\quad   (x_n\wedge (-y_n) >  1/M) $
\v2\n
  as $n\to\infty$.
\end{Thm}

An application  of  Theorem \ref{thm6} leads to the next result which  exhibits  a way  the condition   $C^+<\infty$  is reflected in   the behaviour of   the walk $S^x$, $x>0$: conditioned on   $S^x_n=-x$ it enters $(-\infty,-1]$ without visiting the origin  \lq continuously' or by a very long jump for large $x$  according as $C^+$ is finite or not.  Exactly the same behaviour of the pinned walk is observed in \cite{U1dm_f} in the case $E|X|^2<\infty$  but with the condition (\ref{a_bdd}) replaced by $E[|X|^3; X<0] <\infty$ which is equivalent to $\lim_{x\to-\infty} [a(x)-x/\sigma^2] <\infty$. 

  \begin{Prop}\label{prop2.2}  For each  $M\geq 1$,   under the constraint   $-M\sqrt n <y<0< x< M\sqrt n$  
\beq
&& P[S^x_{\sigma(-\infty,0]}<- R \,|\,\sigma^x_{\{0\}}>n, S^x_n =y] 
\\[2mm]
&&\quad  \longrightarrow  \left\{ 
\begin{array}{ll}
0 \quad\mbox{ as}  \;\;  R \to\infty\;\qquad  \mbox{uniformly for} \;\; x, \,y\; &\mbox{ if}\quad C^+ <\infty, \\[1mm]
1 \quad\mbox{ as}  \; \;  x\wedge (-y) \to\infty \;\; \mbox{ for each}\; R>0\;\; &\mbox{ if}\quad C^+ =\infty.
\end{array} \right.
\eeq
\end{Prop}
\v2

We state  the  following upper bound as a proposition, a unified but partly reduced version obtained  by combining theorems above and the results  in  Section 5.  
\begin{Prop}\label{prop2.3} \,
{\rm (i)} \, For all  admissible $\ga$ and $M>1$, there exists a constant $C_M$ such that for all $n\geq 1$ and $x\in \Z$, 
\[ 
p^n_{\{0\}}(x,y)  \leq C_M  ((|x_n|\vee 1)^{\alpha-1} \wedge |x_n|^{-\alpha}) |y|^{\alpha -1} 
\quad\mbox{if}\quad  |y_n|<M. \]

 {\rm (ii)}  If $\ga =2-\a$, there exists a constant $C$ such that for all $x, y\in \Z$,
\beqn\label{eqL2.1}
p^n_{\{0\}}(x,y) \leq  
C\bigg[\frac{a^\dagger(x)a(-y)}{ n^{2-1/\a}} +\frac{ a^\dagger (x)\{(y_n)_- \wedge 1\}+  a(-y)\{(x_n)_+\wedge 1\}} { n} \bigg]. 
\eeqn
 \end{Prop}
\v2
 (i) is the same as Lemma \ref{lem5.2}. (ii)    follows from Theorems \ref{thm4} and \ref{thm5} in case   $|x_n|\vee |y_n| <1$, 
from Proposition \ref{prop5.2} in case $|x_n|\vee |y_n| \geq 1$ with $xy<0$, from Lemmas \ref{lem5.1} and \ref{lem5.2} in case $|x_n|\wedge  |y_n|\leq 1\leq |x_n|\vee |y_n|$ with $xy\geq 0$ and the bound  $p^n(x) \leq Cn^{1/\alpha}$ (entailed by the local limit theorem) in case $|x_n|\wedge  |y_n|\geq 1$ with $xy\geq 0$. 


\v2
{\bf 2.3. Comparing $p^n_{\{0\}}(x,y)$ and} $p^n_{(-\infty, 0]}(x,y)$. \,  

\v2
 Let  $V_{{\rm as}}$ (resp. $U_{{\rm ds}}$) denote the renewal function of weakly ascending (resp. strictly descending) ladder height process of the walk $S$ and $Q_{t}(\eta)$ and $\hat Q_t(\eta)$, $\eta \geq 0$  the distribution functions of the stable meander of length $t$ at time $t$ for $Y$ and $-Y$, respectively (see (\ref{Bt}) for the definition). 
Doney \cite{D} obtains  an elegant asymptotic formulae of $p^n_{(x,\infty)}(0,y)$ ($x\geq 0$),
  which under the present assumption and with our notation may be rewritten as 
\beqn\label{Doney}
p^n_{(-\infty, 0)}(x,y)  \sim \left\{\begin{array}{lr}
{\displaystyle \frac{U_{{\rm ds}}(x)V_{{\rm as}}(y) \mathfrak{p}_{c_\circ}(0)}{ n^{1+1/\a}} } \quad & (x_n\vee y_n \to 0), \\[3mm]
{\displaystyle  
\frac{V_{{\rm as}}(y)  P[\sigma^0_{[0, +\infty)}>n]  \hat Q'_{c_\circ}(x_n)}{  n^{1/\a}} } \quad  & (y_n \to 0, x_n >1/M ),\\[3mm]
{\displaystyle  
\frac{ U_{{\rm ds}}(x) P[\sigma^0_{(-\infty, -1]}>n] Q'_{c_\circ}(y_n)}{ n^{1/\a}} } \quad  & (x_n \to 0, y_n>1/M),\\[2mm]
\mathfrak{p}^{(-\infty,0]}_{c_\circ n} (x,y) \quad & \;\;   (x_n \wedge y_n\geq 1/M) \end{array} \right.
\eeqn 
 by using the duality relation.

In the regime $x \asymp y \asymp n^{1/\a}$, where  $\mathfrak{p}^{\{0\}}_{c_\circ n} (x,y)/\mathfrak{p}^{(-\infty,0]}_{c_\circ n} (x,y) = \mathfrak{p}^{\{0\}}_{c_\circ } (x_n,y_n)/\mathfrak{p}^{(-\infty,0]}_{c_\circ } (x_n,y_n)  \asymp 1$, we have  $p_{\{0\}}^n(x,y) \asymp p_{(-\infty,0]}^n(x,y)$  for all  values of $\ga$.

If  $|\ga|< 2-\a$,  then $V_{{\rm as}}$ and  $U_{{\rm ds}}$ vary regularly with exponents which are both larger than $\a-1$ and whose sum equals $\alpha$ (cf. \cite{D}) and each of the products  $V_{{\rm as}}( n^{1/\a})  P[\sigma^0_{[0, +\infty)}>n]$ and $U_{{\rm ds}}(n^{1/\a}) P[\sigma^0_{(-\infty, -1]}>n]$ approaches to a positive constant as $n\to\infty$. These lead to an estimate for $p_{(-\infty,0]}^n(x,y)$ analogous  to that for $p_{\{0\}}^n(x,y) $  given in (\ref{eqT4})  and comparing them yields
\[
p_{(-\infty,0)}^n(x,y) / p_{\{0\}}^n(x,y) \to 0  \quad  \mbox{as} \quad x_n\wedge y_n \to 0. 
\] 
In case $|\ga| = |2-\a|$ we need to take a  closer look at  the situation that turns out to be precisely   parallels   the crossover  mentioned in Remark \ref{rem2}(a) as given by (\ref{comp}) below for $\ga=2-\a$. 

 Let  $\ga=2-\a$.  Then  
$p_{\{0\}}^n(x,y) \sim  p_{(-\infty,0)}^n(x,y) $ for $x_n \wedge y_n >1/M$ in view of (\ref{eq_thm4}) and (\ref{Doney}) (see also Remark \ref{rem2}(b)). According to  \cite[Theorems 2 and 9]{R} 
\beqn\label{LH}
U_{{\rm ds}}(x)\sim x L(x)\quad \mbox{and} \quad V_{{\rm as}}(x) \sim \k^{V} x^{\a-1}/L(x)\quad (x\to\infty)
\eeqn 
with a positive constant  $\k^{V}$ and a slowly varying $L(x)$ that tends to zero as $x\to\infty$. 
(More information is found in  Remark \ref{rem3}(b) given below. A condition sufficient for $L$ to be asymptotic to a positive constant is considered in  Remark \ref{rem10}.)  

Let $\hat Z$ be the first  (strictly) descending ladder height, namely $ \hat Z = S_{\sigma(-\infty, 0)}$, and suppose 
 \beqn\label{Z}
E|\hat Z|<\infty.
 \eeqn Then  $L$ may be taken to be  the  constant $1/E|\hat Z|$ owing to the renewal theorem
and letting first $x_n\to \xi>0$ and $y_n\to \eta>0$ and then  $\xi\downarrow 0$ or $\eta\downarrow 0$ in (\ref{eq_thm4}) and  (\ref{Doney})   we  see  that uniformly for  $x,y>0$ and $x\vee y \leq Mn^{1/\alpha}$  as $x\wedge y\to\infty$
\beqn\label{zzz}
 p_{(-\infty,0)}^n(x,y) \sim \left\{\begin{array} {ll}    
 y_n^{\alpha-1} x_n \mathfrak{p}_{c_\circ}(-x_n)/n^{1/\alpha}c_\circ\Gamma(\alpha) \quad &(y_n \to 0), \\[1mm]
  x_n K_{c_\circ}(y_n)/n^{1/\alpha} \quad & (x_n \to 0),
 \end{array} \right.
 \eeqn
(this is confirmed by making an elementary computation (see Remark \ref{rem3}(c))  which however is not needed) and   compare this  with (\ref{eq_thm4}) to deduce that  uniformly for $0\leq x,  y < Mn^{1/\a}$, as $n\to\infty$ 
\beqn\label{comp}
p_{\{0\}}^n(x,y) \sim p_{(-\infty,0)}^n(x,y) + a^\dagger(x) f^0(n) a(-y).  
\eeqn
 According to Kesten \cite{K} (see Remark \ref{rem4C}  of the next subsection) this 
asymptotic relation with $x, y$ fixed (when the first term on the RHS is superfluous)  is valid  for every recurrent walk that is strongly aperiodic and having   $\sigma^2=\infty$.   It is quite plausible that (\ref{comp})   holds 
for $x, y$ subject to the same constraint as above for every such random walk on $\Z$ with $E|\hat Z|<\infty$. If $E|\hat Z|=\infty$,  $ p_{(-\infty,0)}^n(x,y)$ may depend on $L$   in the regime $x_n\wedge y_n\to0$
 while $p_{\{0\}}^n(x,y)$ does not, and (\ref{comp}) must be violated.

 \v2
\begin{rem}\label{rem3} \,   \,(a)\,    Let  $\ga=2-\a$ and $K_t$ be given in (\ref{Kt}).  Then for $x, y>0$ 
 \beqn\label{R2.3a}
  x\mathfrak{p}_t(-x) = t\mathfrak{f}^x(t)  = t^{1/\a} \hat Q_t'(x)/\Ga(1/\a) \quad \mbox{and} \quad
    K_t(y)= \a \mathfrak{p}_t(0) Q'_{t}(y)
    \eeqn
(see Lemma \ref{lem7.4}), and on  rewriting the  first equality and   using  (\ref{R3f1})
\beqn\label{Q}
 \hat Q'_{t}(\eta) =t^{-1/\a} \Ga(1/\a)\mathfrak{p}_{t}(-\eta)\eta \;\; (\eta>0) \quad \mbox{and}\quad Q_t'(\eta) \sim \eta^{\a-1}/t \alpha \Ga(\a) \;\;   (\eta\downarrow 0).
 \eeqn
 Note  that $Q_t(\eta) = Q_1(\eta/t^{1/\a})$, entailing  $Q_n'(y) = Q_1'(y_n)/n^{1/\a} $. 
 \v2
(b)\,    It is known \cite[Eq(15), Eq(31)]{VW} that for some positive constant $b$,
$$ P[\sigma^0_{(-\infty, -1]}>n] \sim  b/U_{dc}( n^{1/\a}),  \; \mbox{and}$$
$$ P[\sigma^0_{(-\infty, -1]}>n]  P[\sigma^0_{[0, +\infty)}>n]  \sim \b/n \quad\mbox{with}\;\;  \b := 1/\Ga(\rho)\Ga(1-\rho),$$
where  $\rho:= \lim_{n} n^{-1}\sum_{k=1}^n P[S_k>0] =\frac12(1-\ga/\a)$ (cf. (\ref{rho})).
We derive in below that if $\ga=2-\a$, then
\beqn\label{b.kappa}
 b= \frac{1}{c_\circ^{1/\a}\Ga(1-1/\a)}  \qquad \mbox{and} \qquad   \k^{V} = \frac1{c_\circ\Ga(\a)}
\eeqn
($\k^{V}$ appears in (\ref{LH})), the former identity (together with $\rho =1-1/\alpha$ and (\ref{LH}))
entailing
\beqn\label{eqR2.3}
P[\sigma^0_{(-\infty, -1]}>n]\sim  \frac{1/ \Ga(1-1/\a)}{(c_\circ n)^{1/\a}L(n^{1/\a})}\quad
\mbox{and}\quad 
P[\sigma^0_{[0, +\infty)}>n ]\sim   \frac{c_\circ^{1/\a}L(n^{1/\a})}{n^{1-1/\a}\Ga(1/\a)}.
\eeqn
 (\ref{b.kappa}) as well as what are stated prior to it  is valid even  if  $X$ belongs to a domain of attraction (instead of a normal domain).
 For the derivation,  employing  (\ref{eq_lem6.4}) that says $\mathfrak{p}^{\{0\}}_{c_\circ}(\xi,\eta)\sim \a \mathfrak{p}_{c_\circ}(0)Q_{c_\circ}'(\eta)\xi$ ($\xi \downarrow 0, \eta>0$)  we deduce from  the third and fourth cases of (\ref{Doney}) that 
 $$bU_{{\rm ds}}(x)/U_{{\rm ds}}(n^{1/\a}) \sim  \a \mathfrak{p}_{c_\circ}(0)x_n \qquad (\xi =x_n \downarrow 0),$$
which immediately leads to  $b= \a\mathfrak{p}_{c_\circ}(0) = 1/c_\circ^{1/\a}\Ga(1-1/\a)$ (see Lemma  \ref{lem3.2} for the second equality).  In a similar way employing (\ref{eq_lem6.41}) 
and the second formula of (\ref{Doney}) we derive  $\k^{V}=1/c_\circ\Ga(\a)$ as required. %

\v2
 (c) \, By (\ref{R2.3a})  and  (\ref{eqR2.3})  we observe that as $n\to\infty$
\beqn\label{V/a} \begin{array}{lr}
P[\sigma^0_{[0, +\infty)}>n]  \hat Q'_{c_\circ}(x_n)  \sim L(n^{1/\alpha})x_n \mathfrak{p}_{c_\circ}(-x_n)/n^{1-1/\alpha} \quad &( x>0),  \\[2mm]
P[\sigma^0_{(-\infty, -1]}>n] Q'_{c_\circ}(y_n) \sim K_{c_\circ}(y_n)/ n^{1/\alpha}L(n^{1/\alpha}) \quad &(y>0),
\end{array}
\eeqn
which together directly derive (\ref{zzz}) from (\ref{Doney}) as noted before.

\end{rem}

 \v2
{\bf 2.4. Extension to the process killed on a finite set}
\v2
Let $A$ be a finite subset of $\Z$. Suppose for simplicity that for some $M>1$
\beqn\label{HA}
g_A(x,y)>0\quad \mbox{if} \quad  |x|\wedge |y| >M, 
\eeqn
where for a non-empty  $B\subset \Z$, $g_B$ denotes the Green function for the walk killed on $B$:
\beqn\label{green}
g_B(x,y) = \sum_{n=0}^\infty p_B^n(x,y).
\eeqn
Under the condition $\sigma^2=\infty$ there exists 
 \beqn\label{def_u}
u_A(x) =\lim_{|y|\to\infty} g_A(x,y)
\eeqn
\cite[T30.1]{S}.  $u_A$ is positive and  harmonic for the killed walk:   $u_A(x)= \sum_{z\notin A}p(z-x)u_A(z)>0$ for all $x\in \Z$. 
Put  
\[
f^x_A(n) = P[\sigma^x_A=n].
\]
In order to obtain the asymptotic form of $f^x_A(n)$ we may  simply replace $f^x (n)$ and  $a^\dagger(\cdot)$ by $f^x_A(n)$ and  $u_A(\cdot)$, respectively,  in Theorems \ref{thm2} and \ref{thm3}, the resulting formula being valid  in the same range of variables so that uniformly for $|x|<M n^{1/\alpha}$, as $n\to\infty$
\beqn\label{f/A}
f_A^x(n) \sim\left\{ \begin{array}{ll}  u_A(x)f^0(n)\quad &\mbox{if} \quad x =o(n^{1/\a})  \;\, \mbox{and} \;\, |\gamma|<2-\alpha,\\
{\displaystyle  u_A(x)f^0(n) +\frac{(x_n)_\pm\mathfrak{p}_{c_\circ}(-x_n)}{n} } \quad &\mbox{if} \quad x_\mp =o(n^{1/\a})  \;\, \mbox{and} \;\, |\gamma|=2-\alpha,
 \end{array}\right.
\eeqn
where  in the symbols  $\pm$ and $\mp$  the upper  (resp. lower) sign prevails if  $\ga>0$  (resp. $\gamma <0$). 
 After  virtually the same replacement [replace $f^\cdot (n)$  by $f_A^{\cdot} (n)$, $C^+$ by $C_A^+ :=  \lim_{x\to\infty}  u_A(x)$ and $a(-y)$ by $u_{-A} (-y)\1(y\notin A)$] Theorems \ref{thm4} to \ref{thm6}  with
$p^n_A(x,y)$ in place of  $p^n_{\{0\}}(x,y)$  remain true.  In case  $\ga =2-\a$ in particular, 
 the result  corresponding to   Theorem \ref{thm5} read 
 \beqn\label{p/A} 
p^n_{A}(x,y) \sim \left\{\begin{array}{lr}{\displaystyle   f^x_A(n)u_{-A}(-y)  } \quad  & (y_n \to 0),
\\[1mm]
{\displaystyle   u_A(x)f^{-y}_{-A}(n) + \frac{(x_n)_+ K_{c_\circ}(y_n)}{n^{1/\a} }} \quad\;\;  & (x_n \to 0),\\[2mm]
\mathfrak{p}^{\{0\}}_{c_\circ n} (x,y) \quad & \;\;   (|x_n| \wedge y_n\geq 1/M).
\end{array} \right.\eeqn 
uniformly for $|x_n|<M$ and $-M<y<Mn^{1/\a}$.

Note that from the definition of $g_A(x,y)$ it follows  that $u_A$ is the probability distribution of  the hitting place of $A$ by the dual walk \lq started at infinity'. 
From  (\ref{f/A})  we therefore  infer
 $$\sum_{z\in A}f_A^z(n) \sim f^0(n),$$
  which relation is observed in \cite{K} when $\ga=0$.
By a similar consideration  or  by the identity
\beqn\label{HD_st}
 P[\sigma_A^x =n, S^x_n =y] = \sum_{z\notin A} p^{n-1}_{A}(x,z)p(y-z) \qquad (y\in A)
\eeqn
one deduces   from (\ref{p/A})  the following asymptotic form of space-time hitting distribution.
\begin{Cor}\label{cor3}  \, Uniformly for $|x_n|<M$, as $n\to\infty$
$$ P[\sigma_A^x =n, S^x_n =y]  \sim f_A^x(n)u_A(-y)  \quad (y\in A).$$
\end{Cor}
\v2
\begin{rem}\label{rem4C}\, As mentioned in Introduction Kesten \cite{K} obtained  asymptotic formulae of  $p_A^n(x,y)$ with $x, y$ fixed for a large class of random walks on multidimensional lattices $\Z^d$, which if specialized to  one-dimensional recurrent walk may  read in the present notation
$$\lim_{n\to\infty}p_A^n(x,y) \Big/ \sum_{z\in A} f_A^z(n) = u_A(x)u_{-A}(-y) \quad (y\notin A),$$ 
provided the walk is  strongly aperiodic and  having infinite variance. 
\end{rem}

The rest of the paper is organized as follows.The proofs of Theorems \ref{thm2} and \ref{thm3} are given  in Section  3  and those of Theorems \ref{thm4} and \ref{thm5}  in section 4.  In Section 5 some estimations of $p^n_{\{0\}}(x,y)$ are made in case $xy <0$ and, for this purpose,  beyond the regime $|x|\vee |y|=O(n^{1/\a})$:  Propositions \ref{prop5.1} and \ref{prop5.2} given there  provide a lower and upper bound, respectively and Theorem \ref{thm6} and Proposition \ref{prop2.2} are   proved after them;  Proposition \ref{prop2.1} is proved at the end of Section 5 where we provide to this end an upper bound in case  $|y|=O(n^{1/\alpha})$ and $|x_n|\to\infty$. In Section 6 the results  are extended to those for an arbitrary finite set instead of the single point set $\{0\}$. In Section 7 we deal with   the limit stable process and present some properties of $\mathfrak{f}^\xi(t)$ and $\mathfrak{p}_t^{\{0\}}(x,y)$.
In the last section we give
miscellaneous consequences  of the present assumption on the random walk that are derived from  the general theory: they are  (1) condition (\ref{f_hyp})  expressed  in terms of  the tails of $F$ and some related facts, (2)   some upper bounds of $p^n(x)$  for $|x| > n^{1/\a}$  and  (3) 
\lq escape probabilities' from the origine.

\section{Estimation of $f^x(n)$}

{\bf 3.1.} 
In several places in this subsection we  shall apply following  identity 
\beqn\label{eq3.000}
\int_0^\infty \left\{\begin{array}{c}\cos u\\
\sin u \end{array}\right\} \frac{du}{u^\nu}
 =  \left\{\begin{array}{l}  \Ga(1-\nu) \sin\frac12\pi\nu \qquad 0<\nu<1 \\[1mm]
 \Ga(1-\nu) \cos\frac12 \pi\nu \qquad 0<\nu<2,
 \end{array}\right.
 \eeqn
where  for $\nu=1$, $ \Ga(1-\nu) \cos\frac12 \pi\nu$ is understood to be $\frac12 \pi$,  its limit value  (cf. \cite[pp.10, 68]{E},  \cite[p.260]{WW}).

\begin{lem}\label{lem3.1}
\, Put  $\kappa^a_{\a,\ga,\pm} = -\Ga(1-\a) \pi^{-1} \sin [\tst12\pi (\a\pm\ga)]$. Then 
\v2
{\rm (i)}\qquad \;
${\displaystyle    \lim_{x\to \pm\infty} \frac{c_\circ a(x)}{|x|^{\a-1}} = -\kappa^a_{\a,\ga,\pm},}$
\v2
{\rm (ii)} \qquad ${\displaystyle   \lim_{x\to \pm\infty} c_\circ \{a(x+1)-a(x)\}|x|^{2-\a} = \pm (\a-1)\kappa^a_{\a,\ga,\pm}; }$
\v2\n
 $\kappa^a_{\a,\ga,\pm} >0$ if $|\ga| <2-\a$,  and  $\kappa^a_{\a,\ga,\pm}  = 0$ or $1/\Ga(\a)$  according as $\pm \ga >0$ or $<0$  if $|\ga| =2-\a$; 
 in particular
if $\ga=2-\a$, then $c_\circ a(x) \sim (-x)^{\a-1}/\Ga(\a)$ as $x\to -\infty$ and $= o(|x|^{\a-1}) $   as $x\to +\infty$. 
\end{lem}
\n
\pf\, Although a little more general result of (i)  is given in \cite[Section 6.1]{Upot}, we give its proof, which is partly used in the proof of (ii).   From   the second of formula  (\ref{eq3.000})   one obtains  
\beqn\label{eq3.00} 
\int_0^\infty \left\{\begin{array}{c}1-\cos u\\
\sin u \end{array}\right\} \frac{du}{u^\a}
 =  \left\{\begin{array}{l}  -\Ga(1-\a) \sin\frac12 \pi\a, \\[1mm]
 \Ga(1-\a) \cos\frac12\pi\a,
 \end{array}\right.
 \eeqn
 where for the first formula one takes $\nu=\alpha -1$ and performs  integration by parts). 
 In the representation 
 $a(x)= \frac1{2\pi}\int_{-\pi}^\pi (1-e^{ix \th})(1-\phi(\th))^{-1}d\th$ we replace $1-\phi(\th)$ by $c_\circ\psi(\th)$,  its principal part about zero,  
 and compute the resulting integral. 
 Changing a variable we have
\beqn\label{L3.1}
\int_{-\pi}^\pi\frac{1-e^{ix\th}}{c_\circ \psi(\th)} d\th =\frac{|x|^{\a-1}}{c_\circ}\int_{-\pi|x|}^{\pi |x|}
\frac{1-\cos u \mp i\sin u}{\cos \frac12 \pi\ga +i(u/|u|) \sin \frac12 \pi\ga} |u|^{-\a}du\quad (\pm = x/|x|),
\eeqn
which  an easy computation with the help of  (\ref{eq3.00}) shows  to be asymptotically equivalent  as $|x|\to \infty$ to
$$
 \frac{-2\Ga(1-\a)  |x|^{\a-1}}{ c_\circ}\Big[\cos\frac{\pi\ga}{2} \sin \frac{\pi\a}{2} \pm \sin\frac{\pi\ga}{2} \cos\frac{\pi\a}{2} \Big]. $$
The combination of the  sine's and cosine's in the square brackets being equal to  $\sin[\frac12\pi(\a\pm \ga)]$ we find the equality  (i), provided  that the replacement mentioned at the beginning causes only a negligible term of the magnitude  $o(|x|^{\a-1})$, but this   is
assured from the way of computation carried out above since  the integrand in the RHS integral in (\ref{L3.1}) is summable on $\R$.  

   For the proof of (ii) it suffices to show that
\beqn\label{L3.11}
\int_{-\pi}^\pi e^{ix \th}(1-e^{i\th})\bigg[\frac1{1-\phi(\th)}- \frac1{\psi(\th)}\bigg]d\th =o(|x|^{\a-2}),
\eeqn
since  $a(x+1)-a(x)= \frac1{2\pi}\int_{-\pi}^\pi e^{ix \th}(1-e^{i\th})(1-\phi(\th))^{-1}d\th$ and  this integral with $\psi(\th)$ replacing  $1-\phi(\th)$ is asymptotically equivalent to $ \pm [(\a-1)\kappa^a_{\a,\ga,\pm}/c_\circ] |x|^{\a-2}$ as one sees by looking at the increment of the RHS of  (\ref{L3.1}). Because of  the fact that if $\psi(\th) = \{1-\phi(\th)\}(1+\de(\th))$ then  $\de'(\th) \th \to 0$ ($\th \to 0$) (cf. (\ref{F/psi})), the relation 
(\ref{L3.11}) is shown in a usual way. 
\qed

\v2
{\bf 3.2.}  For evaluation of $f^x(n)$ we follow \cite{Uht}, although therein  the walk  is assumed to have finite variance. 
Set
$$\pi_x(\tau)=\frac1{2\pi}\int_{-\pi}^\pi \frac{e^{-ix \th}}{1-e^{i\tau}\phi(\th)}d\th~~~~~~(\tau\neq 0, x\in \Z)$$
and
$$ f^\wedge_x(\tau)=\sum_{n=1}^\infty f^x(n)e^{in\tau}.$$
Since  $|\phi(\th)|<1$ for $0<|\th|\leq \pi$ by aperiodicity of the walk, the function $1-e^{i\tau}\phi(\th)$
does not vanish in $\R\times [-\pi,\pi]$ except at $\tau=\th =0$.
We have the following identities
\beqn\label{eq3.0}
\pi_x(\tau)=\lim_{r\uparrow 1}\frac1{2\pi}\int_{-\pi}^\pi \frac{e^{-ix \th}}{1-r e^{i\tau}\phi(\th)}d\th=\lim_{r\uparrow 1}\sum_{n=0}^\infty p^n(x)e^{i\tau n}r^n;\qquad\; 
\eeqn
\beqn\label{hat_f}
 f^\wedge_x(\tau)=-\frac{\1(x=0)}{\pi_0(\tau)}+\frac{\pi_{-x}(\tau)}{\pi_{0}(\tau)},\; \mbox{especially} \quad  f^\wedge_0(\tau)=1-\frac1{\pi_0(\tau)}.
\eeqn
Note that $\Re \,\pi_x(\tau)$ and $\Re  f^\wedge_x(\tau)$  are  even,  while $\Im\, \pi_x$,  $\Im f_x^\wedge(\tau)$ are odd. Obviously $\Re  f^\wedge_x(\tau) $ equals the cosine series 
$\sum_{n=1}^\infty f^x(n) \cos nx$ and
we shall use the following inversion formulae
\beqn\label{inv_f}
 f^x(n) =\frac{1}{\pi}\int_{-\pi}^\pi f^\wedge_x(\tau)\, \cos n\tau d\tau= \frac{2}{\pi}\int_{0}^\pi \Re \, f^\wedge_x(\tau)\, \cos n\tau d\tau.
 \eeqn
 Note that $\pi_0(\tau)$  vanishes nowhere on $[-\pi,\pi]$  and is smooth off the origin, and that
\beqn\label{eq3.2}
1- e^{i\tau}\phi(\th) = -i\tau + c_\circ \psi(\th) + O(\tau^2+|\tau||\th|^\a) + o(|\th|^\a)
\eeqn
and for some constant $C>0$
$$|1- e^{i\tau}\phi(\th) | \geq C^{-1} (|\tau| + |\eta|^\a) \qquad (-\pi<\tau<\pi, \, -\pi<\th<\pi).$$
We shall compare $\pi_x(\th)$ with the corresponding function for the limit stable process given by 
$$\pi_x^\infty(\tau) := \frac1{2\pi}\int_{-\infty}^\infty \frac{e^{-ix\th}d\th}{-i\tau + c_\circ \psi(\th)}$$
($-i\tau + c_\circ \psi(\th)= \int_0^\infty E[e^{i\th Y_t}]e^{i\tau t}dt$ ($\th\neq 0$) corresponding to $1-e^{i\tau}\phi(\th) = \sum E[e^{i\th S_n}]e^{i\tau n}$).
\begin{lem}\label{lem3.2} \,  As $\tau \to 0$
\beqn\label{eq3.3}1/\pi_0(\tau) \sim 1/\pi^\infty_0(\tau),\quad   [1/\pi_0]'(\tau) \sim [1/\pi^\infty_0]'(\tau) \quad\mbox{and} \quad  [1/\pi_0]''(\tau) \sim [1/\pi^\infty_0]''(\tau)
\eeqn
where $'$ indicates the differentiation; and for $\tau\neq 0$,
$$\pi_0^\infty (\tau)  =  \mathfrak{p}_1(0)c_\circ^{-1/\a} \Ga(1-1/\a)ie^{-i\pi/2\a}\frac{|\tau|^{1/\a}}{\tau}, \quad \mathfrak{p}_1(0) = \frac{\Ga(1/\a)}{\pi\a} \sin \frac{\pi(\a-\ga)}{2\a};$$
in particular, $(d/d\tau)^j [1/\pi_0(\tau)] = O(|\tau|^{-1/\alpha +1-j})$ ($\tau \to 0, j=0,1, 2$) and  if $|\ga|= 2-\alpha$,  $\mathfrak{p}_1(0)= 1/\a \Ga(1-1/\a)$.
\end{lem}
\n
\pf\, In view of  (\ref{eq3.2})  it is easy to deduce from the defining expression of $\pi_0(\tau)$ that  $\pi_0(\tau)\sim \pi_0^\infty(\tau)$, $\pi_0'(\tau)\sim {\pi_0^\infty}'(\tau)$ and $\pi_0''(\tau)\sim {\pi_0^\infty}''(\tau)$, which show (\ref{eq3.3}). The expression of $\mathfrak{p}_1(0)$ is  obtained by specializing the series expansion  of  $\mathfrak{p}_1(x)$ as found in 
 e.g. \cite[Lemma 17.6.1]{F}. Direct computation of $\pi_0^\infty(\tau)$ is not hard at all   but here we apply the fact that  $\pi_0^\infty$ is the Fourier transform of $\mathfrak{p}_{c_\circ t}(0) =\mathfrak{p}_{c_\circ}(0)/t^{1/\a} \; (t\geq 0)$ (verified by using the analogue of  (\ref{eq3.0})) so that
$$\pi_0^\infty(\tau)  = \mathfrak{p}_{c_\circ}(0) \int_0^\infty t^{-1/\a} e^{i\tau t}dt. $$
This  integral is written as  $|\tau|^{1/\a}\tau^{-1}\int_0^\infty t^{-{1/\a}} e^{it}dt$ and  is evaluated by applying  (\ref{eq3.000}),  giving the formula of the lemma. \qed
\v2
From the expression of $\pi_0^\infty(\tau)$ given in Lemma \ref{lem3.2} we have
\beqn\label{eq3.4}
\frac{1}{\pi^\infty_0(\tau)} 
 =\frac{c_\circ^{1/\a}\a \pi [\sin (\pi/2\a)- i\cos (\pi/2\a)] } {\Ga(1/\a)\Ga(1-1/\a)\sin [\pi(\a-\ga)/2\a]} \tau^{1-1/\a}  \qquad (\tau>0)
 \eeqn
 with which we compute the integral arising in  (\ref{inv_f}) to show Theorem \ref{thm1}.

\v2
{\bf Proof of Theorem \ref{thm1}.} \, The formula to be shown is 
$
f^0(n) \sim \k_{\a,\ga}c_\circ^{1/\a}/n^{2-1/\a}.
$
In  case $\ga=0$  (i.e., the  walk is centered) this is obtained by Kesten \cite{K} and  the same proof applies. Here we  
proceeds somewhat differntly as follows.  On making trivial decomposition  $1/\pi_0 = 1/\pi_0^\infty + [1/\pi_0-1/\pi_0^\infty]$ an integration by parts transforms $f^0(n) =\frac2{\pi}  \int_0^\pi \Re \,[ -1/\pi_0(\tau)]  \cos n\tau \,d\tau $ into 
$$ \frac2{\pi n} \int_0^\pi \Re \, [1/\pi_0^\infty]'(\tau) \sin n\tau d\tau + \frac2{\pi n}
\int_0^\pi  \Re \,[1/\pi_0 -1/\pi_0^\infty]'(\tau) \sin n\tau d\tau.
$$
The first term, easily evaluated by (\ref{eq3.000}) owing to (\ref{eq3.4}), gives the asymptotic form asserted by the lemma. The second integral restricted on $[0,1/n]$ is shown to be  $o(1/n^{1-1/\a})$ by using
$[1/\pi_0 -1/\pi_0^\infty]'(\tau) = o(|\tau|^{-1/\a})$ and that on $(1/n,\pi)$ is dealt with by integrating by parts once more. Further details are omitted.
 \qed
\begin{lem}\label{lem3.7} \, For any $M>1$, uniformly for $M^{-1} < x/n^{1/\a}< M$, as $n\to\infty$
$$f^x(n) \sim \mathfrak{f}^x(n).$$
\end{lem}
\n\pf\,
Bring  in the functions  $R_1(\tau,\th)$ and $R_2(\tau,\th)$ by
$$
R_1=\frac1{1-e^{i\tau}\phi(\th)}-\frac{1}{-i\tau +1-\phi(\th)}~~~\mbox{and}~~~R_2=\frac{1}{-i\tau +1-\phi(\th)}-\frac{1}{-i\tau + c_\circ\psi(\th)}
$$
so that
\beqn\label{eq8}
\frac1{1-e^{i\tau}\phi(\th)}=\frac{1}{-i\tau +c_\circ\psi(\th)}+R_1+R_2.
\eeqn
It is easily observed that
$$f^x(n) =\frac1{\pi}\int_{-\pi}^\pi \frac{\pi_{-x}^\infty(\tau)}{\pi_0(\tau)} \cos n\tau \,d\tau + \frac{1}{2\pi^2}\int_{-\pi}^\pi \frac{\cos n\tau \,d\tau}{\pi_0(\tau)} \int_{-\pi}^\pi (R_1+R_2)e^{ix\th}d\th.$$
 Using  Lemma \ref{lem3.2} we can readily deduce that the first term on the RHS is 
asymptotically equivalent to $\mathfrak{f}^x(n)$ in the same sense as in the lemma and that the second term is $o(1/n)$, which shows the assertion of the lemma since $\mathfrak{f}^1(t)$ is positive (because of a Huygens-like property) and continuous  on $t>0$ and hence  $\mathfrak{f}^x(n)=  \mathfrak{f}^1(n/x^{1/\a})/x^{1/\a} \geq c_M/n$ for some $c_M>0$ for the range of $x$ specified in the lemma. \qed

\vskip2mm
{\bf 3.3.}\;  In this subsection we prove Theorems \ref{thm2} and \ref{thm3} and Corollary \ref{cor1}. 
Recalling
$ f^\wedge_x(\tau)=[-\1(x=0) +{\pi_{-x}(\tau)}]/{\pi_0(\tau)}$ ($x\neq 0$)   we introduce as in \cite{Uht}
\beqn\label{decm_e}
{\rm e}_x(\tau):=\pi_{-x}(\tau)-\pi_0(\tau)+a(x)
\eeqn
so that 
$$ f^\wedge_x(\tau)=\frac{{\rm e}_x(\tau)}{\pi_0(\tau)}+1-\frac{a^\dagger(x)}{\pi_0(\tau)}.$$
The integral representation  $a(x)=(2\pi)^{-1}\int_{-\pi}^\pi (1-\phi)^{-1}(1-e^{ix\th})d\th$ yields
$${\rm e}_x(\tau)=\frac1{2\pi}\int_{-\pi}^\pi \bigg(\frac1{1-e^{i\tau}\phi(\th)}-\frac1{1-\phi(\th)}\bigg)(e^{ix\th}-1)\,d\th.$$
We make the decomposition ${\rm e}_x(\tau)={\rm c}_x(\tau)/2\pi+i\,{\rm s}_x(\tau)/2\pi$, where
$${\rm c}_x(\tau)=\int_{-\pi}^\pi \bigg(\frac1{1-e^{i\tau}\phi(\th)}-\frac1{1-\phi(\th)}\bigg)(\cos x\th -1)d\th$$
$${\rm s}_x(\tau)=\int_{-\pi}^\pi \bigg(\frac1{1-e^{i\tau}\phi(\th)}-\frac1{1-\phi(\th)}\bigg)\sin x\th\,d\th.$$

The  computations    the present approach necessitates  are carried out in the proofs of the  succeeding two lemmas.
\v2
\begin{lem}\label{lem3.4}~ For some constants $C_1$ and $C_2$
$$\bigg|\int_{0}^\pi \frac{{\rm s}_x(\tau)}{\pi_0(\tau)}\cos n\tau \,d\tau \bigg| \leq \frac{C_1|x|}{n^{1+1/\a}}
\quad \mbox{and}\quad 
\bigg|\int_{0}^\pi \frac{{\rm c}_x(\tau)}{\pi_0(\tau)}\cos n\tau \,d\tau \bigg| \leq \frac{C_2|x|^{1+\e}}{n^{1+(1+\e)/\a}},$$
where in the second bound $\e$ is  any  constant not larger than unity  such that  $0\leq \e < 2\a -2$  and $C_2$ may depend on $\e$.
\end{lem} 
\n
\pf \,  First we claim
\beqn\label{eq3.6}
| (d/d\tau)^j{\rm s}_x(\tau)| \leq C|x| |\tau|^{2/\a -1-j}\qquad (j=0,1,2,3).
\eeqn
Writing
\beqn\label{s1}
{\rm s}_x(\tau)=\int_{-\pi}^\pi \frac{(e^{i\tau}-1)\phi(\th)}{1-e^{i\tau}\phi(\th)}\cdot\frac{\sin x\th }{1-\phi(\th)}d\th
\eeqn
and using   $|1-e^{i\tau}\phi(\th) |\geq C^{-1}(|\tau| + |\th|^\a)$  we see that
\beqn\label{s10}
|{\rm s}_x(\tau)|\leq C'\int_0^\pi \frac{|\tau x|\th^{1-\a}}{|\tau| + |\th|^\a}d\th\leq C'|x||\tau|^{2/\a-1}\int_0^\infty \frac{u^{1-\a}}{1+u^\a}du,
\eeqn
hence  the claimed bound of ${\rm s}_x$.  Differentiating the defining expression of ${\rm s}_x$  we  have
\beqn\label{eq3.7}
{\rm s}'_x(\tau) =ie^{i\tau} \int_{-\pi}^\pi \frac{\phi(\th)\sin x\th}{\{1-e^{i\tau}\phi(\th)\}^{2}}\,d\th,
\eeqn
 which  yields the claimed bound for $j=1$ in the same way as   above.  Those for $j=2,3$ are similar
 and the claim has been  verified.

 Now integrating by parts gives
\beqn\label{eq3.8}
\int_{0}^\pi \frac{{\rm s}_x(\tau)}{\pi_0(\tau)}\cos n\tau \,d\tau = \frac{-1}{n}\int_{0}^\pi  
[{\rm s}_x/\pi_0]'(\tau) \sin n\tau \,d\tau. 
\eeqn
By Lemma \ref{lem3.2}   $|(d/d\tau)^j(1/\pi_0(\tau))|\leq |\tau|^{1-1/\a -j}$, which together with
(\ref{eq3.6}) shows that the integral restricted to $\tau<1/n$ is $O(xn^{-1-1/\a})$.  On integrating by parts once more  the remaining integral  admits the same bound, showing the first one of the lemma.

Following the proof of (\ref{eq3.6})  performed above but by using   the bound $1-\cos x\th
\leq |x\th|^{1+\e}$ in place of $|\sin x\th|\leq |x\th|$ (so that the integral corresponding to the last one in (\ref{s10}) is finite) we obtain
\beqn\label{eq3.81}
| (d/d\tau)^j{\rm c}_x(\tau)| \leq C|x|^{1+\e} |\tau|^{(2+\e)/\a -1-j}\qquad (j=0,1,2,3).
\eeqn
The rest of the proof is the same as above. \qed
\v2

By the same computation as in the preceding proof  we obtain the following bounds
\beqn\label{3.40}
\bigg|\int_{0}^\pi  {\rm s}_x(\tau)\cos n\tau \,d\tau \bigg| \leq \frac{C_1|x_n|}{n^{1/\a}},
\quad \mbox{and} \quad  
\bigg|\int_{0}^\pi {\rm c}_x(\tau)\cos n\tau \,d\tau \bigg| \leq \frac{C_2|x_n|^{1+\e}}{n^{1/\a}};
\eeqn
 expanding ${\rm e}_x(\tau)  {\rm e}_{-y}(\tau) =(-{\rm s}_x{\rm s}_{-y} +i{\rm s}_x{\rm c}_{_y} +i {\rm c}_x{\rm s}_{-y} + {\rm c}_x{\rm c}_{-y})(\tau)/4\pi^2$
we also have   
 \beqn\label{eqL3.4}
 \bigg| \int_{0}^\pi \frac{{\rm e}_x(\tau)  {\rm e}_{-y}(\tau)}{\pi_0(\tau)} \cos n\tau d\tau\bigg| \leq \frac{C_3}{n^{1/\a}} \{|x_n|\vee 
|x_n|^{1+\e}\} \{|y_n|\vee 
|y_n|^{1+\e}\},
 \eeqn
which are  used not in this but in the next section. Here  $\e$ is  chosen as in Lemma \ref{lem3.4}.

\v2
\begin{lem}\label{lem3.5}\,  There exists a   constant $\La$ such that for each $\e>0$  there exists $\de>0$ such that 
$$\bigg|\frac2{\pi}\int_{0}^\pi \Re\, \frac{i{\rm s}_x(\tau)}{\pi_0(\tau)}\cos n\tau \,d\tau 
- \frac{ \La x_n}{n} \bigg| <\e \frac{|x_n|}{n} \quad
\mbox{if}\quad  |x_n|<\de,\; |x|\wedge n>1/\de,
$$
where $x_n=x/n^{1/\a}$.
\end{lem} 


\n
\pf\,  We evaluate the RHS of (\ref{eq3.8}). Take $M>1$ such  that $\cos M=0$. Then, on   integrating by parts and  applying $|({\rm s}_x/\pi_0)''(\tau)| \leq C |x| \tau^{1/\a-2}$
\begin{eqnarray}\label{eq3.9}
 \frac{1}{n}\bigg|\int_{M/n}^\pi \Re\, [i{\rm s}_x/\pi_0]'(\tau) \sin n\tau \,d\tau\bigg| &=&
 \frac{1}{n^2}\bigg|\int_{M/n}^\pi \Re\, [i{\rm s}_x/\pi_0]''(\tau) 
 \cos n\tau \,d\tau\bigg| \nonumber\\
&\leq& C(1-1/\a)^{-1}M^{1/\a -1}|x_n|/n.
\end{eqnarray}

Here we have applied the fact that  $\Re\,  [i{\rm s}_x/\pi_0]'(\tau)$ vanishes at $\pi$ since
it is odd and periodic with period $2\pi$, hence attains the same value for $\tau= \pm \pi$.
 [To see that $\Re\,[ i{\rm s}_{x}/\pi_0]'(\tau) $ is odd, it suffices to show that ${i\rm s}_{x}$  (as well as $\pi_x(\tau)$) has the even real and odd imaginary parts, which may be verified, e.g., by observing that $\int_{-\pi}^\pi (1-re^{i\tau}\phi(\th))^{-1}\sin x\th\,d\th$ is represented by a Fourier series (with real coefficients).]

 For the integral over $0<\tau <M/n$  let $c_\circ =1$ for simplicity.   We   replace 
$1-e^{i\tau}\phi(\th)$ by $-i\tau + \psi(\th)$ and $1-\phi(\th)$ by $ \psi(\th)$ in the  integral defining  ${\rm s}_x(\tau)$    as in the proof of Lemma  \ref{lem3.7},  the replacement  being justified without difficulty in view of  (\ref{eq3.2}).  We  further replace $\sin x\th$ by $x\th$ and  extend the range of integration to the whole real line, which we shall show to cause only a negligible error  (see   the end of this proof).  
 In any case these modifications  of ${\rm s}_x(\tau)$ together  result in the function
$${\rm s}^\circ_x(\tau) := \int_{-\infty}^\infty \bigg(\frac1{-i\tau + \psi(\th)}-\frac1{  \psi(\th)}\bigg) x\th\,d\th = \int_{-\infty}^\infty 
 \frac{i x\tau \th}{\{-i\tau +c_\circ \psi(\th)\} \psi(\th)  } d\th.$$
 In view of (\ref{eq3.9}) it will  suffice to show that for any $\e>0$ there exists $\de>0$ such that for each $M$, $ |x|$ and $n$ large enough, if $ |x_n|<\de$, then 
\beqn\label{eq3.10}
\bigg|\frac2{\pi}\int_{0}^{M/n} \Re\, [-i {\rm s}^\circ_x/\pi_0]'(\tau)\sin n\tau \,d\tau 
-  \La x_n \bigg| < \e |x_n|. 
\eeqn\label{circ}
  After substitution of $\psi(\th) = e^{\pm i\pi\ga/2}|\th|^\a$ and the change of variable  $u = \th/|\tau|^{1/\a}$ we have
\beqn 
-i{\rm s}^\circ_x(\tau) = \frac{x|\tau|^{2/\a}}{\tau}\int_{-\infty}^\infty \frac{u }{(-i e^{\pm i\pi\ga/2}\,\sgn\,\tau + e^{\pm i\pi\ga}|u|^\a)|u|^\a}du,
\eeqn
where the upper or lower sign in $\pm$ prevails according as $u>0$ or $u<0$. By Lemma \ref{lem3.2} or (\ref{eq3.4})
$$- \Re\, [\, i{\rm s}^\circ_x/\pi^\infty_0]'(\tau) =  \La_1 x\tau^{1/\a-1} \quad \mbox{for}\quad \tau >0$$
 for a constant $\La_1$, and with the help of  $ (d/d\tau)^j [{\rm s}^\circ_x(\tau) (1/\pi_0 -1/\pi^\infty_0)(\tau)] =  o(|\tau|^{1/\a-j})$ ($j=1, 2$)
we deduce that 
 $$-\int_{0}^{M/n} \Re\, [i{\rm s}^\circ_x/\pi_0]'(\tau)\sin n\tau \,d\tau \sim  x_n \La_1  \int_0^{M} s^{1/\a-1} \sin s \, d s $$
and  on using  (\ref{eq3.000})  we conclude that (\ref{eq3.10}) holds with $\La= 2\pi^{-1}\La_1\Ga(1/\a)\sin  (\pi/2\a)$. 

It remains  to show that the error caused by the replacement of $s_x$ by $s^\circ_x$  is negligible.  The range  $|\th|>   1/x$ in the integral defining $s_x$, which  corresponds to
 $u> 1/x\tau^\alpha$ in the integral on the RHS of (\ref{circ}) that (absolutely) converges, 
 is negligible since  for $|\tau| < M/n$,  $x\tau^\alpha \to 0$ as  $x_n\to 0$. The same is true for the derivative 
 $$(s_x^\circ)'(\tau)=i\int_{-\infty}^\infty \{-i\tau +\psi(\th)\}^{-2}\sin x\th\,d\th = i\int_{-\infty}^\infty \{-i +\psi(u)\}^{-2}\sin (x\tau^{1/\alpha}u)\,du,$$
  the same integral as obtained by replacing $\{1-e^{i\tau}\phi(\th)\}$ and $e^{i\tau}\phi(\th)$ 
by $\{-i\tau +\psi(\th)\}$ and $1$, respectively,  in the RHS of  (\ref{eq3.7}), so that  $(s_x^\circ)'(\tau) \sim s'_x(\tau)$ uniformly for $|\tau|<M/n$ as  $x_n \to 0$. This finishes the proof of Lemma \ref{lem3.5}.
\qed

\v2
{\bf Proof of Theorem \ref{thm2}.}  Let $c_\circ =1$ for simplicity. According to  the decomposition (\ref{decm_e}) we have 
\beqn\label{pf_T2}
f^x(n) = \frac{1}{\pi^2}\int_0^\pi \Re\,\frac{{\rm c}_x(\tau) +i{\rm s}_x(\tau)}{\pi_0(\tau)}\cos n\tau\, d\tau 
+ a^\dagger(x) \frac{2}{\pi}\int_0^\pi \Re\,\frac1{\pi_0(\tau)}\cos n\tau \,d\tau.
\eeqn
By Lemmas \ref{lem3.4} and \ref{lem3.5} it follows that as $x_n\to 0$
\beqn\label{eq3.19}
f^x(n) =a^\dagger(x)f^0(n) + \frac{\La x_n}{\pi^2 n}\{1 +o(1)\}.
\eeqn
The first term on the RHS is the leading term and we have the first formula of (\ref{eq_thm1}). Indeed this is evident from Theorem \ref{thm1} when  $x$ remains in a bounded set since   $a^\dagger(x)>0$, while  applying  and Lemma \ref{lem3.1}(i)  in addition we have
$a^\dagger(x)f^0(n)  \sim  \kappa_\pm |x_n|^{\alpha-1}/n$ as $x_\pm\wedge n \to\infty$ with some $\kappa_\pm>0$, showing that the second term of (\ref{eq3.19}) 
is negligible as $x_n \to 0$.
   The second formula of  (\ref{eq_thm1})  follows from Lemma \ref{lem3.7}. \qed

\v2
{\bf Proof of Theorem \ref{thm3}.} Let $\ga=2-\a$.  First note that for  the regime  $1/M \leq |x_n| \leq M$ the result follows from Lemma  \ref{lem3.7}. In case $x_n\to 0$ we apply
relation (\ref{eq3.19})   (valid for all $\gamma$). For $x<0$  $a(x)$ behave in a similar way to the case $|\ga| <2-\a$, so that  the preceding proof works well. For $x>0$,  it follows that  $a(x) = o(x^{\a-1})$ as $x\to\infty$, hence, on the  one hand,  taking limit in  (\ref{eq3.19}) we obtain 
$$ nf^x(n)/x_n \to   \La/\pi^2 \qquad\mbox{as $x_n \to \xi >0$ and $\xi\downarrow 0$ in this order}.$$
 On the other  hand,  owing to the identity $c_\circ \mathfrak{f}^x(c_\circ n) = x_n\mathfrak{p}_{c_\circ}(-x_n)/n$  it follows that
  $$n c_\circ \mathfrak{f}^x(c_\circ n)/x_n \to  \mathfrak{p}_{c_\circ}(0)$$ in the same way of taking the limit  as  above.  By the result for the case $x_n \asymp 1$ this leads to $\La/\pi^2 = \mathfrak{p}_{c_\circ}(0)$, which allows us to replace the second term on the RHS of (\ref{eq3.19}) by $x_n\mathfrak{p}_{c_\circ}(0)/ n$, thus
  concludes the proof, the case $\ga =-2+\a$ being dealt with in the same way. \qed
 
  \begin{rem}\label{rem5}
  In view of (\ref{pf_T2})---recall  ${\rm c}_x +i{\rm s}_x =2\pi {\rm e}_x$---what is shown in the proofs above is paraphrased as  follows: If $\ga=2-\a$, then uniformly for $|x_n|<M$, as  $n\to\infty$
 \beqn\label{pf_T3}
 \frac{1}{\pi}\int_{-\pi}^\pi \frac{{\rm e}_{x}(\tau)}{\pi_0(\tau)}\cos n\tau\, d\tau 
  = \left\{\begin{array}{lr} x_n\mathfrak{p}_{c_\circ}(-x_n)/{n} + o(a(x)/n^{2-1/\a})
 &(x>0),\\
  o(a(x)/n^{2-1/\a}) &(x<0,  x_n \to 0),
\end{array}\right.
  \eeqn
and if $|\ga|<2-\a$,  the integral on the LHS is $o(a(x)/n^{2-1/\a})$ as $x_n\to0$.
\end{rem}

\v2
{\bf Proof of Corollary \ref{cor1}.}   The first expression of  $\k^{\mathfrak{f}}_{\a,\ga}$ as well as the equivalence relation in case $\ga\neq 2-\a$ follows from Lemma \ref{lem7.1}. For $\ga=2-\a$ the equivalence relation follows from what is mentioned in the paragraph  preceding  the corollary. As in the last part of the  proof of Theorem \ref{thm3} given above,    by Lemma \ref{lem3.7} (with $c_\circ =1$) and scaling relation of $\mathfrak{f}^x(t)$  it follows that 
$$nf^x(n) \sim  n\mathfrak{f}^x(n) \sim \mathfrak{f}^1(1/x_n^\a)/x_n^\a, $$
which together with Theorems \ref{thm2} and \ref{thm3} shows that if $\ga <2-\a$, then
 $\mathfrak{f}^1(t)\sim \k^{\mathfrak{f}}_{\a,\ga}/t^{2-1/\a}$ with $\k^{\mathfrak{f}}_{\a,\ga}$ determined by  $\k_{\a,\ga}a(x)/n^{1-1/\a} \sim  \k^{\mathfrak{f}}_{\a,\ga} x_n^{\a-1}.$
 By Lemma  \ref{lem3.1}(i) and the expression defining $\k_{\a,\ga}$, this leads to the second expression of $\k^{\mathfrak{f}}_{\a,\ga}$. \qed

\v2
Because of the similarity of the proof to that of Lemma \ref{lem3.5} we here give the following lemma that is used  in  the next section. 

\begin{lem}\label{lem3.6} \, There exists a   function  $D_n(y), y\in \Z$ such that  $|D_n(y)|\leq C |y_n|$ for $|y_n|<M$ and  for each $\e>0$  there exists $\de>0$ such that 
if $  |x_n|<\de,\; |x|\wedge n>1/\de, \, |y_n|<1/\e$,
$$\bigg|\frac2{\pi} \int_{0}^\pi \Re\, \frac{{\rm e}_x(\tau)  {\rm e}_{-y}(\tau)}{\pi_0(\tau)} \cos n\tau \, d\tau 
- \frac{ D_n(y) x_n}{n^{1/\a}} \bigg| <\e\frac{|x_ny_n|}{n^{1/\a}}.$$
\end{lem} 
\n
\pf\, On recalling the derivation of (\ref{eqL3.4})  the terms $|x|$ and  $|x|^{1+\e}$ on the RHS of it correspond to  ${\rm s}_x$ and ${\rm c}_x$, respectively and similarly for $|y|$ and $|y|^{1+\e}$, and one sees it suffices to show that
\beqn\label{eqL3.6}
\bigg| \frac2{\pi}\int_{0}^\pi \Re\, \frac{i {\rm s}_{x}(\tau){\rm e}_{-y}(\tau)}{\pi_0(\tau)} \cos n\tau \,d\tau 
-  \frac{ D_n(y) x_n}{n^{1/\a}} \bigg| <\e\frac{|x_ny_n|}{n^{1/\a}}
\eeqn
provided $|x_n|<\de,\; |x|\wedge n>1/\de, \, |y_n|<1/\e$. 
First suppose  $3/2 < \a<2$ so that $-1< 3/\a -2 < 0$. By (\ref{eq3.3}), (\ref{eq3.6}), (\ref{eq3.81}) it follows that $|[s_x {\rm e}_{-y}/\pi_0]'(\tau)| \leq C |xy||\tau|^{3/\a-2}$, and on 
integrating by parts 
\beqn\label{eq3.80}
\frac{2}{\pi}\int_{0}^\pi  \Re\,\frac{i{\rm s}_{x}(\tau){\rm e}_{-y}(\tau)}{\pi_0(\tau)}
\cos n\tau \,d\tau 
 =  \frac{2}{\pi n}\int_{0}^{\pi} \Re\, \bigg[\frac{ {-i\rm s}_{x}{\rm e}_{-y}}{\pi_0}\bigg]' (\tau)\sin n\tau \,d\tau.
 \eeqn
Integrating by parts once more we observe  that the 
contribution from $|\tau|>M/n$ to the integral on the RHS   becomes negligibly small as $M$ is taken  large and then  ${\rm s}_{x}(\tau)$ may be replaced by ${\rm s}^\circ_{x}(\tau)$ as in the proof of Lemma \ref{lem3.5}. [Here we have applied the fact that  $\Re\,[ -i{\rm s}_{x}{\rm e}_{-y}/\pi_0]'(\tau) $ is odd, hence vanishes at $\tau =\pi$.]
 Now define   
$$\omega(\tau) = - i{\rm s}^\circ_{x}(\tau)/x 
$$
 and 
 \beqn\label{D_n}
 D_n(y) := \frac{ 2n^{2/\a-1}}{\pi}\int_0^{\pi} \Re\,\bigg[\frac{ \omega{\rm e}_{-y}}{\pi_0}\bigg]' (\tau)\sin n\tau \,d\tau.
 \eeqn
 By the same reason as above  the contribution from $\tau>M/n$ to the integral on the RHS becomes negligible as $M$ gets large, and we see that (\ref{eqL3.6}) is satisfied.
Noting   $|[\omega {\rm e}_{-y}/\pi_0]'(\tau)| \leq C_1 |y||\tau|^{3/\a-2}$ we also deduce that   $|D_n(y)| = O(y_n)$.

 In case $1<\a < 3/2$ we can further  integrate the RHS of (\ref{eq3.80})    by parts  to  have
\beqn\label{eqL3.61}
\int_{0}^\pi  \Re\, \frac{i {\rm s}_{x}(\tau){\rm e}_{-y}(\tau)}{\pi_0(\tau)} \cos n\tau \,d\tau 
 = \frac{1}{ n^2}\int_{0}^{\pi} \Re\, \bigg[\frac{-i {\rm s}_{x}{\rm e}_{-y}}{\pi_0}\bigg]''(\tau) \cos n\tau \,d\tau
 \eeqn
 and accordingly putting
$
 D_n(y) = \frac{2n^{2/\a-2}}{\pi }\int_0^{\pi} \Re\,\Big[\frac{ \omega{\rm e}_{-y}}{\pi_0}\Big]'' (\tau)\cos n\tau \,d\tau
$
and making a  similar argument  to the above we obtain  (\ref{eqL3.6}). 
 
 Let  $\a =3/2$. This is a critical case when $[{\rm s}_{x} {\rm e}_{-y}/\pi_0]'(\tau)$ tends to a constant multiple of $xy$ as $\tau \to 0$, and by the very this fact  we have $[{\rm s}_{x}  {\rm e}_{-y}/\pi_0]''(\tau) = o(1/\tau)$.  Split the range of integral on the LHS of (\ref{eq3.80}) at $\tau =2\pi N/n$ with  a positive integer $N$ and  denote by $I$ and $I\!I$  the integrals over  $(0,M/n]$ and $[M/n,\pi]$, respectively,  where $M=2\pi N$. Then on  integrating by parts  once more 
 \[
 I =
  \frac1{n^2}\int_{0}^{M/n} \Re\, \bigg[\frac{ -i{\rm s}_{x}{\rm e}_{-y}}{\pi_0}\bigg]''(\tau) (\cos n\tau \, -1) \,d\tau, \quad \mbox{and}
\]
 $$
I\!I = \frac1{n^2}\Re\,\bigg[\frac{ -i{\rm s}_{x}{\rm e}_{-y}}{\pi_0}\bigg]'(M/n)
+ \frac1{n^2}\int_{M/n}^\pi \Re\,\bigg[\frac{ -i{\rm s}_{x}{\rm e}_{-y}}{\pi_0}\bigg]''(\tau) \cos n\tau \,d\tau. 
 $$
Since the integrand of the first integral above is at most $o(1/\tau)\times n\tau$,   we see that $I= o(1/n^2)$.   The second integral  which we further integrate by parts  is dominated by a constant multiple of  $|xy|/Mn$, thus negligible since  $M$ can be chosen arbitrarily large, while 
$\Re\, [-i{\rm s}_x{\rm e}_{-y}/\pi_0](M/n) \to \k xy$ ($n\to\infty$) with some $\kappa \in \R$.  Finally recalling  $n^2 = n^{3/\alpha}$,
we find  that   (\ref{eqL3.6}) holds with $D_n(y)= (2/\pi) \kappa y_n$, and hence conclude the proof of the lemma. 
   \qed

\section{Estimates of $p_{\{0\}}^n(x,y)$}

In this section we prove Theorems \ref{thm4} and \ref{thm5}, the proofs being given at the end of the section.   
We continue to use the notation $\pi_x(\tau)$ introduced in the preceding section. 

The  arguments that follow  are based on  the representation
\beqn\label{eq3.1}
p_{\{0\}}^n(x,y)=p^n(y-x)-\sum_{k=1}^n f^x(n-k)p^k(y)
\eeqn
or,  to say more exactly,   its Fourier version: from  (\ref{eq3.0}) one can easily deduce that $p^n(x) = (1/2\pi)\int_{-\pi}^\pi \pi_x(\tau) e^{-i n\tau} \,d\tau$ and, on combining this with (\ref{hat_f}),  (\ref{eq3.1}) may be written as
\beqn\label{eq4.2}
p_{\{0\}}^n(x,y)=\frac1{2\pi}\int_{-\pi}^\pi \bigg[\pi_{y-x}(\tau)- \frac{\pi_{-x}(\tau)\pi_y(\tau)}{\pi_0(\tau)}\bigg]e^{-in\tau}d\tau~~~~~(x\neq 0)
\eeqn
and
\beqn\label{eq3.20}
p_{\{0\}}^n(0,y)=\frac1{2\pi}\int_{-\pi}^\pi \frac{\pi_y(\tau)}{\pi_0(\tau)}e^{-in\tau}d\tau.
\eeqn
Note that for $y\neq 0$, $p_{\{0\}}^n(0,y)=f^{-y}(n)$ by duality (or by coincidence of the Fourier coefficients),  so that  in the case $x=0$ the required estimate  is immediate from 
Theorems \ref{thm1} to \ref{thm3} that have been verified in the preceding section.
 
 \begin{lem}\label{lem4.1}~  Uniformly for $x, y\in \Z$, as $n\to\infty$
\begin{eqnarray*}
{\rm (i)} \; &&p^n(y-x)-p^n(-x)-p^n(y)+p^n(0) \\
&& \quad =\mathfrak{p}_{c_\circ n}(y-x)-\mathfrak{p}_{c_\circ n}(-x)-\mathfrak{p}_{c_\circ n}(y)+\mathfrak{p}_{c_\circ n}(0)+o({xy}/n^{3/\a}); \mbox{and}\qquad \qquad\qquad\\[3mm]
{(\rm ii)}\; &&   p^n(y-x)-p^n(-x)
=\mathfrak{p}_{c_\circ n}(y-x)-\mathfrak{p}_{c_\circ n}(-x)+ o(y/n^{2/\a}).\qquad \qquad
\end{eqnarray*}
\end{lem}
\n
\pf\, Put $c= c_\circ \cos (\pi\ga/2) (>0)$, choose a positive constant $\e$ so that $1-|\phi(\th)| \geq  |\th|^\a c/2$ for $|\th|<\e $ and put $\eta=\sup_{\e\leq |\th|\leq \pi}|\phi(\th)| (<1)$.  Then the error in the first relation  (i)  is written as
\beq
 \,\frac1{2\pi} \int_{-\e}^\e \Big([\phi(\th)]^n-e^{-nc_\circ\psi(\th)}\Big)K_{x,y} (\th)d\th+O(\eta^n \vee e^{-nc\e^\a}),
\eeq
where  and 
$K_{x,y} (\th)=e^{-i(y-x)\th}-e^{ix\th}-e^{-iy\th}+1.$ By (\ref{f_hyp}) $\log[\phi(\th)e^{c_\circ \psi(\th)}]= o(|\th|^\alpha)$ as $\th\to 0$. Since $K_{x,y} (\th)=(e^{ix\th}-1)(e^{-iy\th}-1)$, we have $|K_{x,y} (\th)|\leq |xy|\th^2$ and, scaling $\th$ by $n^{1/\a}$ and applying the dominated convergence theorem,  we deduce that the integral above is $o(xy/n^{3/\a})$, showing  (i). 

The proof of (ii) is similar, rather simpler.  One may only to use  $|e^{-i(y-x)\th}-e^{ix\th}|\leq |y\th|$ in place of the bound of $K_{x,y}(\th)$.  \qed

\begin{lem}\label{lem4.2} \,  Given $M>1$, if either  $\ga=2-\a$ and  $y > 0$ or $|\ga|<2-\a$, then as $n\to\infty$ and $|x_n|\vee y_n \to 0$ 
\beqn\label{eq_lem4.2}
p_{\{0\}}^n(x,y)\sim a(-y) f^x(n).  
\eeqn
If either  $\gamma =2-\alpha$ and $x_+/a^\dagger(x) = o(n^{2/\a-1})$   or $|\ga|<2-\a$,  then the expression on the RHS is asymptotically equivalent to $\k_{\a,\ga}c_\circ^{1/\a} a^\dagger(x)a(-y)/n^{2-1/\a}$ as $n\to\infty$ and $|x_n| \to 0$. 
\end{lem}
\n
\pf \, 
Of the integrand in
 (\ref{eq4.2}) we make  the decomposition
\begin{eqnarray}\label{eqL4.2}
\pi_{y-x}(\tau)- \pi_{-x}(\tau)\pi_y(\tau)/{\pi_0(\tau)}&=&\pi_{y-x}-\pi_{-x}-\pi_{y}+\pi_{0}
- a(x)a(-y)/\pi_0 \nonumber\\
&&+[ -{\rm e}_{x}\,{\rm e}_{-y}+a(x)\,{\rm e}_{-y}+a(-y)\,{\rm e}_{x}]/\pi_0,
\end{eqnarray}
where we recall
$\,{\rm e}_{x}=\,{\rm e}_x(\tau)=\pi_{-x}(\tau)-\pi_0(\tau)+a(x) =[\,{\rm c}_x(\tau)+i\,{\rm s}_x(\tau)]/2\pi.$
Noting that 
$$\mathfrak{p}_t(y-x)-\mathfrak{p}_t(-x)- \mathfrak{p}_t(y)+\mathfrak{p}_t(0)= -xy \mathfrak{p}_{1}''(0)t^{-3/\a}\{1+ o(1)\} \quad \mbox{ as} \quad |x_n|\vee |y_n| \to 0,$$
 we apply  Theorem \ref{thm1} and Lemma \ref{lem4.1}(i) to see
\begin{eqnarray}
&&\frac1{2\pi}\int_{-\pi}^\pi \bigg[\pi_{y-x}-\pi_{-x}-\pi_{y}+\pi_{0}
- \frac{a(x)a(-y)}{\pi_0}\bigg]e^{-in\tau}d\tau \nonumber\\
&&=p^n(y-x)-p^n(-x)-p^n(y)+p^n(0)+a(x)a(-y)f^0(n)   \nonumber\\
&&= \k_{\a,\ga}\frac{c^2_\circ a(x)a(-y)}{(c_\circ n)^{2-1/\a}} - \frac{\mathfrak{p}_{1}''(0) xy}{(c_\circ n)^{3/\a}} +o\bigg(\frac{xy}{n^{3/\a}}\bigg).
\label{eq3.21}
\end{eqnarray}
By Lemma \ref{lem3.4}  and   (\ref{eqL3.4}) 
the integral  $\int_{-\pi}^\pi[-\,{\rm e}_{x}\,{\rm e}_{-y}+a(x)\,{\rm e}_{-y}+a(-y)\,{\rm e}_{x}] \cos n\tau \,d\tau/\pi_0$  is 
dominated in absolute value by a constant multiple of 
\beqn\label{eq3.211}
\frac{|x_ny_n|}{n^{1/\a}} + \frac{a(x)|y_n| +a(-y)|x_n|}{n}     \quad\mbox{for}\quad  |x_n|\vee |y_n|\leq M.
\eeqn

 If  either  $\ga=2-\a$, $x<M $ and  $y>0$  or $\ga<2-\a$, both ratios in (\ref{eq3.211})  are   $ o(a(x)a(-y)/n^{2-1/\a})$ (as $|x_n|\vee y_n\to0$),  so that  the first term on the right most member of (\ref{eq3.21}) is dominant over the others and in view of   Theorem \ref{thm2}   formula (\ref{eq_lem4.2}) follows. 

In the other case   $\ga=2-\a$, $x \geq M$ and  $y>0$, the term $a(x)|y_n|/n$ in (\ref{eq3.211})   is negligible while we have to take account of $a(-y)|x_n|/n$  and turn back to the integral $\pi^{-1}\int_{-\pi}^\pi a(-y)\,{\rm e}_{x}(\tau) \cos n\tau \,d\tau/\pi_0(\tau)$ which is asymptotic to  
$ a(-y)\mathfrak{p}_{c_\circ}(0)x_n/n$
as we have noted   in Remark \ref{rem5} (after the proof of Theorem \ref{thm3}) so that  the terms of order at most $|x_ny_n|/n^{1/\alpha}$ are negligible, and we see that 
the combination 
 $$ \k_{\a,\ga}\frac{c^2_\circ a(x)a(-y)}{(c_\circ n)^{2-1/\a}} + \frac {a(-y)\mathfrak{p}_{c_\circ}(0)x_n}{n} \sim   a(-y)f^x(n)$$
constitutes the leading term.  See (\ref{crs_ov}) for the second half of the theorem. \qed

\v2
 From (\ref{eq3.21}) and (\ref{eq3.211}) (with a simple amplification for the  case $|x_n|\wedge |y_n|> 1/M$) we have the following upper bound:  For some constant $C$ depending only on  $F$,
$$p_{\{0\}}^n(x,y) \leq C\bigg[\frac{ a^\dagger(x)a(-y)}{ n^{2-1/\a}}+\frac{a^\dagger (x)|y| +a(-y)|x|}{n^{1+1/\a}}\bigg] \quad \mbox{if}\quad  |x|\vee |y|< n^{1/\a}. $$
In the next section we shall remove the  restriction $|x|\vee |y| < n^{1/\a}$  and improve the estimate in  case $|\ga|=2-\a$ and $\ga x< 0, xy<0$ (cf. Proposition \ref{prop5.2}).
\begin{lem}\label{lem4.3} \,Uniformly for $ |x_n|, |y|_n \in  [1/M, M]$, as $n\to\infty$ 
$$p_{\{0\}}^n(x, y) \sim \mathfrak{p}^{\{0\}}_{c_\circ n}(x,y) =\mathfrak{p}^{\{0\}}_{c_\circ}(x_n,y_n)/n^{1/\a} .$$
\end{lem}
\n\pf\, Let $c_\circ=1$ for simplicity.
In view of identity (\ref{eq3.1}) It suffices to show that for $\e>0$, 
$$\sum_{\e n\leq k\leq (1-\e)n}f^x(n-k)p^k(y) \sim \frac1{n^{1/\a}}\int_\e^{(1-\e)} \mathfrak{f}^{x_n}(t-s)\mathfrak{p}_s(y_n)ds\quad (n\to\infty)$$
and  the sum over $k\in [0, \e n]\cup [(1-\e)n,n] $  and the corresponding integral are both negligible as $n\to\infty$ and $\e\downarrow 0$ in this order. The first requirement is easily deduced  from the asymptotic form of $f^x(k)$ given in Theorems \ref{thm2} and \ref{thm3} (see also (\ref{eqR(b)}) in case $\ga=2-\a$) and the local limit theorem \cite{GK}, according to which  uniformly for $y\in \Z$, as $k\to\infty$ 
\beqn\label{l_l_thm}
p^k(y) = \mathfrak{p}_{k}(y) + o(1/k^{1/\a}).
\eeqn

 For the second one, we address only the sum, the integral being similarly treated.  Denoting  the sums over $k\in [0, \e n]$ and $[(1-\e)n,n]$ by $\Sigma_{<\e}$ and $\Sigma_{>(1-\e)}$, respectively,  we must show that 
$$\lim_{\e\downarrow 0}\lim_{n\to\infty} (\Sigma_{<\e} +\Sigma_{>(1-\e)})n^{1/\a}= 0.$$
The sum  $\Sigma_{>(1-\e)}$ is immediately disposed  of by  the fact that  $p^{k}(y) =O(n^{-1/\a})$ for $k>(1-\e)n$ and   $\sum_{k>(1-\e)n} f^x(n-k) = P[\sigma^x_{\{0\}} \leq \e n] \to 0$ 
  in the present scheme of passing to the limit.  As for $\Sigma_{<\e}$ we use the bound   $f^x(n-k) \leq O(1/n)$ ($k <\e n$) as well as  (\ref{l_l_thm}) to see that 
$$ \Sigma_{<\e} \leq \frac{C}{n}\sum_{k < \e n}\mathfrak{p}_k(y) = \frac{C}{ny}\sum_{k < \e n} \mathfrak{p}_{k/y^\a}(1)\sim \frac{Cy_n^{\a-1}}{n^{1/\a}}\int_0^{\e/y_n^{\a}}\mathfrak{p}_t(1)dt, $$
showing $n^{1/\a} \Sigma_{<\e}  \to 0$ as required. \qed
\v2

 Theorems \ref{thm4} and  \ref{thm5} follow from 
 Lemmas \ref{lem4.2} and \ref{lem4.3} when either  $|x_n|\vee |y_n| \to 0$ or $|x_n|\wedge |y_n|$ is bounded away from zero.  We need to deal with the case
when    $|x_n|\vee |y_n|$  is  bounded away from zero and  $|x_n|\wedge |y_n| \to 0$. 
\begin{lem}\label{lem4.4} \, For any $M>1$, uniformly for   $1/M<|x_n|\vee |y_n| <M $, it holds that if   $\ga=2-\a$, then  as $ n\to\infty$ and $|x_n|\wedge y_n \to 0$ under $y>0$
\beqn\label{eq_lem4.4} 
p_{\{0\}}^n(x,y) \sim  \left\{ \begin{array} {ll} a(-y) f^x(n) \quad &  y_n \to 0,\\
{\displaystyle   a^\dagger(x)f^{-y}(n)+ \frac{x_+ K_{c_\circ}(y_n)}{n^{2/\a} } } 
\quad & x_n \to 0, \end{array} \right.
\eeqn
where  $K_{t}(\eta)$  is given by (\ref{Kt}); 
and   if $|\ga|<2-\a$, then  $p_{\{0\}}^n(x,y) \sim  a(-y) f^x(n)$ or  $a^\dagger(x) f^{-y}(n)$  according as $y_n \to 0$ or  $ x_n \to 0, y\neq 0$.
\end{lem}
\n
\pf \, 
 As before the proof rests on the Fourier representation (\ref{eq4.2}).  Let $\gamma >-2+\alpha$.

  First consider the case $y_n\to 0$.  This time we employ the decomposition
$$
\pi_{y-x}- \pi_{-x}\pi_y/\pi_0=\pi_{y-x}-\pi_{-x}+a(-y)\,{\pi}_{-x}/\pi_0+\,{\pi}_{-x}\,{\rm e}_{-y}/\pi_0.
$$
Owing to Lemma \ref{lem4.1}(ii)  and the present assumption on $x, y$, 
$$p^n(y-x)-p^n(-x)= n^{-1/\a}[\mathfrak{p}_{c_\circ}'(-x_n)y_n +  o(y_n)].$$
Hence,  by (\ref{hat_f})
$$
\frac1{\pi}\int_{-\pi}^\pi \bigg[\pi_{y-x}-\pi_{-x}
+\frac{a(-y)\,{\pi}_{-x}}{\pi_0}\bigg](\tau) \cos n\tau \,d\tau
=a(-y)f^x(n) + O(y_n/n^{1/\a}).
$$
Now suppose $1/M\leq |x_n|\leq M$ and  $ y> 0$, which  imply   $f^x(n) \asymp 1/n$ and  $a(-y) \asymp |y|^{\a-1}$, respectively.  Using   Lemma \ref{lem3.4},  (\ref{3.40}) (both with $y$ in place of $x$),  (\ref{eqL3.4}), and  the identity $\pi{_{-x}} = {\rm e}_x  + \pi_0 - a(x)$  we then  deduce 
$$\int_{0}^\pi \frac{\,{\pi}_{-x}(\tau)\,{\rm e}_{-y}(\tau)}{\pi_0(\tau)} \cos n\tau \, d\tau =O(y_n/n^{1/\a}).$$ 
By  Theorems \ref{thm2} and \ref{thm3}, $y_n/n^{1/\a}$  is negligible in comparison to $a(-y)f^x(n) \asymp a(-y)/n$ 
as $y_n\to 0$,   hence the first relation in (\ref{eq_lem4.4})  follows.  

 If $|\ga|<2-\a$,  the first case of it is proved by the arguments above which are valid without the condition  $y>-M$, while  the second case  follows from the first by duality.

 Let $\ga =2-\a$, $1/M\leq y_n \leq M$ and $x\neq 0$.  We follow  the proof of Lemma \ref{lem4.2}, employing  the decomposition  (\ref{eqL4.2}) and  applying the estimates given there.
 On the one hand   by the first equality of (\ref{eq3.21}) 
\beqn\label{eq3.22}
\frac2{\pi}\int_{0}^\pi \Big[\pi_{y-x}-\pi_{-x}-\pi_{y}+\pi_{0}
\Big] \cos n\tau \,d\tau =  - [\mathfrak{p}_{c_\circ}'(y_n) -\mathfrak{p}_{c_\circ}'(0)]x_n n^{-1/\a}\{1+o(1)\}  
\eeqn 
  as $x_n \to 0$  (uniformly for $y_n \in [1/M, M]$) and by 
  $  [{\rm e}_{-y}(\tau)- a(-y)]/\pi_0(\tau) = f_{-y}^{\wedge}(\tau)-1$, 
$$
\frac1{\pi}\int_{-\pi}^\pi \bigg[\frac{a(x){\rm e}_{-y}(\tau)}{\pi_0(\tau)}
+\frac{a(x)a(-y)}{\pi_0(\tau)}\bigg]\cos n\tau \,d\tau = a(x)f^{-y}(n)
$$  
On the other hand  by Lemmas \ref{lem3.4} and \ref{lem3.5}     
 $$\frac2{\pi}\int_{0}^\pi \Re\, \frac{a(-y)\,{\rm e}_{x}(\tau) }{\pi_0(\tau)} \cos n\tau \,d\tau  = \frac{a(-y)x_n \{\Lambda +o(1)\}}{n},$$ 
which together with Lemma \ref{lem3.6}  shows  that  uniformly for $y_n \in [1/M,M]$, as $x_n \to \xi>0$ and $\xi\to 0$ in this order
$$\frac2{\pi}\int_{0}^\pi\Re\,\frac{-\,({\rm e}_{x}\,{\rm e}_{-y})(\tau)+a(-y)\,{\rm e}_{x}(\tau)}{\pi_0(\tau)} \cos n\tau \,d\tau = \frac{[ - D_n(y)+ \Lambda' y_n^{\a-1} ]x_n + o(x_n)}{n^{1/\a}}
$$
with $\Lambda'=\Lambda/c_\circ \Ga(\alpha)$ ($\Lambda$ and $D_n(y)$ are given in the proof of Lemmas \ref{lem3.5} and \ref{lem3.6}), respectively). These together yield
\begin{eqnarray}
n^{1/\a} p^n_{\{0\}}(x,y) = a(x)f^{-y}(n)
 + \{\mathfrak{p}_{c_\circ}'(0) -\mathfrak{p}_{c_\circ}'(y_n) 
 - D_n(y) +  \Lambda y_n^{\a-1} \}x_n + o(x_n).  
\label{4.11}
\end{eqnarray}
Since  $a(x)=o(x^{\a-1})$, $f^0(n)=O(1/n^{2-1/\a})$   and  $R_n(y) =O(y_n/n)$,
 the second term on the RHS of (\ref{4.11}) tends to zero as $x_n\to \xi>0$, and hence
in  view of Lemma \ref{lem4.3} letting  $x_n \to \xi>0$ yield
     $$\mathfrak{p}_{c_\circ}^{\{0\}}(\xi,y_n) 
   =\{\mathfrak{p}_{c_\circ}'(0) -\mathfrak{p}_{c_\circ}'(y_n)  
  - D_n(y)  + \Lambda' y_n^{\a-1}\} \xi   + o(\xi) \quad (\xi\to0). $$
Thus  dividing both sides by $\xi$  and passing to the limit we find  
 $$   K_{c_\circ}(y_n)  =   \mathfrak{p}_{c_\circ}'(0) -\mathfrak{p}_{c_\circ}'(y_n) 
 - D_n(y)  +  \Lambda' y_n^{\a-1} + o(1)$$
($n\to\infty$) uniformly for $y_n\in [1/M,M]$,     which  together with (\ref{4.11}) shows the second relation of (\ref{eq_lem4.4}),  the term 
$O(x_n/n^{1/n})$ being negligible as compared with  $a^\dagger(x)f^{-y}(n)$ for $x<0$.     Lemma \ref{lem4.4} has been proved. \qed

\v2

   

 
 \v2
 {\bf Proof of Theorems \ref{thm4} and \ref{thm5}.} 
Both Theorems \ref{thm4} and \ref{thm5} follow from Lemmas \ref{lem4.2}, \ref{lem4.3} and \ref{lem4.4}. Note that  the case $x=0$ is dealt with by  (\ref{eq3.20}). (Note that $|x_n| <\!< a(x)/n^{1-1/\alpha}$ as $x_n\uparrow 0$.) \qed

 

\section{Esitmates of $p^n_{\{0\}}(x,y)$ in case $xy< 0$ and proof of Theorem \ref{thm6}}

Here we derive estimates of $p^n_{\{0\}}(x,y)$ for $x, y$ not necessarily confined in $0< |x_n|, |y_n|<M$, that lead to Proposition \ref{prop2.3} and are useful  for the proof of   Theorem \ref{thm6}. We  assume $\ga =2-\a$ throughout this section except for Lemma \ref{lem5.2}, the case $xy<0$ for $|\ga|<2-\a$ being included in Theorem \ref{thm4}.     Sometimes we suppose $E|\hat Z|<\infty$ (see (\ref{Z})), which entails $\ga=2-\a$.
[In case  $E|\hat Z|=\infty$ and  $\ga=2-\a$ there arises  a troublesome question caused by the obscure nature of $L(x):= U_{{\rm ds}}(x)/x$ (cf. \cite{VW}; see also Remark \ref{Rprop1}(a)).]

\begin{Prop}\label{prop5.1}   Suppose $E|\hat Z|<\infty$. Then, given  $M\geq 1$,  for  $-M  <y_n <0< x_n<M$ 
\beqn\label{eqP1}
p^n_{\{0\}}(x,y) \geq c_M \bigg( \sum_{w=2}^{x\wedge|y|} p(-w)w^{2\a-1}\bigg)
\bigg[  \frac{a(x)+ a(-y)}{n^{2-1/\a}} + \frac{x_n + |y_n|}{ n} \bigg],
\eeqn
where $c_M$ is a positive constant (depending  on $M$ and $F$).
\end{Prop}
\n\pf\,  The walk is supposed to be not left-continuous, otherwise the result being trivial.
  This proof employs   the obvious lower bound 
   $$p^n_{\{0\}}(x,y) \geq  \sum_{\de x^\a\leq k \leq n/2}\; \sum_{0\leq w \leq (k/\de)^{1/\a}}\;\sum_{z=1}^{x}p^{k}_{(-\infty,0]}(x,w)p(-z-w) p^{n-k}_{\{0\}}(z,y)
$$
valid for any constant $\de>0$. We may  and do suppose  $x\leq -y <Mn^{1/\a}$, the case $-y<x$ being dealt with by duality.  $\de$  needs  to be  chosen so small  that $\de x^\a< \eta n$ for some $\eta<1/2$.  To this end we  take, e.g.,
$$\de = 1/2(2 M^2)^{\a/2},$$
entailing $n/2[(x-y)x]^{\a/2}> \de $, which  after substituting from   $x-y \geq 2x$ and multiplying   by $x^\a$ reduces to  $\de x^\a <  n/2^{1+\a/2}$.

 For $k, w, z$ taken from the range of summation above, we have  by Theorem \ref{thm5} (see (\ref{R3f2}))
$$p^{n-k}_{\{0\}}(-z,y) \asymp a(-z)f^{-y}(n-k)\asymp  z^{\a-1}\{a(-y)+ |y_n|n^{1- 1/\a}\}/n^{2-1/\a},$$
and by (\ref{Doney}) 
$$p^{k}_{(-\infty,0]}(x,w) \asymp  U_{{\rm ds}}(x)V_{{\rm as}}(w) /k^{1+1/\a};$$
 by  $k\geq \delta x^\a$ it also follows that    $x\leq (k/\delta)^{1/\a}$. Hence,    putting
\beqn\label{s_m}
m(x) =  \sum_{z=0}^x\;\sum_{w=0}^x p(-z-w)V_{{\rm as}}(w)z^{\a-1}
\eeqn 
we have
$$p^n_{\{0\}}(x,y) \geq c' m(x)U_{{\rm ds}}(x)\bigg\{\frac{a(-y)}{n^{2-1/\alpha}}  +\frac{ |y_n|}{n}\bigg\}\sum_{\de x^\a\leq k \leq n/2}\;\frac1{ k^{1+1/\a}}.$$
Since $\de x^\a \leq n/2^{1+1/\a}$, the last sum is bounded below by  a positive multiple of $1/x$. 
 In the double sum in (\ref{s_m}) restricting the inner summation to $w\leq x-z$, making change of the variable $w=j-z$ and  interchanging the order of summation we obtain 
 \beqn\label{m_lb}
 m(x)\geq  
 \sum_{j=0}^x p(-j)  \sum_{z=0}^{j} V_{{\rm as}}(j-z)z^{\a-1}.
 \eeqn
Now suppose  $E|\hat Z|<\infty$. Then   $V_{{\rm as}}(w)\asymp w^{\a-1}$ and $U_{{\rm ds}}(x)\asymp x$, and  we see that $m(x)\geq c'' \sum_{j=0}p(-j) j^{2\alpha-1}$ and  the required
lower bound follows. \qed

\v2
\begin{rem}\label{Rprop1}\,
(a) \, Even in case $E|\hat Z|=\infty$ we know that  if $\ga=2-\a$,  $V_{{\rm as}}(w)$ varies regularly of index $\a-1$ and $V_{{\rm as}}(x)U_{{\rm ds}}(x) \sim Cx^\a$  as $x\to\infty$ with a positive constant $C$ \cite{R}; hence
from (\ref{m_lb}) we have $m(x) \geq c_1 \sum_{w=1}^x p(-w)V_{{\rm as}}(w)w^{\a}$ and instead of (\ref{eqP1})
 \beqn\label{eqP21}
p^n_{\{0\}}(x,y) \geq c_2 \bigg( \sum_{w=1}^{x\wedge|y|} p(-w)V_{{\rm as}}(w)w^{\a}\bigg)\frac{D_n(x,-y)\vee D_n(-y,x)}{n^{2-1/\a}}
\eeqn
where $D_n(x,z)= U_{{\rm ds}}(x)x^{-1} \{z n^{1-2/\a} + a(z)\}$ for $x, z>0$.

(b) \, Let  $\ga =2-\a$ and  $E|\hat Z|<\infty$ so that $U_{{\rm ds}}(x)\asymp x$. Suppose  that $F(x)$ is regularly varying as $x\to -\infty$ of index $ -\beta$ (necessarily  $\beta\geq \a$).  
Then $a(x) \asymp \sum_{w=1}^{x}\sum_{z=1}^\infty  p(-w-z)[ V_{{\rm as}}(z)]^2$ (\cite[Theorem 2(i), (iii)]{Uladd})
and we deduce  $ \sum_{w=1}^{x} p(-w)w^{2\a-1} \asymp  a(x)$ 
so that in view of Proposition \ref{prop5.2}(i) (given shortly) the lower bound (\ref{eqP1})  is exact. 
 \end{rem}

 \v2

\begin{lem}\label{lem5.1}  Suppose $\ga=2-\a$.  For each $M>1$ there exists a constant $C_M$ such that 
\v2
{\rm (i)} \; ${\displaystyle      p_{\{0\}}^n(x,y) \leq C_M\bigg(  \frac{a^\dagger(x)}{n}  +  \frac{(x_n)_+}{n^{1/\a}} \bigg) (y_n^{\a-1} \wedge y_n^{1-\a})  }$   \qquad  $( |x|\leq Mn^{1/\a}$, $y>0)$,
\v2\n
and that if $E|\hat Z|<\infty$, 
\v2
{\rm (ii)}\;\; ${\displaystyle p_{(-\infty,0]}^n(x,y) \leq C_M n^{-1/\a}x_n  [\, y_n^{\a-1} \wedge y_n^{-1}] }$ \qquad  for\;   $0\leq x\leq Mn^{1/\a}$  and $y>0$.
\end{lem}
\n
\pf\,  Let $ |x_n|<M$.   By Theorem \ref{thm5} (cf.  (\ref{R3f2})) as before we have
 $$
 p_{\{0\}}^n(x,y) \asymp a(-y)f^x(n) \asymp  [a^\dagger(x) + (x_n)_+n^{1-1/\a}]n^{-1} y_n^{\a-1} \;\;\mbox{for }\;\; 0< y_n \leq 3M$$
 and   for the proof of  (i)   it  therefore suffices  to show that
for some constant  $C$,
\beqn\label{y>>n}
 p_{\{0\}}^n(x,y) \leq C [a(x) + (x_n)_+n^{1-1/\a}]n^{-1}/ y_n^{\a-1}   \quad \mbox{for} \quad    y_n> 3M .
\eeqn 
Putting $R=\lfloor y/2 \rfloor +1$, $N= \lfloor n/2 \rfloor$ we make the decomposition.
\begin{eqnarray}\label{eq_L5.1}
p_{\{0\}}^n(x,y) &=&
\sum_{k=1}^{N} \sum_{z\geq R}  P[S_{k}^x=z,  \sigma^x_{[R,\infty)} =k >\sigma^x_{\{0\}}]p^{n-k}_{\{0\}}(z,y) \nonumber\\
&&+ \sum_{z < R}P[  \sigma^x_{[R,\infty)}\wedge \sigma^x_{\{0\}} > N, S^x_{N} =z]p_{\{0\}}^{n- N}(z,y)  \nonumber\\
&=& J_1+ J_2\quad \mbox{(say)}. 
\end{eqnarray}  

By the  bound  $p^n(w) =O(n^{-1/\a})$ (valid for all $w\in \Z$)  it follows that  
\begin{equation}\label{0i0}
J_1\leq C P[ \sigma^x_{[R,\infty)} <\sigma^x_{\{0\}}] /n^{1/\a}.
\end{equation}
On using  Lemma \ref{lem7.8}
\begin{eqnarray}\label{(A)} 
P[ \sigma^x_{[R,\infty)} <\sigma^x_{\{0\}}] \leq C P[\sigma^x_{\{R\}} <\sigma^x_{\{0\}}] & =&C \frac{a^\dagger(x)+a(-R)-a(x-R)}{a(R)+a(-R)} \\
&\leq& C'[a^\dagger(x)R^{-\a+1} + x_+ R^{-1}], \nonumber
\end{eqnarray}
where Lemma \ref{lem3.1}(ii) is applied to  estimate   the increment of $a$  for the inequality (as for the equality see (\ref{id_h_p})).
These  together lead to 
$$J_1 \leq C''[a^\dagger(x)y^{-\a+1}/n^{1/\alpha} + (x_n)_+ /y] \leq C(a^\dagger(x) + (x_n)_+ n^{1 - 1/\a})/ny_n^{\a-1} .$$
On the other  hand   by employing  the bound $p^n(x)\leq Cn^{1-1/\alpha}/|x|^\alpha$ (see Lemma \ref{lem7.6})
$$J_2 \leq P[\sigma^x_{\{0\}} \geq \tst 12 n]\sup_{z\leq R} p_{\{0\}}^{n-N}(z,y) \leq C [nf^x(n)] n^{1-1/\a}/y^\a \leq C' [a^\dagger(x) + (x_n)_+n^{1-1/\a}]/y^\a, $$
and hence  (\ref{y>>n}) is obtained.  Thus (i) has been proved. 

(ii) is derived in a similar way; we define $J_1$ and $J_2$ analogously. From $\gamma =2-\alpha$  we have  $\lim P[S_n^0>0] =1/\alpha$ which together with  $E|\hat Z|<\infty$ entails  $P[\sigma^x_{(-\infty, 0]}\geq \frac12 n] \leq C xn^{-1/\alpha}$  so that $J_2 \leq C' n^{-1/\alpha} x_n /y_n^\alpha$. For the estimation of  $J_1$ 
we use, instead of  (\ref{(A)}), 
$$P[\sigma^x_{[R,\infty)}<\sigma^x_{(-\infty,0]} ] \leq C [V_{{\rm as}}(R)- V_{{\rm as}}(R-x)]/V_{{\rm as}}(R)$$ 
 as $R\to\infty$ uniformly for $0<x<R$ (cf. \cite[Remark 5.1]{Uladd} for the first relation).  With this we take an average in the bound corresponding to (\ref{0i0})  to see 
 \[
 J_1\leq C' \int_R^{2R}\frac{V_{{\rm as}}(r)- V_{{\rm as}}(r-x)}{ n^{\alpha}V_{{\rm as}}(r)} \cdot \frac{dr}{R} \leq \frac{C'}{ n^{\alpha} R^{\alpha-1}} \bigg(\int_{2R-x}^{2R}+\int_{R-x}^R\bigg) r^{\alpha-1} \frac{dr}{R} \leq \frac{C'' x}{ n^{\alpha}y},
 \]
showing the bound of (ii). \qed


\v2
In the next lemma  $\gamma$ may be any admissible constant.
\begin{lem}\label{lem5.2}\, For each  $M>1$, there exists a constant $C_M$ such that for all  $n\geq 1$, $x\in \Z$,
\[
p^n_{\{0\}}(x,y)  \leq C_M |y|^{\alpha -1}/|x|^\alpha 
\quad\mbox{if}\quad |x_n|>1 \;\; \mbox{and} \;\;   |y_n|<M. \]
\end{lem}
\n
\pf\, We prove the bound of the lemma in  the dual
form which  is given as
\beqn\label{eqL5.4}
p^n_{\{0\}}(x,y) \leq C_M |x|^{\alpha-1 }/y^\alpha  \quad  \mbox{for} \quad |x_n|<M,  |y_n| >1.
\eeqn
The proof is carried out  by examining  the proof of Lemma \ref{lem5.1}. We can suppose $y_n> 3M$ by symmetry and Theorems \ref{thm2} and \ref{thm3} and   let $R=\lfloor y/2\rfloor$, $N=\lfloor n/2\rfloor$, and
$J_1$ and $J_2$ be defined as in  the proof of Lemma \ref{lem5.1}.   We have  shown that $J_2$ admits   the same  upper bound as required for $p^n_{\{0\}}(x,y)$ in (\ref{eqL5.4})  which though presented in case $\gamma =2-\alpha$ applies to the other case too.  As for $J_1$ we  first recall
 \[J_1 =\sum_{k=1}^{N} \sum_{z\geq R}  P[S_{k}^x=z,  \sigma^x_{[R,\infty)} =k >\sigma^x_{\{0\}}]p^{n-k}_{\{0\}}(z,y).
 \]
 The double sum with the  additional restriction $|z-y| > \frac12 R$ to the inner sum is dominated by a constant multiple of
 \[
P[  \sigma^x_{[R,\infty)}  >\sigma^x_{\{0\}}]\frac{n^{1-1/\alpha}}{y^{\alpha}}\leq  C\frac{|x|^{\alpha-1}}{y^{\alpha-1}}\cdot \frac{n^{1-1/\alpha}}{y^{\alpha}} =C\frac{ |x|^{\alpha-1}}{y_n^{\alpha-1}y^\alpha},
 \]
 where we have used Lemmas \ref{lem7.6} and \ref{lem7.9}.  
On writing down the probability under the  double summation sign  by means of transition probabilities   it suffices to
 show that
 \beqn\label{eqL5.41}
 \sum_{w<R} \sum_{\frac32 R\leq z \leq 3R} \sum_{k=0}^{N-1} p^k_{\{0\}}(x,w)p(z-w)p^n(y-z) 
 \leq C  |x|^{\alpha-1}/ y^{\alpha }.
 \eeqn  
By the trivial bound $\sum_{k=0}^Np^k_{\{0\}}(x,w)\leq g_{\{0\}}(x,w) \leq C|x|^{\alpha-1}$ and  $\sum_{w<R} p(z-w) \leq C R^{-\alpha}$ for $z>\frac32 R$   the above triple sum is bounded by
\[
\frac{C'|x|^{\alpha-1}}{y^{\alpha} } \sum_{\frac32 R\leq z \leq 3R} p^n(y-z) 
\]
On using Lemma \ref{lem7.6} again the sum above is bounded by a constant multiple of 
\[
\sum_{z: |y-z|\leq n^{1/\alpha}} n^{-1/\alpha} +\sum_{z: |y-z|>n^{1/\alpha}} \frac{n^{1-1/\alpha}}{|y-z|^\alpha}\leq 2+ 2n^{1-1/\alpha}\int_{n^{1/\alpha}-1}^\infty u^{-\alpha}du \leq C'',
\]
showing (\ref{eqL5.41}) as required. \qed



\v2
\begin{lem}\label   {cor5.1}  Suppose  $\ga = \a-2$ and define $\omega_{n,x,y}$ for $x\neq 0$ and $y>0$ via
\beqn \label{iv}
 p_{\{0\}}^n(x,y)= a(-y)f^x(n)\om_{n,x,y}.
 \eeqn
Then,  $\omega_{n,x,y}$ is dominated by a constant multiple of  $1 \wedge  y_n^{-2\a+1}$ (in particular uniformly bounded), and tends to unity as $y_n \to 0$   and $n\to\infty$ uniformly for  $0<x < Mn^{1/\a}$  for each  $M>1$. 
\end{lem}
\n
\pf\, The convergence of $\omega_{n,x,y}$  to zero follows from Theorems \ref{thm4} and
 \ref{thm5} (the first  case) and  the stated of $\om_{n,x,y}$ is derived from Lemma \ref{lem5.1}(i) with a simple manipulation. 
 \qed
\v2
In the sequel we use the notation $H_B^x(y) = P[S^x_{\sigma_B}=y] $, $B\subset \Z$. It follows that
\beqn\label{H_B}
H_B^x(y) =  \sum_{z\notin B}g_B(x,z) p(y-z) \quad \mbox{for}\quad y\in B.
\eeqn

\begin{Prop}\label{prop5.2}\, Suppose $\ga=2-\a$. Then for some constant $C$
\v2
{\rm (i)}  \quad  ${\displaystyle  p^n_{\{0\}}(x,y) \leq  
C\bigg[\frac{a^\dagger(x)a(-y)}{ n^{2-1/\a}} +\frac{ a^\dagger (x)(|y_n|\wedge 1) +  a(-y)(x_n\wedge1)} { n} \bigg]  } \qquad  (x\geq 0, y\leq -1) $,
\v2
{\rm (ii)}  \quad  ${\displaystyle  p^n_{\{0\}}(x,y) \leq  C n^{-1}\Big[a(x)(y_n^{\a-1}\wedge 1) + a(-y)(|x_n|^{\a-1}\wedge1) \Big]  } \qquad (x \leq -1, y\geq 1)$.
\v2
\end{Prop}

\n
\pf\,  
 First we prove (ii). 
The proof  is based on  Lemmas \ref{lem5.1}(i) and \ref{lem5.2}, that entail for $k\leq n/2$
 \begin{eqnarray} \label{P2eq1}
&p^{n-k}_{\{0\}}(z,y)  \leq  C[a(z)/n+z/n^{2/\a}] (y_n^{\a-1}\wedge y_n^{1-\a}) \;   &(0<z<n^{1/\a},  y >0)\;\;  \mbox{and} \qquad \\
&p^{n-k}_{\{0\}}(z,y)  \leq  Cy^{\a-1}/z^{\a}  \quad  & (z>n^{1/\a}, 0<  y <n^{1/\a}),
\label{P2eq2}
\end{eqnarray}
respectively.  
Let $x\leq -1$ and $y\geq 1$ and consider the RHS of the trivial  inequality
\beqn\label{eqP2_4}
P[\sigma^x_{[0,\infty)} \leq n/2, \sigma^x_{\{0\}} >n, S^x_n= y] 
\leq
\sum_{z>0}H_{[0,\infty)}^x(z)\sup_{k\leq \frac12 n}p^{n-k}_{\{0\}}(z,y). 
\eeqn
By  $g_{[0,\infty)}(x,z) \leq g_{\{0\}}(x,x))$ it it follows that  $H^x_{[0,\infty)}(z) \leq Ca(x)P[X \geq z]$ (cf. (\ref{(C)})), and hence  
$$\sum_{0<z<n^{1/\a}}H_{[0,\infty)}^x(z)z \leq C' a(x) \sum_{1\leq z< n^{1/\a} } z^{1-\a}  
\leq C_1 a(x) n^{-1+2/\a}.$$
Since $H^x_{[0,\infty)}\{a\} = O(a(x))$ (cf. \cite{Uladd}),  this together with  (\ref{P2eq1}) shows
$$\sum_{0< z< n^{1/\a}}H_{[0,\infty)}^x(z)p^{n-k}_{\{0\}}(z,y)\leq C'' n^{-1}a(x)(y_n^{\a-1}\wedge y_n^{1-\a}).$$
In a similar way
$$ \sum_{z\geq n^{1/\a}}H_{[0,\infty)}^x(z)z^{-\a} <C'a(x)n^{-2+1/\a} \quad\mbox{and}\quad  \sum_{z\geq n^{1/\a}}H_{[0,\infty)}^x(z)  <C'a(x)n^{-1+1/\a}$$
and by (\ref{P2eq2}) and the bound $p^{n-k}_{\{0\}}(z,y)  \leq C/n^{-1/\a}$ (following the local limit theorem)
$$\sum_{z\geq  n^{1/\alpha}}H_{[0,\infty)}^x(z)\sup_{k\leq n/2}p^{n-k}_{\{0\}}(z,y)\leq C'' a(x)(y^{\alpha-1}n^{-2+1/\alpha} \wedge n^{-1}) = C'' n^{-1}a(x)(y_n^{\a-1}\wedge 1).$$
Thus the RHS of (\ref{eqP2_4}) is bounded by a constant multiple of $ n^{-1}a(x)(y_n^{\a-1}\wedge 1).$

Let $\hat S^x$ and $\hat \sigma_B^x$ denote the dual walk and its hitting time, respectively. It then follows that
\begin{equation}\label{eq_3.1}
\begin{array}{ll}
p^n_{\{0\}}(x,y) -P[\sigma^x_{[0,\infty)} \leq n/2, \sigma^x_{\{0\}} >n, S^x_n= y]  \\[2mm]
 \qquad\qquad  \leq 
 P[\hat \sigma^y_{(-\infty,0]} \leq n/2, \hat \sigma^y_{\{0\}} >n, \hat S^y_n= x].
 \end{array}
\end{equation}
By duality relation the probability on the RHS is the same as what we have just estimated but with $x$ and $y$ replaced by $-y$ and $-x$, respectively, and hence dominated by a constant multiple of  $ n^{-1}a(-y)(|x_n|^{\a-1}\wedge 1)$.
This concludes   (ii). 

For the proof of (i)  we apply (\ref{P2eq1})  with $z, y$ replaced by $-y, -z$, in which  
 we may replace $(-z_n)^{\a-1}\wedge (-z_n)^{1-\a}$ by   
$a(z)/n^{1-1/\a}$  ($z <0$)  to obtain
\beqn\label{P2_1}
p^{n-k}_{\{0\}}(z,y) =p^{n-k}_{\{0\}}(-y,-z) \leq Ca(z)\frac{a(-y)+(|y_n|\wedge 1)n^{1-1/\a}}{n^{2-1/\a}}  \quad  (z<0, k\leq n/2).
\eeqn
This is valid at least  for all $-n^{1/\a} <y<0$ and  can be extended to $y\leq -n^{1/\a}$.  For the proof of the extension we have only to observe that if  $y\leq -n^{1/\a}$, then  the RHS is not less than  $c/n^{1/\a}$ if $z<-n^{1/\a}$ with some $c>0$  while for $z\leq -n^{1/\a}$, (\ref{P2_1})  follows from Lemma \ref{lem5.2} (note  that $1/y^\alpha= O(1/n)$).
  Since $E[S_{\sigma[1,\infty)}]=\infty$ we have  $H^x_{(-\infty,0]}\{a\} = a(x)$ (cf. \cite{Uladd} or (\ref{uA/H})), and  from (\ref{P2_1}) we deduce 
\beq
P[\sigma^x_{(-\infty,0]} \leq n/2, \sigma^x_{\{0\}} >n, S^x_n= y] 
&\leq&
\sum_{z<0}H_{(-\infty,0]}^x(z) \sup_{k\leq  n/2}p^{n-k}_{\{0\}}(z,y) \\
&\leq& Ca(x)\frac{a(-y)+(|y_n|\wedge 1)n^{1-1/\a}}{n^{2-1/\a} }. 
\eeq
By the analogue of (\ref{eq_3.1}) we conclude (i) by  duality relation as above. \qed




\v2
The  next lemma concerns the hitting distribution of the negative half line defined by
$$h^x(n,y) = P[\sigma^x_{(-\infty,0]} =n, S^x_n=y].$$
 
\begin{lem}\label{lem5.4} Suppose $E|\hat Z|<\infty$. Then,

{\rm (i)} for $M>1$ and $\e>0$,  uniformly for $0\leq x_n<M$ and $y\leq 0$
$$h^x(n, y) = \frac{x_n \mathfrak{p}_{c_\circ}(-x_n)}{n} \Big[ H^\infty_{(-\infty,0]}(y)\{1+o_\e(1)\} + r(n,y)\Big]$$
where $o_\e(1)$ is bounded and tend to zero  as $n\to\infty$ and $\e\to 0$ in this order and
$$|r(n,y)| \leq  C_M\sum_{z>\e n^{1/\a}} z^{\a-1}p(y-z)$$
for a constant $C_M$   depending  only on  $M$ and $F$; and

 {\rm (ii)} there exists a constant $C$ such that for all $x\geq 1, y <0$ and  $n\geq 1$,
$$h^x(n, y) \leq C   (n^{-1}\wedge x^{-\a})x_n  H^\infty_{(-\infty,0]}(y).$$
 \end{lem}
\n
\pf\,  Let $\e>0$ and  in the expression  
\beqn\label{++}
h^x(n+1,y) =\sum_{z=1}^\infty p^{n}_{(-\infty,0]}(x,z)p(y-z)
\eeqn
 we divide the sum into two parts,  the sum on  $z< \e n^{1/\a}$  and 
the remainder which  are denoted by $\Sigma_{<\e n^{1/\a}}$ and    $\Sigma_{\geq\e n^{1/\a}}$, respectively.  
By Doney's result (\ref{Doney}) and (\ref{V/a}) it follows that
$$p^{n}_{(-\infty,0]}(x,z) =(E|\hat Z|)^{-1}V_{{\rm as}}(z) x_n\mathfrak{p}_{c_\circ}(-x_n)n^{-1}\{1+o_\e(1)\}
\quad \mbox{uniformly for}  \;\; z\leq \e n^{1/\a},$$
(note $L(n^{1/\alpha})\to 1/E|\hat Z|$ in  (\ref{V/a})) and substituting  this and  using
\beqn\label{bd_H}
      \frac1{E|\hat Z|} \sum_{z=1}^\infty V_{{\rm as}}(z-1)p(y-z) = H^{+\infty}_{(-\infty,0]}(y)  \quad (y\leq 0),
 \eeqn
(cf. \cite[Eq(2.7)]{Uladd}) we deduce that  
$$\Sigma_{<\e n^{1/\a}} 
= \frac{x_n \mathfrak{p}_{c_\circ}(-x_n)}{n}\bigg[H^\infty_{(-\infty,0]}(y) - \frac1{E|\hat Z|}  \sum_{z>\e n^{1/\a}} V_{{\rm as}}(z)p(y-z)  \bigg] \{1+o_\e(1)\}. $$
By Lemma \ref{lem5.1}(ii)
$$\Sigma_{\geq\e n^{1/\a}} \leq  C\frac{x_n }{n} \sum_{z>\e n^{1/\a}} z^{\a-1}p(y-z), $$ 
and on noting $V_{{\rm as}}(z) \asymp z^{\a-1}$ the assertion  (i)  follows. It in particular follows that 
\beqn\label{bd_h1}
h^x(n, y) \leq C   n^{-1}x_n  H^\infty_{(-\infty,0]}(y)\quad \mbox{for}\quad  0\leq x\leq M n^{1/\a},  y\leq 0,
\eeqn 
since   $|r(n,y)| \leq C_M'H_{(-\infty,0]}^{+\infty}(y)$ in view of (\ref{bd_H}). 

For the proof of  (ii)  we replace $p_{(-\infty,0]}^{n}(x,z) $ by $p_{\{0\}}^{n}(x,z)$ in (\ref{++}) to have an upper bound and  verify that for $y\leq 0$ and $x\geq 2n^{1/\a}$,
\beqn\label{upb-h2}
\sum_{z=1}^\infty p_{\{0\}}^{n}(x,z)p(y-z)\leq
\frac{C}{x^{\a}}H_{(-\infty,0]}^{+\infty}(y) +\frac{C}{n^{1/\a}}F(y-{\textstyle \frac12}x).
\eeqn
For verification of (\ref{upb-h2}) we
break the range of summation  into three parts $0<z\leq n^{1/\a}$,  $n^{1/\a} <z\leq x/2$ and $z>x/2$,
and denote the corresponding sums by $J_1$, $J_2$ and $J_3$, respectively.
It is immediate  from Lemma \ref{lem5.2} and (\ref{bd_H})  that  $J_1 \leq Cx^{-\alpha}H_{(-\infty,0]}^{+\infty}(y)$.  
By the  bound $p^n(z) \leq C n^{1-1/a}/|z|^\a$ (cf. Lemma \ref{lem7.6}) 
it follows  that $p_{\{0\}}^{n}(x,z)\leq  Cn^{1-1/\a}x^{-\a}$ for  $n^{1/\a}<z\leq x/2$, which combined with  the  bound
\beqn\label{eq4.-1} \sum_{z> n^{1/\a}}p(y-z) \leq \sum_{z=1}^\infty \bigg(\frac{z}{n^{1/\alpha}}\bigg)^{\a-1} p(y-z) \leq \frac{2E|\hat Z|}{ n^{1-1/\a}}\bigg[\sup_{z\geq 1}\frac{z^{\a-1}}{V_{{\rm as}}(z)}\bigg] H_{(-\infty,0]}^{+\infty}(y)
\eeqn
yields $J_2= C'x^{-\a} H_{(-\infty,0]}^{+\infty}(y)$. Finally 
$J_3 \leq C n^{-1/\a}F(y-x/2)$.
 These
estimates together verify  (\ref{upb-h2}). 
As in (\ref{eq4.-1})  we derive $F(y-{\textstyle \frac12}x)\leq C_1 H^{+\infty}_{(-\infty,0]}(y)/x^{\a-1}$.  
Hence
$$h^x(n,y) \leq C_1 H_{(-\infty,0]}^{+\infty}(y)/x^{\a-1}n^{1/\a}\quad \mbox{for}\quad x> 2n^{1/\a},  y<0, n\geq 1, $$
which  combined with (\ref{bd_h1})   shows the bound in  (ii).  The proof of Lemma \ref{lem5.4} is complete.
\qed
\v2

{\bf Proof of Theorem \ref{thm6}.} If either $x$ or $y$ remains in a bouded set,  the formula  (i) of Theorem \ref{thm6} agrees with  that of Theorem \ref{thm5}, so that
we may and do suppose both $x$ and $-y$ tend to infinity. Note that  the second ratio on the RHS of (i)  is  then asymptotically equivalent to the ratio in (ii), hence  (i) and (ii)  of Theorem \ref{thm6} is written as a single formula. Let $c_\circ =1$ for simplicity and put
$$\Phi(t;\xi) = t^{-1}\xi \mathfrak{p}_t(-\xi).$$
Then  what  is to be shown may be stated as follows:  as  $n\to\infty$ and $x \vee (-y)\to\infty$
\beqn \label{eq_thm5}
p_{\{0\}}^n(x,y) \sim  \k_{\a,\ga} a(x)a(-y)n^{-2+1/\a} + C^+ \Phi (n; x-y)  
\eeqn 
 uniformly for $-M< y_n< 0< x_n <M$, provided $0< C^+= \lim_{z\to -\infty} a(z)<\infty$. 

We follow the proof in \cite{U1dm} to the corresponding result.   We employ the representation
\beqn\label{eqT5}
p_{\{0\}}^n(x,y)=\sum_{k=1}^n \sum_{z<0}h^x(k,z)p_{\{0\}}^{n-k}(z,y).
\eeqn
Break  the RHS into three parts by partitioning the range of the first summation as follows
\beqn\label{eq4.5}
1\leq  k< \e n; ~~ \e n \leq k \leq (1-\e)n;~~ (1-\e)n<k\leq n
\eeqn
and call the corresponding sums $I,~II $ and $I\!I\!I$, respectively. Here $\e$ is a positive constant that will be chosen small. 
The proof is divided into two cases corresponding to (i) and (ii).
\v2

{\it Case  $x_n\wedge y_n \to 0$}:
By duality one may suppose that $x_n\to 0$. From  $E Z=\infty$ and (\ref{a_bdd})  it follows \cite[Theorems 1 and 2]{Uladd} that
\beqn\label{Ha}
\sum_{z\leq 0} H_{(-\infty,0]}^x(z)a(z) = a(x) \quad\mbox{and}\quad C^+= \sum_{z\leq 0} H_{(-\infty,0]}^{+\infty}(z)a(z) <\infty,
\eeqn
respectively.
From the latter bound above and Lemma \ref{lem5.4} (or (\ref{bd_h1}))  one deduces, 
\beqn\label{hh}
(*) \qquad \sum_{k\geq \e n} \sum_{z<0} h^x(k,z)a(z)\leq    C_{\e} x_n
\eeqn
with a constant $C_\e$ depending on $\e$.
As the dual of  (\ref{iv}) of Lemma  \ref{cor5.1} we have 
\beqn\label{pf}
p_{\{0\}}^n(z,y) = a(z)f^{-y}(n)\{1+r_{n,z,y}\} \quad  (z<0, -Mn^{1/\a} <y < 0)
\eeqn
where $r_{n,z,y}$ is uniformly bounded and 
tends to zero as $z/n^{1/\a}\to 0$ and $n\to \infty$  uniformly for  $y$,
which together with   (\ref{hh}) shows 
$$II \leq C_{\e,M} x_n f^{-y}(n).$$


Similarly on using  (\ref{pf})  above
$$I=\sum_{1\le k< \e n}\,\sum_{z=-\infty}^{-1}h^x(k,z) a(z)f^{-y}(n-k)\{1+r_{n-k,z,y}\}.$$
For the evaluation of the last double sum we may replace $f^{-y}(n-k)$ by
 $f^{-y}(n) (1+O(\e))$, and the contribution to it of $r_{n-k,z,y}$ is negligible since  $\sum_{z> N} H^x_{(-\infty,0]}(x) a(z) \to 0$ ($N\to\infty$)  uniformly in $x$ in view of the second relation of (\ref{Ha}). By (\ref{hh}) the summation over $z$ may be extended to the whole half line $k\geq 1$. Now  applying  the first relation of (\ref{Ha})  we find
$$ I= a(x) f^{-y}(n)\{1+O(\e) +o(1)\}.$$

As for $ I\!I\!I$   first observe that by (\ref{pf}) and Theorem \ref{thm3}
 $$\sum_{k=1}^{\e n} p_{\{0\}}^k(z,y)=g(z, y)-r_n\leq C(a(z)\wedge a(y)) \quad \mbox{with}\quad 0\leq r_n\leq C_\e a(z)f^{-y}(n)n$$
 ($y, z<0$).
 If $y_n$ is bounded away from zero so that $x/y\to 0$, then  $ I\!I\!I=O(x_n/n)=o(y_n/n)$. On the other hand, applying Lemma \ref{lem5.4} we  see that if $y_n\to 0$,
 $$ I\!I\!I= x_n\mathfrak{p}_{1}(x_n) n^{-1}\sum_{z<0}H_{(-\infty,0]}^\infty(z)g(z, y)(1+O(\e))+ O(x_n f^{-y}(n) ),$$
whereas  by (\ref{Ha})  and the  subadditivity of $a$ we infer that  $\sum_{z\leq 0} H_{(-\infty,0]}^\infty(z)g(z, y) \to  C^+$ as $y\to\infty$. Hence
  $$ I\!I\!I= x_n\mathfrak{p}_{1}(x_n)n^{-1}( C^+ +o(1) +O(\e))+  O(x_n f^{-y}(n)).$$
Adding these expressions of $I$, $II$ and $I\!I\!I$ yields the desired formula, because of  arbitrariness of $\e$ as well as  the identity $x_n\mathfrak{p}_{1}(x_n)/n =\Phi(n;x)$. 
\v2

{\sc Case $x_n\wedge (-y_n) \geq 1/M$.}
 By  Lemma \ref{lem5.4}(ii) and (\ref{pf}) it follows that   in this regime 
$$I\leq  \sum_{1\leq k <\e n}\frac C {k^{1/\a}\,x^{\a-1}}\sum_{z<0}  H_{(-\infty,0]}^\infty(z) a(z)f^{-y}(n) \leq C' \frac{\e^{1-1/\a}}n.$$

 For evaluation of $I\!I\!I$  change the variable $k$ into $n-k$ and  apply Lemma \ref{lem5.2} to $p^k_{\{0\}}(-y,-z)$ (with $(-y,z)$ in place of $(x,y)$)  to see that   for any $\de>0$ 
\beq
\sum_{k=1}^{\e n} p^k_{\{0\}}(z,y) &\leq& C \sum_{k\leq \de |z|^\a} k^{-1/\a} + C_\de \sum_{\de |z|^\a<k <\e n} |z|^{\a-1}/|y|^{\a}\\
&\leq&
C(\de |z|^\a)^{1-1/\a} + C_\de (\e n) |z|^{\a-1}/|y|^{\a},
\eeq
where $C_\de$ may depend on $\de$ but $C$ does not. Then by Lemma \ref{lem5.4}(ii)
$$I\!I\!I \leq C' n^{-1} [C\de^{1-1/\a} + C_\de \e] \sum_{z<0} H_{(-\infty,0]}^{+\infty}(z) |z|^{\a-1} \leq C'' [C\de^{1-1/\a} + C_\de \e] /n, $$
hence for any $\e'>0$ we can choose $\e>0$ and $\de>0$ so that $I\!I\!I  \leq \e'/n$.

 By  Lemma \ref{lem5.4}(i), (\ref{pf})  and (\ref{Ha})
$$II=\sum_{ \e n\leq k \le (1-\e)n}\frac{ x_k \mathfrak{p}_{1}(-x_k)}{k}\,\sum_{z=-x}^{-1}H_{(-\infty,0]}^\infty(z) p_{\{0\}}^{n-k}(z,y)(1+o_\e(1))+o\bigg( \frac{f^{-y}(n)}{\e^{1/\a}}\bigg).$$
Here (and  in the rest of the proof) the estimate indicated by $o_\e$ may depend on $\e$ but is  uniform in the passage to the limit under  consideration once $\e$ is fixed.


Since $-y_n$ is bounded away from zero as well as from infinity,
  we may replace  $p_{\{0\}}^{n-k}(z,y)$ by $a(z)y_{n-k}\mathfrak{p}_1(y_{n-k})/ (n-k)$ to  see that
\[
II=\sum_{ \e n\leq k \le (1-\e)n}\frac{x_k|y_{n-k}|\mathfrak{p}_{1}(-x_k)\mathfrak{p}_{1}(y_{n-k})}{ k(n-k)}\,\sum_{z=-x}^{-1}H_{(-\infty,0]}^\infty(z)a(z)(1+o_\e(1))+\frac{o(1/n)}{\e^{1/\a}} \nonumber.
\]
 noting the identity $x_k\mathfrak{p}_1(-x_k) =x_n\mathfrak{p}_{k/n}(-x_n) = \Phi(k/n;x_n)k/n$ and similarly for $y_{n-k}\mathfrak{p}_1(y_{n-k})$ and
$$\sum_{ \e n\leq k \le (1-\e)n}\frac{x_k|y_{n-k}|\mathfrak{p}_{1}(-x_k)\mathfrak{p}_{1}(y_{n-k})}{ k(n-k)}\,=\frac{1+o(1)}{n}\int_0^1 \Phi(t;x_n)\Phi(1-t; y_n)dt + O\bigg(\frac{\e}{n}\bigg).$$
Here we have used the fact that $\int_0^\e\mathfrak{p}_t(\xi) \xi dt/t= \int_{\xi/\e^{1/\a}}^\infty = O(\e/\xi^\a)$.
 Since  for $\xi>0$, $\Phi(dt; \xi)dt$ is the distribution  of the  hitting-time  to zero by $\xi+Y$, we have
$$\int_0^1 \Phi(t; x_n)\Phi(1-t; -y_n)dt=\Phi(1; x_n - y_n).$$
Hence
\beqn\label{II}
II= \frac1{n} C^+ \Phi(1; x_n - y_n) \{1+o(1)\}+ O\bigg(\frac{\e}{n}\bigg)+ \frac{o(1/n)}{\e^{1/\a}}.
\eeqn
  (as well as  $nI+n I\!I\!I \to 0$)  as $n\to\infty$ and $\e\to 0$ in this order.  Thus  (\ref{eq_thm5}) is obtained, the first term on the RHS of it being negligible.
  
\v2
{\bf Proof of Proposition \ref{prop2.2}.} The case $C^+=0$ is trivial. If    $0<C^+ <\infty $,  by noting that   Theorem \ref{thm6} and Lemma \ref{lem5.1}(i) (in the dual form ( \ref{P2_1})) together yield
$$
\frac{p_{\{0\}}^{n-k}(z,y)}{p^n_{\{0\}}(x,y)} \leq C \frac{a(z)[1 + |y|n^{1-{2/n}}]}{1 + |y|n^{1-{2/n}} + xn^{1-{2/n}}} \leq Ca(z)  \quad (z<0, k<n/2)$$
 and  that $H_{(-\infty,0]}^x(z)\leq CH_{(-\infty,0]}^\infty(z)$, we deduce  that  $\sum_{k<n/2}\sum_{z<-R} h^x(k,z)p^{n-k}_{\{0\}}(z,y)$ is at most a constant multiple of $\sum_{z<-R}H_{(-\infty,0]}^\infty(z)a(z)$ which approaches zero   as $R\to\infty$; for the sum over $n/2\leq k \leq n$,  one  uses the bound  $\sum_{n/2\leq k\leq n}p^{n-k}_{\{0\}}(z,y) \leq Ca(z)$ as well as Lemma \ref{lem5.4}(ii) to obtain 
 the same bound in a similar way. This verifies the first half of the asserted formula.

The second half   obviously follows  if $E|\hat Z|=\infty$ so that $H^\infty_{(-\infty,0]}$ vanishes.  
Let $E|\hat Z|<\infty$. Then we can apply  Lemma \ref{lem5.4}(ii) as well as Theorem \ref{thm5} (in the dual form  (\ref{eqC2})) 
to see that  the contribution to the sum (\ref{eqT5}) from   $- R \leq z\leq 0$ is dominated by a positive multiple of 
\[
\sum_{-R\leq z\leq 0} \sup_{k<n/2}\Big[H_{(-\infty,0]}^x(z) p^{n-k}_{\{0\}}(z,y) + h^x(k,z) g_{\{0\}}(z,y)\Big]  \leq C R^{\alpha-1}\bigg[ \frac{a(-y)}{n^{2-1/\a}} +  \frac{x\vee |y|}{n^{1+1/\a}}\bigg],
\]
which is negligible as compared with the lower bound of $p^{n}_{\{0\}}(x,y)$ given by Proposition \ref{prop5.1} provided that  $C^+=\infty$ or,  equivalently,  $\sum_{w\geq 1} w^{2\alpha-1}p(-w)=\infty$. \qed

\v2
{\bf Proof of Proposition \ref{prop2.1}.}
\, In case $|x_n|\leq 3$ the assertion follows from Theorems \ref{thm2} and \ref{thm3}. We let $x_n>3$, the case $x_n<-3$ being treated in the same way.
In the obvious identity 
\beqn\label{f/=}
f^x(n) = \sum_{y} p^{n-1}_{\{0\}}(x,y)p(-y)
\eeqn 
the sum on the RHS over $|y|\leq n^{1/\alpha}$ is bounded by a constant multiple of $x^{-\alpha}$ by virtue of Lemma \ref{lem5.2}.  Since $p^{n}_{\{0\}}(x,y)\leq p^n(y-x)$,  it suffices to show that
\beqn\label{eqP2.1}
\sum_{|y| > n^{1/\alpha}} p^{n}(y-x)p(-y) \leq  Cn^{-1/\alpha} x^{-\alpha}.
\eeqn
We break the sum into three parts by splitting the range of summation at $y=x\pm n^{1/\alpha}$ and denote them by $\Sigma_{|x-y|<n^{1/\alpha}}$, $\Sigma_{n^{1/\alpha} <y \leq x -n^{1/\alpha} }$ and  $\Sigma_{y \geq x +n^{1/\alpha} }$.
The first sum is estimated as follows:
\[\Sigma_{|x-y|<n^{1/\alpha}} \leq C  n^{-1/\alpha}\sum_{|y-x| < n^{1/\alpha}} p(y) = n^{-1/\alpha} x^{-\alpha} \times o(1).
\]
For the second sum  we  further split its range of summation at  $y= x/2$ and apply  Lemma \ref{lem7.6}  to see that $\Sigma_{n^{1/\alpha} <y \leq  x -n^{1/\alpha} }$ is at most a constant multiple of 
\[
 \sum_{n^{1/\alpha} \leq y <x/2}p(-y)\frac{n^{1-1/\alpha} }{ x^{\alpha}}  +  \sum_{x/2<y \leq  x -n^{1/\alpha}} p(-y)  \frac{ n^{1-1/\alpha} }{(x-y)^\alpha} \leq  \frac{Cn^{-1/\alpha}}{x^\alpha}.
 \]
  The third sum is evaluated to be $o(1/x^\alpha)$ in a similar  way. Thus (\ref{eqP2.1}) and hence Proposition \ref{prop2.1} has been verified. \qed
 



\section{Extension to an arbitrary finite set}

Let $A$ be a finite non-empty subset of $\Z$. The function $u_A(x), x\in \Z$  defined  in  (\ref{def_u}) may be  given by 
\beqn\label{u/g}
u_A(x) = g_A(x,y) + a(x-y) - E_x[a(S_{\sigma(A)} -y)]
\eeqn
(whether (\ref{HA}) is assumed or not), the RHS being independent of $y\in \Z$ (cf. \cite[Lemma 3.1]{U1dm_f}, \cite{PS}) and  the difference of the last two terms in it tending to zero as $|y|\to \infty$.
Taking an arbitrary $w_0\in A$ for $y$ it in particular follows that
\beqn\label{uA/a}
u_A(x) = a^\dagger(x-w_0) - E[a(S^x_{\sigma_A} -w_0)]
\eeqn
so that $u_A(x)\sim a(x)$ as $x\to \pm \infty$ if $a(x)\to\infty$ as $x\to \pm \infty$.
The function $u_A$ is harmonic for the semi-group $p_A^n$ as noted previously, and  $u_A(S^x_{n})\1( n <  \sigma^x_A)$ is accordingly a martingale for each $x\in \Z$.
Put $f_A^x(n) = P[\sigma_A =n]$.  We state  only the extensions corresponding to those given in  Theorem \ref{thm3} (restricted to the  case $\ga=2-\a$) and Theorem  \ref{thm5}. %
In the following theorem we include the case of  periodic walks (i.e., the condition 2) stated in Section 1 may be violated).  What is stated about  (\ref{u/g}) also  holds for the periodic walk. 
  
\begin{Thm}\label{thm7}\,  Let $\nu\geq 1$ denote  the  period of the walk, which amount to assume  (in addition to (\ref{f_hyp})) that  $p^{\nu n}(0) >0$  and $p^{\nu n+j }(0) =0$ ($1\leq j < \nu)$ for all sufficiently large  $n$.  Let  $\ga=2-\a$  and    $M$ be any number greater than $1$. Then, 

{\rm (i)}  for $x$ with $P[S^x_n\in A]>0$,  as $n\to\infty$
 \beqn\label{A1-0}
f_A^x(n) \sim \left\{ \begin{array} {lr}
 u_A(x) f_A^0(n) + \nu x_n\mathfrak{p}_{c_\circ}(-x_n)/n   &(0 \leq x_n<M), \\[1mm]
 u_A(x) f_A^0(n) \quad &(x<0, x_n\to 0), \\[1mm]
\nu c_\circ\mathfrak{f}^{x_n}(c_\circ)/n \quad &(-1/M\leq  x_n< -M)
\end{array} \right.
\eeqn
 and
 \beqn\label{f/nu}
f_A^0(n) \sim   f^0(\lfloor n/\nu\rfloor\nu)\sim \nu \kappa_{\alpha,\gamma}c_\circ^{1/\alpha}/n^{2-1/\alpha};
 \eeqn

 {\rm (ii)}   uniformly  for $|x|<Mn^{1/\a}$ and $-M<y<Mn^{1/\a}$ with $u_{-A}(-y)>0$, $p^n(y-x)>0$,
 as $n\to\infty$
\beqn\label{A1-1}
p^n_{A}(x,y) \sim \left\{\begin{array}{lr}{\displaystyle   f^x_A(n)u_{-A}(-y)  } \quad  & (y_n \to 0),
\\[1mm]
{\displaystyle   u_A(x)f_{-A}^{-y}(n) + \nu \frac{(x_n)_+ K_{c_\circ}(y_n)}{n^{1/\a} }}  & (x_n \to 0),\\[2mm]
 \nu \mathfrak{p}^{\{0\}}_{c_\circ n} (x,y)  &  (|x_n| \wedge y_n\geq 1/M).
\end{array} \right.\eeqn
\end{Thm}
\begin{rem}\label{R8}\,
(a)  If  (\ref{HA}) is violated,  then    $u_{-A}(-y) =0$ for $y\leq \min A$ and (\ref{A1-1}) says nothing about   $p^n_{A}(x,y)$ which is positive for $x\leq \min A$ and whose asymptotic form is   found in   the dual version of (\ref{A1-1}) deduced by using $ p_A^{n}(x,y)=p_{-A}^{n}(-y, -x) $.  

(b)\, 
 The results for the periodic walks  are derived from those of the aperiodic ones.
  The process  $\tilde S_n= S_{\nu n}/\nu$, $n=0, 1,\ldots$ is a strongly aperiodic walk   on $\Z$ such that  $1-E[e^{i \th \tilde S_n}] \sim \nu^{1-\alpha}[1-\phi(\th)]$; hence   in case  $A=\{0\}$,
     the results restricted on $ \nu\Z$ follow immediately from  those of the aperiodic walks and the extension to $\Z$ is then readily performed by using $E[a(S_1^x)]=a^\dagger(x)$.   
  The  general case is reduced to the case  $A=\{0\}$
  in the same way as for  aperiodic walks as is described below.  


\end{rem}
\v2
  
The basic idea of proof is the same as in \cite{U1dm_f}, the details are rather  simpler because of uniqueness of  positive harmonic function for the killed walk. In the sequel we may and do assume $\nu=1$  (see Remark \ref{R8}(b)).  

Take an integer   $R>M$ and  let $\tau^x_R$ be the first exit time of $S^x$ from the interval  $(-R, R)$: 
$$\tau_R^x= \sigma^x_{(-\infty,-R] \cup [R,\infty)} =\inf\{n\geq 1: |S^x_n|\geq R\}.$$
Put   $N = mR^\a\lfloor(1+ \lg n)\rfloor$ with  a positive  integer $m$ determined shortly and  decompose
\begin{eqnarray}\label{decomp}
p^n_A(x,y) &=& \sum_{k=1}^{N-1} \sum_{|z|\geq R}P[\tau^x_{R}=k<\sigma^x_A, S^x_k=z]p_A^{n-k}(z,y) \\
&& +\,\, \e(x,y;R),\nonumber
\end{eqnarray}
provided that $1<N< n/2$. 
Here
\beqn\label{e}
\e(x,y;R) = \sum_{z}P[ \tau^x_{R}\wedge \sigma^x_A \geq N, S^x_N =z]p_A^{n-N}(z,y).
\eeqn
Using  the  fact that the process  $Y_t^n:= S_{\lfloor nt\rfloor}/n^{1/\a}$ converges to a stable process   we deduce  that there exists  a constant $\la>0$ such that $\sup_{x: |x|\leq R}P_x[\tau^x_{R} \geq R^\a] < e^{-\la}$ for all sufficiently large $R$, by which we deduce  (cf. \cite[(XI.3.14)]{F}) that for all sufficiently large $k$
\beqn\label{l_dev}
P[ \tau^x_{R} > k] \leq e^{-\la k/R^\a}.
\eeqn 
Hence  
\beqn\label{N_R}
 \e(x,y;R)  \leq  Ce^{-\la N/R^\a} /n^{1/\a} =  O(n^{-\la m}/n^{1/\a}),
\eeqn
so that  $\e(x,y;R) $ is negligible if $m> 2/\la$ and our task reduces to 
the evaluation of  the double sum in (\ref{decomp}) with an appropriate choice of $R=R_n$. 

It is  easily  seen  that  at least within $|x_n|\vee |y_n| < M$
\begin{eqnarray}
0\leq p_{\{0\}}^{n}(x,y) - p_{A}^n(x,y) 
&\leq& C\sup_{k\leq n/2}[\,  p^{n-k}_{\{0\}}(0,y) +f_A^x(n-k)] \nonumber\\
&\leq& C'[f^{-y}(n) + f^x(n)]
\label{0-A}
\end{eqnarray}
(cf. the proof of \cite[Lemma 5.1]{U1dm_f} for the first inequality and Theorem \ref{thm5} for the second).   If $C^+ = \lim_{x\to\infty}a(x) = \infty$, it therefore follows that  within $|x|\vee |y| < M n^{1/\a}$
\beqn\label{A/0}  p_A^{n}(x,y) \sim p_{\{0\}}^n(x,y) \quad \mbox{as}\quad |x|\wedge |y| \wedge n \to \infty
\eeqn
and hence  
 both (\ref{A1-0}) and (\ref{A1-1}) hold for $\ga=2-\a$  if $ |x|\wedge y\to\infty$  in view of  Theorems  \ref{thm1} through \ref{thm5}.  
  In the sequel we suppose $\ga=2-\a$  (entailing $u_A(-y) \sim y^{\a-1}/\Ga(\a)$; the other case being similarly dealt with)   
and  verify that the restriction $|x|\wedge y \to\infty$ can be removed in the above. In case  $C^+ <\infty$  the situation is not much different and rather simpler.  At the end of the section  we shall advance certain remarks about the extension of Theorem \ref{thm6}

 In the sequel we shall tacitly suppose $|x|\vee |y| =O(n^{1/\a})$. 
\v2 

{\sc Let $C^+= \infty$ so that (\ref{A/0}) holds.}   First of all we  observe that 
in  the  case  $x_n\to 0$ of (\ref{eq_thm4}) the second term on its RHS  is negligible relative to the first so that $p_{\{0\}}^{n}(x,y) \sim a^\dagger(x)f^{-y}(n)$,    if $|x|= o(n^{2/\a-1})$ (since $f^{0}(n) \asymp 1/n^{2-1/\a}$). This together with
 (\ref{A/0})  and (\ref{eq_thm1}) implies  that as $  |x| \wedge y \to \infty$ 
  \begin{eqnarray}\label{A10}
p_A^{n}(x,y) \sim \left\{ \begin{array} {ll} u_A(x) f^{0}(n) u_{-A}(-y) \quad & (y_n \to0)\\
u_A(x) c_\circ\mathfrak{f}^{-y_n}(c_\circ)/n \quad &(1/M\leq  y_n< M)
\end{array}\right.  \nonumber \\
  \quad \mbox{under}\;\; x =o(n^{2/\a-1}), 
\end{eqnarray}
and, in view of duality,  for the proof of (\ref{A1-1}) it suffices  to  show that (\ref{A10})  remains true for each $x$ fixed. 
To this end we prove  two lemmas, Lemmas \ref{lem6.1} and \ref{lem6.2}; the proof of  (\ref{A10}) will be given after that of Lemma \ref{lem6.2}.

\begin{lem}\label{lem6.1}\,   {\rm (i)}  For $|x|<R$
\beqn\label{L6.1}
u_A(x)=  E[ u_A( S^x_{\tau_R}); \tau^x_R< \sigma_A].
\eeqn

{\rm (ii)} \, Let $\ga=2-\a$ and $b>1$. Then uniformly  for   $|x|<R$,  as $R\to\infty$
$$ E[ u_A( S^x_{\tau_R});\tau^x_R =\sigma^x_{(- bR, -R]}< \sigma^x_A] = u_A(x) +\{a^\dagger(x) +a(-x) + x_+/R^{2-\a}\}\times o(1).$$
\end{lem}
\n\pf\,  The process  $m_n:=
u_A(S^x_{n}) \1(  n < \sigma^x_A)$ 
is a martingale for each $x\in \Z$.  Noting  that  $m_n=u_A(S^x_{n\wedge \sigma_A}) \1( S_{n\wedge \sigma_A}^x\notin A)$ for $n\geq 1$ and
using  the optional stopping theorem we see 
 $$u_A(x)=  E[ u_A( S^x_n);  n <\sigma^x_A\wedge \tau^x_R] + E[ u_A( S^x_{\tau_R});\tau^x_R  \leq n \wedge \sigma^x_A]. $$
The first expectation approaches  zero as $n\to\infty$ since $u_A$ is bounded on $(-R,R)$
so that the monotone convergence shows  (\ref{L6.1}).
  Turning to the proof of (ii) let $\ga=2-\a$ and  $B(R)= (-\infty,-R]\cup A\cup [R,\infty)$.  We may suppose $0\in A$ for simplicity. Putting $\bar a^\dagger(x):= [a^\dagger(x) +a(-x)]/2$ so that $g_{B(R)}(x,\cdot)  \leq g_{\{0\}}(x,x) =2\bar a^\dagger(x)$  
we see 
\begin{eqnarray}\label{CR}
P[S^x_{\tau_R}  \leq  z, \tau^x_R< \sigma^x_A]&=& \sum_{w\notin B(R)} g_{B(R)}(x, w) F(z-w) \nonumber\\
& \leq& 4R\bar a^\dagger(x) F(z+R)\qquad\qquad  (z\leq -R)
\end{eqnarray}
 and making summation by parts we deduce that for any $b>1$,
 \beq 
 &&E[ u_A( S^x_{\tau_R});  \tau^x_R =\sigma^x_{(-\infty, -bR]}< \sigma^x_A\,] \\
 &&\leq C \sum_{z\leq - bR} |z|^{\a-1} P[S^x_{\tau_R} =z, \tau^x_R< \sigma^x_A] \\
 && \leq C'\bar a^\dagger(x)R \Big((bR)^{\a-1} F(-bR +R) +  \sum_{z\leq - bR} |z|^{\a-2} F(z+R) \Big),
 \eeq
 of which  the last member divided by $\bar a^\dagger(x)$ tends to zero since $F(z)= o(|z|^{-\a})$ as $z\to -\infty$. By virtually  the same way we derive a bound analogous to (\ref{CR}) and make summation by parts again to obtain  
$$E[ u_A( S^x_{\tau_R}); \tau^x_R =\sigma^x_{[bR,\infty)]}< \sigma^x_A] \leq C \sum_{z\geq bR} z^{\a-1} P[S^x_{\tau_R} =z, \tau^x_R<\sigma_A^x]\times o(1) \leq \bar a^\dagger(x) \times o(1), $$
where we have  $o(1)$ since $u_A(z) =o(z^{\a-1})$ as $z\to\infty$. 
It holds that  for  $|x|<R$,
$$
P[ \tau^x_R = \sigma^x_{[R,\infty)} <\sigma^x_A] \leq P[ \sigma^x_{[R,\infty)} <\sigma^x_A] 
 \leq C\{a^\dagger(x)/R^{\a-1}+x_+/R\}
 $$
(see Lemma \ref{lem7.9})),  which together with $u_A(z) =o(z^{\a-1})$ shows
$$E[ u_A( S^x_{\tau_R}); \tau^x_R =\sigma^x_{[R, bR)]} <\sigma_A] = (a^\dagger(x)+x_+R^{\a-2})\times o(1) $$
uniformly for  $|x|<R$.  Now the assertion of Lemma \ref{lem6.1} is easy to verify.  \qed



\begin{lem}\label{lem6.2}\,Suppose $C^+= \infty$.   For each $x$, as $|y|\wedge n\to\infty$ under $|y|<Mn^{1/\a}$
\[
p_A^n(x,y) \sim u_A(x)f^{-y}(n).
\]
\end{lem}
\n\pf\,  In (\ref{decomp})  take $R=R_n\sim  n^{2/\a-1}/\lg n$. Then,  by virtue of  (\ref{0-A}) and Corollary \ref{cor2}, uniformly for $-2R< z <- R$, $k\leq N$ and $|y|\leq Mn^{1/\a}$ as $n\wedge |y|\to\infty$
\beqn\label{eqL6.30}
p_A^{n-k}(z,y) \sim p_0^{n-k}(z,y) \sim a(z)f^{-y}(n). 
\eeqn 
 We can replace $a(z)$  by $u_A(z)$ in the right-most member for obvious reason. Then by  the exponential  bound  (\ref{l_dev}) and (\ref{N_R}) 
\beq
p_A^{n}(x,y) &=&  E[\, p^{n-\tau^x_R}_A(S^x_{\tau_R},y); \tau^x_{R} <\sigma^x_A\wedge N]+o(1/n^{2-1/\alpha})\\
&\geq& E[u_A(S^x_{\sigma(- 2R,-R]}); \tau^x_{R} =\sigma^x_{(- 2R,-R] }<\sigma^x_A\wedge N] 
f^{-y}(n)\{1+o(1)\} \\
&=& u_A(x)f^{-y}(n) \{1+o(1)\},
\eeq
where the last equality follows from Lemma \ref{lem6.1}.
Comparing this applied with $A=\{0\}$ and   the formula of Corollary \ref{cor2}  we see that  the inequality sign  above must be replaced by the equality sign for $A=\{0\}$, and hence   
\beqn\label{eqL6.3}
E[\, p^{n-\tau^x_R}_{\{0\}}(S^x_{\tau_R},y);  \tau^x_{R} \neq \sigma^x_{(- 2R,-R] }, \tau^x_R<\sigma^x_{\{0\}}\wedge N]  =o(f^{-y}(n)),
\eeqn
 which shows that the same replacement of the inequality sign is valid for $A$ itself, concluding   the asserted relation. \qed

\v2

Note that the case  $|x|\to \infty$ and $|y|<M$ is included in  Lemma \ref{lem6.2}  by duality relation so that  for each $y$ with $u_{-A}(-y)\neq 0$,
 \beqn\label{eqL6.31}
 p_A^{n}(x,y)=p_{-A}^{n}(-y, -x)  \sim  f^{x}(n)u_{-A}(-y). 
\eeqn
 
  It remains to deal with the case $|x|\vee |y| <M$, but  now having   (\ref{eqL6.31}) available  we may replace  the RHS  of (\ref{eqL6.30})  by $f^z(n)u_A(-y) \sim u_A(z)f^0(n)u_A(-y)$  for $|y|<M$
and repeating
the same argument made after (\ref{eqL6.30})   leads to
  the required relation.  Thus we have shown the formula  for $p_A^{n}(x,y)$ of Theorem \ref{thm7} in case $C^+=\infty$.  
    That for $f^x_A(n)$ follows from it  in view of the expression of  $P[\sigma^x_A=n, S^x_n=y]$ given in  (\ref{HD_st})
with the help of $\sum_{z\notin A} u_{-A}(-z)p(y-z) = u_{-A}(-y)$ and $\sum_{y\in A} u_{-A}(-y)=1$.
  \v2

{\sc Case $C_+<\infty$.}  
  Recalling  (\ref{0-A}), namely   $0\leq p^n_{\{0\}}(x,y) - p^n_A(x,y) \leq C\{f^{-y}(n) + f^x(n)\}$ we infer from Theorem \ref{thm5} (see also  Corollary \ref{cor2}) that \beqn\label{A/01}
p^n_A(x,y) \sim p^n_{\{0\}}(x,y)  \quad  \mbox{as} \quad  (-x) \wedge y \to\infty,
 \eeqn
 which allows us to follow the same arguments made above
to  verify  (\ref{A1-1}) except for the case $x>-M,  |y| <M$.  The completion of the proof is 
performed as follows. By duality we may consider the case  $|x|<M, y<M$. By virtue of (\ref{A/01}) (applied with $S^x_{\tau_R}\leq -R$ in place of $x$)  the argument deriving  (\ref{eqL6.3}) is valid for such $x,y$ and hence using Lemma \ref{lem6.1}
 we deduce that for $ |x| < M$, $y<M$ as  $R\wedge (-y)\to\infty$,
\beq
p^n_A(x,y) 
&=&  E[ p^{n-\tau}_{A}(S^x_{\tau_R},y); \tau^x< N\wedge \sigma_A] \{1+o(1)\}\\
&=& E[ a^\dagger(S^x_{\tau_R}); \tau_R=\sigma_{(-\infty,-R]} < \sigma_A] f^{-y}(n) + o(f^{-y}(n)), \\
&=& u_A(x)f^{-y}(n)\{1+o(1)\}, 
\eeq
hence  as before we have (\ref{eqL6.31}) as $x\to\infty$ for $|y|<M$ with $u_{-A}(-y)\neq 0$,  and can repeat the same argument to conclude  the formula asserted in Theorem \ref{thm7}.  $f^x_A(n)$ is  dealt with  as in the case $C^+=\infty$.  The proof of Theorem \ref{thm7} is complete.


\v2
We conclude this section with a  short  remark about the extension of Theorem \ref{thm6}. On letting $x\to \infty$ in (\ref{uA/a})  it follows that 
$$C_A^+ := \lim_{x\to +\infty} u_A(x) = C_+ - \sum_{z\in A}H^{+\infty}_A(z) a(z -w_0) \leq \infty$$
independently of the choice of $w_0\in A$, where  $H^{+\infty}_A(z) :=\lim_{x\to\infty} H^x_A(z)$  (cf. \cite[Theorem 30.1]{S} for  existence of the limit).
 We may suppose
 $A\subset (-\infty,0]$. 
Then for $x\geq 1$ and $y < \min A$, we have 
$g_A(x,y) = \sum_{z\notin A, z\leq 0}H_{(-\infty,0]}^x(z)g_A(z,y)$
and noting $g_A(z,y) \leq g_{\{0\}}(z,z)$  let first $y\to-\infty$ and then  $x\to +\infty$  to see that
\beqn\label{uA/H}
u_A(x)=  \sum_{z\notin A, z\leq 0}H_{(-\infty,0]}^x(z)u_A(z)\quad{and} \quad C_A^+=  \sum_{z\notin A, z\leq 0}H_{(-\infty,0]}^{+\infty}(z)u_A(z).
\eeqn
With these identities we can follow the proof of Theorem \ref{thm6} word for word except for trivial modifications to obtain the corresponding formula for $p^n_A(x,y)$.

\section{Some properties of $\mathfrak{f}^\xi$ and $\mathfrak{p}^{\{0\}}_t$ }
We have stated the asymptotic form of $\mathfrak{f}^1(t)$ in Corollary \ref{cor1}. In the following lemma we obtain it for $\gamma \neq 2-\alpha$ by  direct computation concerning the limit stable process.

\begin{lem}\label{lem7.1} \, As $t\to \infty$ 
$$\mathfrak{f}^1(t) \sim  \bigg[\frac{\sin (\pi/\a)}{\pi \mathfrak{p}_1(0)}\int_0^\infty u^{1-\a} \mathfrak{p}_1'(-u)du\bigg]\frac{1}{t^{2-1/\a}}.$$
\end{lem}
\n
\pf\, According to   \cite[Lemma 8.13]{Bt}  
$$\int_0^t \mathfrak{f}^1(s) ds=\frac{\sin (\pi/\a)}{\pi \mathfrak{p}_1(0)}\int_0^t(t-s)^{-1+1/\a}\mathfrak{p}_s(-1)ds.$$
Substitution from $\mathfrak{p}_s(-1)= s^{-1/\a}\mathfrak{p}_1(-s^{-1/\a})$ and the change  of variable $u=s/t$ transform  the integral on the RHS into
$$\int_0^1(1-u)^{-1+1/\a} u^{-1/\a} \mathfrak{p}_1(-(tu)^{-1/\a})du.$$
On  noting that $\int_0^1u^{-1/\a-1}|\mathfrak{p}_1'(-u^{-1/\a})|du = \a\int_1^\infty |\mathfrak{p}_1'(-s)|ds<\infty$ differentiation leads to
\beqn\label{eq3.16}
\mathfrak{f}^1(t)= \frac{\sin (\pi/\a)}{\pi \mathfrak{p}_1(0)}\cdot\frac{1}{\a t^{1+1/\a}}\int_0^1(1-u)^{-1+1/\a}u^{-2/\a}\mathfrak{p}_1'(-(tu)^{-1/\a})du,
\eeqn
 After the change of variable $u=1/ts^\alpha$ this becomes
$$\mathfrak{f}^1(t)= \frac{\sin (\pi/\a)}{\pi \mathfrak{p}_1(0)}\cdot \frac{1}{ t^{2-1/\a}}\int_{1}^\infty \Big(1-\frac{1}{ts^\alpha}\Big)^{-1+1/\a} s^{1-\a}\mathfrak{p}_1'(-s)ds,$$
which shows  the relation of the lemma,   the integral above   being  asymptoically equivalent to  $\int_{1}^\infty  s^{1- \a}\mathfrak{p}_1'(-s)ds$ as $t\to\infty$. \qed


\begin{lem}\label{lem7.2} If  $\fa(t)$ is a continuous function on $t\geq 0$, then for  $T>0$
$$\lim_{y\to \pm 0}\int_0^T \frac{\mathfrak{p}_t (y) -\mathfrak{p}_t(0)}{|y|^{\a-1}}\fa(t)dt  =b^\pm_{\a,\ga} \fa(0) \quad  \mbox{with}\quad b^\pm_{\a,\ga}=\a\int_0^\infty \frac{\mathfrak{p}_1 (\pm u) -\mathfrak{p}_1 (0)}{u^{\a}}du. $$
\end{lem}
\n
\pf\,  Let $w_y(t) = [\mathfrak{p}_t (y) -\mathfrak{p}_t(0)]/|y|^{\a-1}$. For any $\e>0$,
$$\int_0^\e w_y(t)dt =\int_0^\e [\mathfrak{p}_1(y/t^{1/\a}) -\mathfrak{p}_1(0)] \frac{|y|dt}{|y|^{\a} t^{1/\a}} = \a\int_{|y| /\e^{1/\a}}^\infty \frac{\mathfrak{p}_1(\pm u)-\mathfrak{p}_1(0)}{|u|^\a}du,$$
where $\pm$ accords to the sign of $y$. The last member converges to $b^\pm_{\a,\ga}$ 
and   $w_y(t) = O(|y|^{2-\a}) \to 0$  ($y\to0$) uniformly for $t>\e$ (since  $\mathfrak{p}_t'$ is  bounded), and hence the result follows.  \qed

\v2
\begin{lem}\label{lem7.3} Let  $b_{\alpha,\gamma}^\pm$ be given as in Lemma \ref{lem7.2}.  Then for $x>0$
$$\lim_{y\to \pm 0}\mathfrak{p}_t^{\{0\}}(x,y)/|y|^{\a-1} = b^\pm_{\a,\ga} \mathfrak{f}^x(t);
$$
and  
$$b^\pm_{\a,\ga}  = - \pi^{-1}\Ga(1-\a)\sin [\pi(\a\mp \ga)/2],$$
in particular if  $\ga =2-\a$,  $b^-_{\a,\ga} =0$ (the trivial case) and  $b^+_{\a,\ga} = 1/\Ga(\a)$.
\end{lem}
\n
\pf\, Although the result follows from Theorems \ref{thm4} and \ref{thm5}, we use them only for the identification of $b^\pm_{\a,\ga}$  in this proof that  is based  on the identity 
$$\mathfrak{p}_t^{\{0\}}(x,y) =\mathfrak{p}_t(y-x) -\int_0^t \mathfrak{f}^x(t-s)\mathfrak{p}_s(y)ds.$$
 On subtracting from this equality that for $y=0$ when  the LHS  vanishes, and then dividing by $|y|^{\a-1}$
$$\frac{\mathfrak{p}_t^{\{0\}}(x,y)}{|y|^{\a-1}} =\frac{\mathfrak{p}_t(y-x) - \mathfrak{p}_t(-x)}{|y|^{\a-1}} - \int_0^t \frac{\mathfrak{p}_s(y) -\mathfrak{p}_s(0)}{|y|^{\a-1}} \mathfrak{f}^x(t-s)ds.$$
 As $y\to 0$, the first term on the RHS tends to zero and Lemma \ref{lem7.2} applied to the the second term yields  the equality of the lemma.   
By applying Theorems \ref{thm4} and \ref{thm5} with $c_\circ =1$ it follows that
$$b^\pm_{\a,\ga}  =\frac1{ \mathfrak{f}^x(1)}\lim_{y_n \to \pm 0} \frac{\mathfrak{p}_1^{\{0\}}(x,y_n)}{|y|^{\a-1}/n^{1-1/\alpha}} = \lim_{y\to\pm \infty}\frac{ a(-y)}{|y|^{\a-1}}$$
 (see Remark \ref{rem2}(e)), of which the last limit  is evaluated in Lemma \ref{lem3.1}(i) as asserted.
 \qed
 
\v2
Let $Q_t(y)$ denote the distribution function of a stable meander, which may be expressed as
\beqn\label{Bt}
Q_t(y) = \lim_{\e\downarrow 0} P[ Y_t \leq y\,|\, \sigma_{(-\infty,-\e]} > t ]
\eeqn
(cf. \cite[Theorem 18]{Bt}) and  satisfies the scaling relation $Q_t(y) =Q_1(y/t^{1/\a})$.
\begin{lem}\label{lem7.4}  If $\ga=2-\a$, then  for $y>0$
\beqn\label{eq_lem6.4}
K_t(y) :=\lim_{x\downarrow 0}\mathfrak{p}_t^{\{0\}}(x,y)/x = \a \mathfrak{p}_t(0)Q'_t(y)
\eeqn
and
\beqn\label{eq_lem6.41}
\lim_{y\downarrow 0} \frac{\mathfrak{p}_t^{\{0\}}(x,y)}{y^{\a-1}} =\frac{\mathfrak{f}^x(t)}{\Ga(\a)} =  \frac{\hat Q'_t(x)}{\Ga(\a)\Ga(1/\a) t^{1-1/\a}}.
\eeqn
 \end{lem}
\n
\pf\, First we show (\ref{eq_lem6.4}).  In the proof of  Lemma \ref{lem4.4} (that if adapted to the stable process   is much simplified)  it is in effect shown that  the convergence in (\ref{eq_lem6.4}) is locally uniform in $y>0$, and our task is to identify  the limit,  for which it suffices to show
\beqn\label{6.3}
\lim_{x\downarrow 0}\frac1{x} \int_\delta^y \mathfrak{p}_t^{\{0\}}(x,z)dz =  \a \mathfrak{p}_t(0)[Q_t(y)-Q_t(\delta)]
\eeqn
for any $0<\delta<y$.
For  $\ga=2-\a$,   $\sigma_{(-\infty,-\e]}^Y$ agrees with $\sigma_{\{-\e\}}^Y$ a.s. Hence for $x>0$,  the integral in (\ref{6.3}) which equals 
$  P[\delta-x <Y_t\leq y-x, \sigma_{\{-x\}}^Y>t ]$ (since $\sigma^{x+Y}_{\{0\}} = \sigma^Y_{\{-x\}}$) is expressed as
$$P[\delta-x< Y_t\leq y-x\,|\, \sigma_{(-\infty,-x]}^Y > t] P[ \sigma_{(-\infty,-x]}^Y > t]. $$
The first factor  converges to $Q_t(y)-Q_t(\delta)$ as $x\downarrow 0$. For the second one, recalling  $\mathfrak{f}^x(s)= xs^{-1}\mathfrak{p}_s(x)= xs^{-1-1/\a}\mathfrak{p}_1(xs^{-1/\a})$ and making a change of variable we have
$$P[ \sigma^Y_{(-\infty,-x]} > t] =  \int_t^\infty \mathfrak{f}^x(s) ds = \a\int_0^{x/t^{1/\a}} \mathfrak{p}_1(u)du.$$
Thus dividing by $x$ and passing to the limit conclude the required formula (\ref{6.3}) since $\mathfrak{p}_1(0)t^{-1/\a}=\mathfrak{p}_t(0)$. 

 As for (\ref{eq_lem6.41}) we make use of the duality relation and write (\ref{Bt})
as
$$\hat Q_t(x) = \lim_{\e\downarrow 0}
\frac{\int_0^{x+\e}\mathfrak{p}_t^{[0, \infty)}(-\e,-\xi)d\xi}{P[\sigma_{[\e,\infty)}>t]}
=\lim_{\e\downarrow 0}
\frac{\int_0^{x+\e}\mathfrak{p}_t^{(-\infty,0]}(\xi,\e)d\xi}{P[\sigma_{[\e,\infty)}>t]}.
$$
The first equality of (\ref{eq_lem6.41}) follows from the preceding lemma  and is written as  
$\mathfrak{p}_t^{(-\infty,0]}(\xi,\e) = \mathfrak{p}_t^{\{0\}}(\xi,\e)\sim \mathfrak{f}^\xi(t)\e^{\a-1}/\Gamma(\alpha)$ ($\xi>0$).  By $\gamma= 2-\alpha$ we have $P[Y_t>0] =1-1/\alpha$   (cf. (\ref{rho})) which entails 
 $P[\sigma_{[\e,\infty)}>t]= P[\sigma_{[1,\infty)}>t/\e^\a] \sim C_* (t/\e^\a)^{-1+1/\a}$ \cite[Proposition VIII.2]{Bt} and  accordingly obtain 
$$\hat Q_t(x) = \frac{t^{1-1/\a}}{C_*\Ga(\a)} \int_0^x \mathfrak{f}^\xi(t)d\xi.$$
We derive $C_* = 1/\Ga(\a)\Ga(1/\a)$   from  $\hat Q_t(+\infty)=1$ with the help of the next lemma (cf. Remark \ref{R9}).  
Finally  differentiation concludes  the second equality of  (\ref{eq_lem6.41}). \qed

\v2
\begin{lem}\label{lem7.5}
$$\int_{-\infty}^\infty \mathfrak{p}_1(x)|x|dx= \frac{2 t^{1/\a} }{\pi }\Ga(1-1/\a)\sin[\tst12 \pi (\a-\ga)/\a],$$
in particular  if $\ga=2-\a$, $\int_0^\infty \mathfrak{f}^x(t)dx = t^{-1}\int_0^\infty \mathfrak{p}_1(-x)xdx = t^{-1+1/\a}/\Ga(1/\a) $.
\end{lem}
\n
\pf\, Put $\chi_\la(x) =|x|e^{-\la |x|}$ ($\la>0, -\infty <x<\infty$). 
 By Parseval equality
$$\int_{-\infty}^\infty \mathfrak{p}_t(x)|x|dx = \lim_{\la\downarrow 0}\int_{-\infty}^\infty \mathfrak{p}_t(x)\chi_\la(x)dx
= \frac1{\pi}\lim_{\la\downarrow 0} \int_{-\infty}^\infty e^{- t \psi(\th)} C_\la(\th) d\th, $$
where $C_\la(\th) =\int_0^\infty \chi_\la(x) \cos \th x \, dx,$ or explicitly 
$$C_\la(\th) =  \frac{\la^2-\th^2}{(\la^2+\th^2)^2}.$$
 Observing $\int_0^\infty C_\la(\th)d\th =0$, we infer that as $\la \downarrow 0$
$$\int_{0}^\infty e^{-t \psi( \th)} C_\la(\th) d\th= \int_{0}^\infty [e^{- t \psi(\th)}-1] C_\la(\th) d\th \; \longrightarrow  \int_{0}^\infty \frac{1- \exp\{-t e^{i\ga\pi/2}\th^\a\}}{ \th^2} d\th.$$
The last integral equals $(t e^{i\ga\pi/2})^{1/\a}\Ga(1-1/\a)$ \cite[p.313 (18)]{E}, and we find the first formula of the lemma obtained. If $\ga=2-\a$,  then $\Ga(1-1/\a)\sin[\tst12 \pi (\a-\ga)/\a]= \pi/\Ga(1/\a)$, which together with $\mathfrak{f}^x(t)=xt^{-1}\mathfrak{p}_t(-x)$ and $\int_{-\infty}^\infty \mathfrak{p}_t(x)xdx=0$ shows the second 
formula.
\qed

\begin{rem} \label{R9}\, We have used Lemma \ref{lem7.5} for identification of the constant factor in (\ref{eq_lem6.41}). Alternatively we could have applied the exact formula for $P[\sup_{s\leq t} Y_s\in d\xi]/d\xi$ obtained in  \cite{BDP} (cf. also  \cite{DS}).
\end{rem}

\section{Auxiliaries}

Here we give
miscellaneous consequences  of the  assumptions 1) and 2)  stated in Section 1 that are derived from  the general theory. 

 \subsection{Condition   (\ref{f_hyp}) in terms of the tails of $F$} 
The assumption (\ref{f_hyp}) on the characteristic function  $\phi(\th)$ is equivalent to the condition
\beqn\label{A7.1}
P[X>x] \sim q^+ Bx^{-\a} \quad\mbox{and}\quad  P[X< -x] \sim  q^-Bx^{-\a}
\eeqn
as $x\to\infty$ with some positive constant $B$ and two non-negative constants $q^+$ and $ q^-$ such that
 $q^++q^-=1$. The L\'evy measure  $M\{dx\}$ is then given by
 $$M\{(-y,x]\} =\frac{\a B}{2-\a}(q^- y^{2-\a} + q^+ x^{2-\a}) \quad (x\geq 0, y\geq 0)$$
 and the characteristic exponent of the limiting stable process by
 $$c_\circ\psi(\th) =  |\th|^{\a}B\Ga(1-\a)\{\cos \tst12 \a\pi - i(\sgn\, \th) (q^+ - q^-)\sin \tst12 \a\pi\}$$
(cf. \cite[(XVII.3.18)]{F}). From this we read off
 $$c_\circ = B\Ga(1-\a)[(\cos \tst12\a \pi)/(\cos \tst12 \ga\pi)] \quad \mbox{and}\quad    \tan\tst12 \ga\pi = (q^+ - q^-)(- \tan \tst12 \a\pi )$$
  (which reduce to $c_\circ =-B\Ga(1-\a)$ and $q^+ =1$, respectively,  if $\ga =2-\a$) and hence
 $$\psi(\th) =  |\th|^\a(\cos \tst12 \ga\pi)\{1+ i (\sgn\, \th)(q^+- q^-) (-\tan \tst12 \a\pi)\};$$
According to Zolotarev \cite{Zol}  (cf. \cite[Section 8.9.2]{BGT},  \cite[Section VIII.1]{Bt}) Spitzer's  constant  $\rho:= \lim_{n\to \infty} n^{-1}\sum_{k=1}^n P[S_k>0]$ is given by
 \beqn\label{rho}
   \rho= \frac12(1-\ga/\a).
   \eeqn
 
 From (\ref{A7.1}) it follows  that
 \beqn\label{F/psi}
 \phi'(\th) \sim -\psi'(\th) =\mp \a c_\circ e^{\pm i\pi\ga/2}|\th|^{\a-1} \quad \mbox{as}\;\; \th \to  \pm0.
 \eeqn
 Indeed, on writing $ \phi'(\th) = i\int_{-\infty}^\infty (e^{i\th t}-1)tdF(t)$ the integration by parts  yields
 \beqn
 \phi'(\th) 
 =  i\int_{-\infty}^\infty \{ e^{i\th t}-1 + i \th t e^{i\th t}\}[-F(t)\1(t<0) + (1-F(t))\1(t>0)]dt,
 \eeqn
and scaling by the factor $1/|\th|$  we find that $ \phi'(\th)  \sim \pm \zeta |\th|^{\alpha-1}$, where
\[
 \zeta=  iB \int_{-\infty}^\infty \{1- e^{\pm iu} \mp i u e^{\pm iu} \}\frac{q^-\1(u<0) - q^+\1(u>0)}{|u|^\a}du.
 \]
 Since $\zeta$ depends on the  regularity of  tails of $F$ only and $-\psi'(\th)$ is given by the above integral with $d F$ replaced by the Levy measure associated with $\psi$,  $\pm \zeta|\th|^{\a-1}$ must be equal to $-\psi'(\th)$. 
 
 \v2
\begin{rem}\label{rem10}\,    Put  $\beta_\pm = \frac12 (\a\pm \ga)$.  If $q^+ < q^-$, then $\a-1\leq \beta_+ <   \beta_- \leq 1$, where the equality in each   extremity holds if and only if $q^+=0$.
 In order to consider the behaviour of $U_{{\rm ds}}$ and $V_{{\rm as}}$ we rewrite condition (\ref{A7.1})  as
\beqn\label{7.1}
P[ X>x] = B x^{-\a}(q^++r_+(x)) \quad\mbox{and} \quad P[ X< -x] = Bx^{-\a}(q^- +r_-(x)),
\eeqn
where $r_\pm(x)\to 0$ as $x\to 0$. 
If $\int_1^\infty (|r_+(x)|+|r_-(x)| )x^{-1}dx <\infty$, then, as $x\to\infty$
 \beqn\label{7.2}
 U_{{\rm ds}}(x)\sim C'\, x^{\beta_+} \quad \mbox{and}\quad V_{{\rm as}}(x)\sim C''\, x^{\beta_-}
 \eeqn
with some positive constants $C'$ and $C''$ such that  $1/C'C''= c_\circ\Ga(1+\beta_+)\Ga(1+ \beta_- )$.
The proof  is carried out by computation based on  Spitzer's expressions of the generating functions of $v^-$ and $v$ given in \cite[P18.7]{S},   the computation  being somewhat involved and omitted.
     In order that  $E Z<\infty$ (entailing (\ref{7.2}) with $\beta_+=\alpha-1$ and $\beta_-=1$), it is necessary and sufficient that   $q^+=0$   and $\int_1^\infty r_+(x)x^{-1}dx<\infty$ according to  Chow's criterion \cite{Cho}. If  $q^+=0$ with  $\int_1^\infty r_+(x)x^{-1}dx=\infty$,  then $U_{\rm{ds}}(x) =o(x)$ so that $L(x)\to 0$ in (\ref{LH})
     and  hence $  U_{{\rm ds}}(x)x^{-\a+1} \to \infty$ ($x\to\infty$). 
 \end{rem}

\subsection{An upper bound of $p^n(x)$ for $|x|>n^{1/\a}$}

Let  $\la(\th) = \sum_{x\in \Z} w_x e^{i\th x}$, the sum of  the trigonometric series with the coefficient such that  $M:=\sum |w_x| <\infty$ and put
$m(r) = \int_0^r dt \sum_{|x|>t} |w_x|$.
Then
$$
|\la(\th) -\la(\th')| \leq 2|\th-\th'|m(1/|\th-\th'|) \quad (\th\neq \th'),
$$
 (Erickson \cite[Lemma 5]{Ec}). Applied with $w_x = xp(x)$, $m(r) = \int_0^r E[|X| ; |X|>t]dt = O(r^{2-\a})$ this yields 
\beqn\label{Erc}
|\phi'(\th) -\phi'(\th')| \leq C |\th-\th'|^{\a-1} 
\eeqn
under (\ref{A7.1}). This same bound is satisfied  by $\psi'(\th)$ as is verified directly (or by a similar reasoning). 
The following  lemma is based  on this bound.
\begin{lem}\label{lem7.6} If  $p$ satisfies  (\ref{f_hyp}),  then  for some constant $C$
\[
p^n(x) \leq C n^{-1/\a}(1\wedge |x_n|^{-\a }).
\]
\end{lem}
\n
\pf\, For $|x_n|\leq1$, the bound follows from the local limit theorem.   Let $|x_n| > 1$.
 The  integration by parts and 
the change of variable $\th = t/n^{1/\a}$ transforms  the expression $p^n(x) =\frac1{2\pi} \int_{-\pi}^{\pi} [\phi(\th)]^ne^{-ix\th}d \th$ into
\beqn\label{6.5}
p^n(x)= \frac1{2\pi ix} \int_{-\pi n^{1/\a}}^{\pi n^{1/\a}} \frac{n}{n^{1/\a}}\phi'(tn^{1/\a})e^{(n-1) \log \phi(t/n^{1/\a})} e^{-ix_nt}dt.
\eeqn
Put
\[ Q_n(t)=- (n-1)\log \phi(t/n^{1/\a}) , \quad R_n(t) =  n^{1-1/\a} \phi'(t/n^{1/\a}).
\]
By  (\ref{f_hyp}),  (\ref{F/psi}) and (\ref{Erc}) it follows that
$$Q_n(t) = \psi(t)\{1+o(1)\}, \quad  R_n(t) = - \psi'(t)\{1+ o(1)\}$$
where $o(1)$ is bounded for $|t|<\pi n^{1/\alpha}$ and  $o(1)\to 0$  as  $t/n^{1/\a} \to 0$,
and 
\[R_n(t) -R_n(t')  =O(|t-t'|^{\a-1}) 
\]
uniformly for $n$.  By periodicity of $\phi$,   $p^n(x)= -\frac1{2\pi}\int_{-\pi}^\pi [\phi(\th+ \pi/x)]^ne^{-ix\th}d \th$,   we accordingly obtain that $p^n(x) = (I+J)/ 4\pi ix$, where
\[
 I= \int_{-\pi n^{1/\a}}^{\pi n^{1/\a}} [R_n(t)-R_n(t+\pi/x_n)] e^{-Q_n(t)}   e^{-ix_nt}dt
\]
and
\[
J=  \int_{-\pi n^{1/\a}}^{\pi n^{1/\a}} R_n(t+\pi/x_n)[e^{-Q_n(t)}- e^{-Q_n(t+\pi/x_n)}]   e^{-ix_nt}dt.
\]
Noting  $|\phi(\th)|<1$ for $0<|\th|\leq \pi$ and $\Re\, \psi(\th) /|\th|^\alpha =\cos \frac12 \gamma \pi >0$, we can choose a constant  $\lambda>0$ so that
$\Re \, Q_n(t)  > \lambda |t|^\alpha$
for $|t|< \pi n^{1/\alpha} + \pi$ for $n$ large enough. Hence if
 $f_n(t) = [R_n(t)-R_n(t+\pi/x_n)] e^{-Q_n(t)} |x_n|^{\alpha-1}$, then $f_n(t)$ is  dominated in absolute value  by $(e^{-|t|^{\alpha}/2})$, and we deduce that
$
I=|x_n|^{1-\alpha} \int_{-\pi n^{1/\alpha}}^{\pi n^{1/\alpha}} f_n(t)e^{-ix_n t}dt 
= O(|x_n|^{1-\alpha}). 
$
In a similar way we obtain $J =O(1/|x_n|)$. Since  $|x_n|^{1-\alpha}/|x| = n^{1/\alpha}|x_n|^{-\alpha}$,
  this concludes the proof. \qed

\subsection{Escape  probabilities from the origin} 
By \cite[Proposition 29.4]{S} 
$$g_{\{0\}}(x,y) = a^\dagger(x)+a(-y)-a(x-y),$$
which entails the subadditivity $a(x+y) \leq a(x)+ a(y)$ and 
 \beqn  \label{id_h_p}   
P[\sigma^x_{\{y\}}<\sigma^x_{\{0\}}]= \frac{g_{\{0\}}(x,y)}{g_{\{0\}}(y,y)} = \frac{a^\dagger(x)+a(-y)-a(x-y)}{a(y)+a(-y)}.
 \eeqn

 Put  $\om(r) = g_{(-\infty,0]}(r,r)$ ($r=1,2,\ldots$). Note that $\omega(r)=g_{[0,\infty)}(-r,-r)$ and $\om(r)  \leq g_{\{0\}}(r,r) = a(r)+a(-r).$
\begin{lem}\label{lem7.7} \,  Given a positive integer $R$, 
put $\tau^x_R =\inf\{n\geq 1: S^x_n\notin (-R,R)\}$. 
There exists a constant $C$ depending only on $F$ such that for any $R>1$, $z\geq 0$  and $|x|<R$, 
$$
 P[S^x_{\tau_R} =  R+z ] \leq \om (R-|x|) P[ X > z ].
$$
\end{lem}
\n
\pf\,  Since  $g_{[R, \infty)}(x,w) = g_{[0,\infty)}(x-R,w-R) \leq \om (R-|x|)$ and similarly for $g_{(-\infty, -R]}(x,w)$,      it follows that 
 $ g_{(-\infty, -R]\cup[R, \infty)}(x,w) \leq \om (R-|x|)$, and hence
  \[
P[S^x_{\tau_R} = R+ z] \leq \om(R-|x|)  \sum_{w: |w|<R}  p(R+ z-w)
 \leq \om (R-|x|)  P[X >  z],
 \] 
 showing the inequality of the lemma.   \qed
\v2
  
The same proof as above shows 
\beqn\label{(C)}
H^x_{[0,\infty)}(z) := P[S^x_{\sigma[0,\infty)}=z] \leq \om(-x)P[X \geq z] \quad (x\leq 0, z>0).
\eeqn  
\begin{lem}\label{lem7.8} \, If $\ga>\a-2$, 
\beqn \label{R/q}
\liminf_{R\to\infty} \inf_{x\in \Z} P[\sigma^x_{\{R\}}<\sigma^x_{\{0\}}\,|\, \sigma^x_{[R,\infty)}<\sigma^x_{\{0\}}] =: q>0;
\eeqn
with $q=1$ for $\ga= 2- \alpha$.  
 \end{lem}
 \v2\n
 \pf\,   In view of     (\ref{id_h_p})  and the decomposition
 $$ P[\sigma^x_{\{R\}}<\sigma^x_{\{0\}}\,|\, \sigma^x_{[R,\infty)}<\sigma^x_{\{0\}}] =\sum_{z\geq R} P[S^x_{\sigma [R,\infty)} =z\,|\, \sigma^x_{[R,\infty)}<\sigma^x_{\{0\}}] P[\bar\sigma^z_{\{R\}}<\sigma^z_{\{0\}}],$$
where $ \bar\sigma^z_{\{R\}}$ is defined to be  zero if $z=N$ and agree with  $\sigma^z_{\{R\}}$ otherwise,  for the first half of the lemma    it suffices  to show that
$$\lim_{R\to\infty} \inf_{z\geq R} \frac{a(z)+a(-R)-a(z-R)}{a(R)+a(-R)} = \frac{\k_{\a,\ga,-}}
{\k_{\a,\ga,-}+ \k_{\a,\ga,+}},$$
the last ratio being positive if $\ga>\a-2$ and equals unity if $\ga=2-\a$. If $|\ga|<2-\a$, by   Lemma \ref{lem3.1}(ii)  $a(z) -a(z-R) >0$ for $R$ large enough and the equality above follows immediately  from Lemma \ref{lem3.1}(i). The case $\ga =2-\a$   also follows from  Lemma \ref{lem3.1}(i) and (ii),  the latter    showing  $\sup_{z\geq R}|a(z)- a(z-R) |= o(R^{\a-1}).$ \qed

\v2

\begin{lem}\label{lem7.9} \,For any $\ga$ there exists  a constant $C$ such that for $R>1$,
\[
P[ \sigma^x_{[R,\infty)} <\sigma^x_{\{0\}}] \leq  C[a^\dagger(x)R^{-\a+1} + x_+/R] \quad (x\leq R).
\]
\end{lem}
\v2\n
\pf\,
For $\ga>\a-2$, on using   Lemma \ref{lem3.1}(ii)  the result is deduced   from the preceding lemma. The case  $\ga= \alpha -2$ follows by \cite[Lemma 5.5]{Upot} (use the fact that $m_n = a(S^x_{\sigma_{\{0\}\cup[R,\infty)}})$ is a martingale).  \qed



\end{document}